\documentclass[aos,preprint]{imsart}

\RequirePackage[OT1]{fontenc}
\RequirePackage{amsthm,amsmath,natbib}
\RequirePackage[colorlinks,citecolor=blue,urlcolor=blue]{hyperref}

\usepackage{color}  
\definecolor{AliceBlue}{rgb}{0.94,0.97,1.00}
\definecolor{AntiqueWhite1}{rgb}{1.00,0.94,0.86}
\definecolor{AntiqueWhite2}{rgb}{0.93,0.87,0.80}
\definecolor{AntiqueWhite3}{rgb}{0.80,0.75,0.69}
\definecolor{AntiqueWhite4}{rgb}{0.55,0.51,0.47}
\definecolor{AntiqueWhite}{rgb}{0.98,0.92,0.84}
\definecolor{BlanchedAlmond}{rgb}{1.00,0.92,0.80}
\definecolor{BlueViolet}{rgb}{0.54,0.17,0.89}
\definecolor{CadetBlue1}{rgb}{0.60,0.96,1.00}
\definecolor{CadetBlue2}{rgb}{0.56,0.90,0.93}
\definecolor{CadetBlue3}{rgb}{0.48,0.77,0.80}
\definecolor{CadetBlue4}{rgb}{0.33,0.53,0.55}
\definecolor{CadetBlue}{rgb}{0.37,0.62,0.63}
\definecolor{CornflowerBlue}{rgb}{0.39,0.58,0.93}
\definecolor{DarkBlue}{rgb}{0.00,0.00,0.55}
\definecolor{DarkCyan}{rgb}{0.00,0.55,0.55}
\definecolor{DarkGoldenrod1}{rgb}{1.00,0.73,0.06}
\definecolor{DarkGoldenrod2}{rgb}{0.93,0.68,0.05}
\definecolor{DarkGoldenrod3}{rgb}{0.80,0.58,0.05}
\definecolor{DarkGoldenrod4}{rgb}{0.55,0.40,0.03}
\definecolor{DarkGoldenrod}{rgb}{0.72,0.53,0.04}
\definecolor{DarkGray}{rgb}{0.66,0.66,0.66}
\definecolor{DarkGreen}{rgb}{0.00,0.39,0.00}
\definecolor{DarkGrey}{rgb}{0.66,0.66,0.66}
\definecolor{DarkKhaki}{rgb}{0.74,0.72,0.42}
\definecolor{DarkMagenta}{rgb}{0.55,0.00,0.55}
\definecolor{DarkOliveGreen1}{rgb}{0.79,1.00,0.44}
\definecolor{DarkOliveGreen2}{rgb}{0.74,0.93,0.41}
\definecolor{DarkOliveGreen3}{rgb}{0.64,0.80,0.35}
\definecolor{DarkOliveGreen4}{rgb}{0.43,0.55,0.24}
\definecolor{DarkOliveGreen}{rgb}{0.33,0.42,0.18}
\definecolor{DarkOrange1}{rgb}{1.00,0.50,0.00}
\definecolor{DarkOrange2}{rgb}{0.93,0.46,0.00}
\definecolor{DarkOrange3}{rgb}{0.80,0.40,0.00}
\definecolor{DarkOrange4}{rgb}{0.55,0.27,0.00}
\definecolor{DarkOrange}{rgb}{1.00,0.55,0.00}
\definecolor{DarkOrchid1}{rgb}{0.75,0.24,1.00}
\definecolor{DarkOrchid2}{rgb}{0.70,0.23,0.93}
\definecolor{DarkOrchid3}{rgb}{0.60,0.20,0.80}
\definecolor{DarkOrchid4}{rgb}{0.41,0.13,0.55}
\definecolor{DarkOrchid}{rgb}{0.60,0.20,0.80}
\definecolor{DarkRed}{rgb}{0.55,0.00,0.00}
\definecolor{DarkSalmon}{rgb}{0.91,0.59,0.48}
\definecolor{DarkSeaGreen1}{rgb}{0.76,1.00,0.76}
\definecolor{DarkSeaGreen2}{rgb}{0.71,0.93,0.71}
\definecolor{DarkSeaGreen3}{rgb}{0.61,0.80,0.61}
\definecolor{DarkSeaGreen4}{rgb}{0.41,0.55,0.41}
\definecolor{DarkSeaGreen}{rgb}{0.56,0.74,0.56}
\definecolor{DarkSlateBlue}{rgb}{0.28,0.24,0.55}
\definecolor{DarkSlateGray1}{rgb}{0.59,1.00,1.00}
\definecolor{DarkSlateGray2}{rgb}{0.55,0.93,0.93}
\definecolor{DarkSlateGray3}{rgb}{0.47,0.80,0.80}
\definecolor{DarkSlateGray4}{rgb}{0.32,0.55,0.55}
\definecolor{DarkSlateGray}{rgb}{0.18,0.31,0.31}
\definecolor{DarkSlateGrey}{rgb}{0.18,0.31,0.31}
\definecolor{DarkTurquoise}{rgb}{0.00,0.81,0.82}
\definecolor{DarkViolet}{rgb}{0.58,0.00,0.83}
\definecolor{DeepPink1}{rgb}{1.00,0.08,0.58}
\definecolor{DeepPink2}{rgb}{0.93,0.07,0.54}
\definecolor{DeepPink3}{rgb}{0.80,0.06,0.46}
\definecolor{DeepPink4}{rgb}{0.55,0.04,0.31}
\definecolor{DeepPink}{rgb}{1.00,0.08,0.58}
\definecolor{DeepSkyBlue1}{rgb}{0.00,0.75,1.00}
\definecolor{DeepSkyBlue2}{rgb}{0.00,0.70,0.93}
\definecolor{DeepSkyBlue3}{rgb}{0.00,0.60,0.80}
\definecolor{DeepSkyBlue4}{rgb}{0.00,0.41,0.55}
\definecolor{DeepSkyBlue}{rgb}{0.00,0.75,1.00}
\definecolor{DimGray}{rgb}{0.41,0.41,0.41}
\definecolor{DimGrey}{rgb}{0.41,0.41,0.41}
\definecolor{DodgerBlue1}{rgb}{0.12,0.56,1.00}
\definecolor{DodgerBlue2}{rgb}{0.11,0.53,0.93}
\definecolor{DodgerBlue3}{rgb}{0.09,0.45,0.80}
\definecolor{DodgerBlue4}{rgb}{0.06,0.31,0.55}
\definecolor{DodgerBlue}{rgb}{0.12,0.56,1.00}
\definecolor{FloralWhite}{rgb}{1.00,0.98,0.94}
\definecolor{ForestGreen}{rgb}{0.13,0.55,0.13}
\definecolor{GhostWhite}{rgb}{0.97,0.97,1.00}
\definecolor{GreenYellow}{rgb}{0.68,1.00,0.18}
\definecolor{HotPink1}{rgb}{1.00,0.43,0.71}
\definecolor{HotPink2}{rgb}{0.93,0.42,0.65}
\definecolor{HotPink3}{rgb}{0.80,0.38,0.56}
\definecolor{HotPink4}{rgb}{0.55,0.23,0.38}
\definecolor{HotPink}{rgb}{1.00,0.41,0.71}
\definecolor{IndianRed1}{rgb}{1.00,0.42,0.42}
\definecolor{IndianRed2}{rgb}{0.93,0.39,0.39}
\definecolor{IndianRed3}{rgb}{0.80,0.33,0.33}
\definecolor{IndianRed4}{rgb}{0.55,0.23,0.23}
\definecolor{IndianRed}{rgb}{0.80,0.36,0.36}
\definecolor{LavenderBlush1}{rgb}{1.00,0.94,0.96}
\definecolor{LavenderBlush2}{rgb}{0.93,0.88,0.90}
\definecolor{LavenderBlush3}{rgb}{0.80,0.76,0.77}
\definecolor{LavenderBlush4}{rgb}{0.55,0.51,0.53}
\definecolor{LavenderBlush}{rgb}{1.00,0.94,0.96}
\definecolor{LawnGreen}{rgb}{0.49,0.99,0.00}
\definecolor{LemonChiffon1}{rgb}{1.00,0.98,0.80}
\definecolor{LemonChiffon2}{rgb}{0.93,0.91,0.75}
\definecolor{LemonChiffon3}{rgb}{0.80,0.79,0.65}
\definecolor{LemonChiffon4}{rgb}{0.55,0.54,0.44}
\definecolor{LemonChiffon}{rgb}{1.00,0.98,0.80}
\definecolor{LightBlue1}{rgb}{0.75,0.94,1.00}
\definecolor{LightBlue2}{rgb}{0.70,0.87,0.93}
\definecolor{LightBlue3}{rgb}{0.60,0.75,0.80}
\definecolor{LightBlue4}{rgb}{0.41,0.51,0.55}
\definecolor{LightBlue}{rgb}{0.68,0.85,0.90}
\definecolor{LightCoral}{rgb}{0.94,0.50,0.50}
\definecolor{LightCyan1}{rgb}{0.88,1.00,1.00}
\definecolor{LightCyan2}{rgb}{0.82,0.93,0.93}
\definecolor{LightCyan3}{rgb}{0.71,0.80,0.80}
\definecolor{LightCyan4}{rgb}{0.48,0.55,0.55}
\definecolor{LightCyan}{rgb}{0.88,1.00,1.00}
\definecolor{LightGoldenrod1}{rgb}{1.00,0.93,0.55}
\definecolor{LightGoldenrod2}{rgb}{0.93,0.86,0.51}
\definecolor{LightGoldenrod3}{rgb}{0.80,0.75,0.44}
\definecolor{LightGoldenrod4}{rgb}{0.55,0.51,0.30}
\definecolor{LightGoldenrodYellow}{rgb}{0.98,0.98,0.82}
\definecolor{LightGoldenrod}{rgb}{0.93,0.87,0.51}
\definecolor{LightGray}{rgb}{0.83,0.83,0.83}
\definecolor{LightGreen}{rgb}{0.56,0.93,0.56}
\definecolor{LightGrey}{rgb}{0.83,0.83,0.83}
\definecolor{LightPink1}{rgb}{1.00,0.68,0.73}
\definecolor{LightPink2}{rgb}{0.93,0.64,0.68}
\definecolor{LightPink3}{rgb}{0.80,0.55,0.58}
\definecolor{LightPink4}{rgb}{0.55,0.37,0.40}
\definecolor{LightPink}{rgb}{1.00,0.71,0.76}
\definecolor{LightSalmon1}{rgb}{1.00,0.63,0.48}
\definecolor{LightSalmon2}{rgb}{0.93,0.58,0.45}
\definecolor{LightSalmon3}{rgb}{0.80,0.51,0.38}
\definecolor{LightSalmon4}{rgb}{0.55,0.34,0.26}
\definecolor{LightSalmon}{rgb}{1.00,0.63,0.48}
\definecolor{LightSeaGreen}{rgb}{0.13,0.70,0.67}
\definecolor{LightSkyBlue1}{rgb}{0.69,0.89,1.00}
\definecolor{LightSkyBlue2}{rgb}{0.64,0.83,0.93}
\definecolor{LightSkyBlue3}{rgb}{0.55,0.71,0.80}
\definecolor{LightSkyBlue4}{rgb}{0.38,0.48,0.55}
\definecolor{LightSkyBlue}{rgb}{0.53,0.81,0.98}
\definecolor{LightSlateBlue}{rgb}{0.52,0.44,1.00}
\definecolor{LightSlateGray}{rgb}{0.47,0.53,0.60}
\definecolor{LightSlateGrey}{rgb}{0.47,0.53,0.60}
\definecolor{LightSteelBlue1}{rgb}{0.79,0.88,1.00}
\definecolor{LightSteelBlue2}{rgb}{0.74,0.82,0.93}
\definecolor{LightSteelBlue3}{rgb}{0.64,0.71,0.80}
\definecolor{LightSteelBlue4}{rgb}{0.43,0.48,0.55}
\definecolor{LightSteelBlue}{rgb}{0.69,0.77,0.87}
\definecolor{LightYellow1}{rgb}{1.00,1.00,0.88}
\definecolor{LightYellow2}{rgb}{0.93,0.93,0.82}
\definecolor{LightYellow3}{rgb}{0.80,0.80,0.71}
\definecolor{LightYellow4}{rgb}{0.55,0.55,0.48}
\definecolor{LightYellow}{rgb}{1.00,1.00,0.88}
\definecolor{LimeGreen}{rgb}{0.20,0.80,0.20}
\definecolor{MediumAquamarine}{rgb}{0.40,0.80,0.67}
\definecolor{MediumBlue}{rgb}{0.00,0.00,0.80}
\definecolor{MediumOrchid1}{rgb}{0.88,0.40,1.00}
\definecolor{MediumOrchid2}{rgb}{0.82,0.37,0.93}
\definecolor{MediumOrchid3}{rgb}{0.71,0.32,0.80}
\definecolor{MediumOrchid4}{rgb}{0.48,0.22,0.55}
\definecolor{MediumOrchid}{rgb}{0.73,0.33,0.83}
\definecolor{MediumPurple1}{rgb}{0.67,0.51,1.00}
\definecolor{MediumPurple2}{rgb}{0.62,0.47,0.93}
\definecolor{MediumPurple3}{rgb}{0.54,0.41,0.80}
\definecolor{MediumPurple4}{rgb}{0.36,0.28,0.55}
\definecolor{MediumPurple}{rgb}{0.58,0.44,0.86}
\definecolor{MediumSeaGreen}{rgb}{0.24,0.70,0.44}
\definecolor{MediumSlateBlue}{rgb}{0.48,0.41,0.93}
\definecolor{MediumSpringGreen}{rgb}{0.00,0.98,0.60}
\definecolor{MediumTurquoise}{rgb}{0.28,0.82,0.80}
\definecolor{MediumVioletRed}{rgb}{0.78,0.08,0.52}
\definecolor{MidnightBlue}{rgb}{0.10,0.10,0.44}
\definecolor{MintCream}{rgb}{0.96,1.00,0.98}
\definecolor{MistyRose1}{rgb}{1.00,0.89,0.88}
\definecolor{MistyRose2}{rgb}{0.93,0.84,0.82}
\definecolor{MistyRose3}{rgb}{0.80,0.72,0.71}
\definecolor{MistyRose4}{rgb}{0.55,0.49,0.48}
\definecolor{MistyRose}{rgb}{1.00,0.89,0.88}
\definecolor{NavajoWhite1}{rgb}{1.00,0.87,0.68}
\definecolor{NavajoWhite2}{rgb}{0.93,0.81,0.63}
\definecolor{NavajoWhite3}{rgb}{0.80,0.70,0.55}
\definecolor{NavajoWhite4}{rgb}{0.55,0.47,0.37}
\definecolor{NavajoWhite}{rgb}{1.00,0.87,0.68}
\definecolor{NavyBlue}{rgb}{0.00,0.00,0.50}
\definecolor{OldLace}{rgb}{0.99,0.96,0.90}
\definecolor{OliveDrab1}{rgb}{0.75,1.00,0.24}
\definecolor{OliveDrab2}{rgb}{0.70,0.93,0.23}
\definecolor{OliveDrab3}{rgb}{0.60,0.80,0.20}
\definecolor{OliveDrab4}{rgb}{0.41,0.55,0.13}
\definecolor{OliveDrab}{rgb}{0.42,0.56,0.14}
\definecolor{OrangeRed1}{rgb}{1.00,0.27,0.00}
\definecolor{OrangeRed2}{rgb}{0.93,0.25,0.00}
\definecolor{OrangeRed3}{rgb}{0.80,0.22,0.00}
\definecolor{OrangeRed4}{rgb}{0.55,0.15,0.00}
\definecolor{OrangeRed}{rgb}{1.00,0.27,0.00}
\definecolor{PaleGoldenrod}{rgb}{0.93,0.91,0.67}
\definecolor{PaleGreen1}{rgb}{0.60,1.00,0.60}
\definecolor{PaleGreen2}{rgb}{0.56,0.93,0.56}
\definecolor{PaleGreen3}{rgb}{0.49,0.80,0.49}
\definecolor{PaleGreen4}{rgb}{0.33,0.55,0.33}
\definecolor{PaleGreen}{rgb}{0.60,0.98,0.60}
\definecolor{PaleTurquoise1}{rgb}{0.73,1.00,1.00}
\definecolor{PaleTurquoise2}{rgb}{0.68,0.93,0.93}
\definecolor{PaleTurquoise3}{rgb}{0.59,0.80,0.80}
\definecolor{PaleTurquoise4}{rgb}{0.40,0.55,0.55}
\definecolor{PaleTurquoise}{rgb}{0.69,0.93,0.93}
\definecolor{PaleVioletRed1}{rgb}{1.00,0.51,0.67}
\definecolor{PaleVioletRed2}{rgb}{0.93,0.47,0.62}
\definecolor{PaleVioletRed3}{rgb}{0.80,0.41,0.54}
\definecolor{PaleVioletRed4}{rgb}{0.55,0.28,0.36}
\definecolor{PaleVioletRed}{rgb}{0.86,0.44,0.58}
\definecolor{PapayaWhip}{rgb}{1.00,0.94,0.84}
\definecolor{PeachPuff1}{rgb}{1.00,0.85,0.73}
\definecolor{PeachPuff2}{rgb}{0.93,0.80,0.68}
\definecolor{PeachPuff3}{rgb}{0.80,0.69,0.58}
\definecolor{PeachPuff4}{rgb}{0.55,0.47,0.40}
\definecolor{PeachPuff}{rgb}{1.00,0.85,0.73}
\definecolor{PowderBlue}{rgb}{0.69,0.88,0.90}
\definecolor{RosyBrown1}{rgb}{1.00,0.76,0.76}
\definecolor{RosyBrown2}{rgb}{0.93,0.71,0.71}
\definecolor{RosyBrown3}{rgb}{0.80,0.61,0.61}
\definecolor{RosyBrown4}{rgb}{0.55,0.41,0.41}
\definecolor{RosyBrown}{rgb}{0.74,0.56,0.56}
\definecolor{RoyalBlue1}{rgb}{0.28,0.46,1.00}
\definecolor{RoyalBlue2}{rgb}{0.26,0.43,0.93}
\definecolor{RoyalBlue3}{rgb}{0.23,0.37,0.80}
\definecolor{RoyalBlue4}{rgb}{0.15,0.25,0.55}
\definecolor{RoyalBlue}{rgb}{0.25,0.41,0.88}
\definecolor{SaddleBrown}{rgb}{0.55,0.27,0.07}
\definecolor{SandyBrown}{rgb}{0.96,0.64,0.38}
\definecolor{SeaGreen1}{rgb}{0.33,1.00,0.62}
\definecolor{SeaGreen2}{rgb}{0.31,0.93,0.58}
\definecolor{SeaGreen3}{rgb}{0.26,0.80,0.50}
\definecolor{SeaGreen4}{rgb}{0.18,0.55,0.34}
\definecolor{SeaGreen}{rgb}{0.18,0.55,0.34}
\definecolor{SkyBlue1}{rgb}{0.53,0.81,1.00}
\definecolor{SkyBlue2}{rgb}{0.49,0.75,0.93}
\definecolor{SkyBlue3}{rgb}{0.42,0.65,0.80}
\definecolor{SkyBlue4}{rgb}{0.29,0.44,0.55}
\definecolor{SkyBlue}{rgb}{0.53,0.81,0.92}
\definecolor{SlateBlue1}{rgb}{0.51,0.44,1.00}
\definecolor{SlateBlue2}{rgb}{0.48,0.40,0.93}
\definecolor{SlateBlue3}{rgb}{0.41,0.35,0.80}
\definecolor{SlateBlue4}{rgb}{0.28,0.24,0.55}
\definecolor{SlateBlue}{rgb}{0.42,0.35,0.80}
\definecolor{SlateGray1}{rgb}{0.78,0.89,1.00}
\definecolor{SlateGray2}{rgb}{0.73,0.83,0.93}
\definecolor{SlateGray3}{rgb}{0.62,0.71,0.80}
\definecolor{SlateGray4}{rgb}{0.42,0.48,0.55}
\definecolor{SlateGray}{rgb}{0.44,0.50,0.56}
\definecolor{SlateGrey}{rgb}{0.44,0.50,0.56}
\definecolor{SpringGreen1}{rgb}{0.00,1.00,0.50}
\definecolor{SpringGreen2}{rgb}{0.00,0.93,0.46}
\definecolor{SpringGreen3}{rgb}{0.00,0.80,0.40}
\definecolor{SpringGreen4}{rgb}{0.00,0.55,0.27}
\definecolor{SpringGreen}{rgb}{0.00,1.00,0.50}
\definecolor{SteelBlue1}{rgb}{0.39,0.72,1.00}
\definecolor{SteelBlue2}{rgb}{0.36,0.67,0.93}
\definecolor{SteelBlue3}{rgb}{0.31,0.58,0.80}
\definecolor{SteelBlue4}{rgb}{0.21,0.39,0.55}
\definecolor{SteelBlue}{rgb}{0.27,0.51,0.71}
\definecolor{VioletRed1}{rgb}{1.00,0.24,0.59}
\definecolor{VioletRed2}{rgb}{0.93,0.23,0.55}
\definecolor{VioletRed3}{rgb}{0.80,0.20,0.47}
\definecolor{VioletRed4}{rgb}{0.55,0.13,0.32}
\definecolor{VioletRed}{rgb}{0.82,0.13,0.56}
\definecolor{WhiteSmoke}{rgb}{0.96,0.96,0.96}
\definecolor{YellowGreen}{rgb}{0.60,0.80,0.20}
\definecolor{aliceblue}{rgb}{0.94,0.97,1.00}
\definecolor{antiquewhite}{rgb}{0.98,0.92,0.84}
\definecolor{aquamarine1}{rgb}{0.50,1.00,0.83}
\definecolor{aquamarine2}{rgb}{0.46,0.93,0.78}
\definecolor{aquamarine3}{rgb}{0.40,0.80,0.67}
\definecolor{aquamarine4}{rgb}{0.27,0.55,0.45}
\definecolor{aquamarine}{rgb}{0.50,1.00,0.83}
\definecolor{azure1}{rgb}{0.94,1.00,1.00}
\definecolor{azure2}{rgb}{0.88,0.93,0.93}
\definecolor{azure3}{rgb}{0.76,0.80,0.80}
\definecolor{azure4}{rgb}{0.51,0.55,0.55}
\definecolor{azure}{rgb}{0.94,1.00,1.00}
\definecolor{beige}{rgb}{0.96,0.96,0.86}
\definecolor{bisque1}{rgb}{1.00,0.89,0.77}
\definecolor{bisque2}{rgb}{0.93,0.84,0.72}
\definecolor{bisque3}{rgb}{0.80,0.72,0.62}
\definecolor{bisque4}{rgb}{0.55,0.49,0.42}
\definecolor{bisque}{rgb}{1.00,0.89,0.77}
\definecolor{black}{rgb}{0.00,0.00,0.00}
\definecolor{blanchedalmond}{rgb}{1.00,0.92,0.80}
\definecolor{blue1}{rgb}{0.00,0.00,1.00}
\definecolor{blue2}{rgb}{0.00,0.00,0.93}
\definecolor{blue3}{rgb}{0.00,0.00,0.80}
\definecolor{blue4}{rgb}{0.00,0.00,0.55}
\definecolor{blueviolet}{rgb}{0.54,0.17,0.89}
\definecolor{blue}{rgb}{0.00,0.00,1.00}
\definecolor{brown1}{rgb}{1.00,0.25,0.25}
\definecolor{brown2}{rgb}{0.93,0.23,0.23}
\definecolor{brown3}{rgb}{0.80,0.20,0.20}
\definecolor{brown4}{rgb}{0.55,0.14,0.14}
\definecolor{brown}{rgb}{0.65,0.16,0.16}
\definecolor{burlywood1}{rgb}{1.00,0.83,0.61}
\definecolor{burlywood2}{rgb}{0.93,0.77,0.57}
\definecolor{burlywood3}{rgb}{0.80,0.67,0.49}
\definecolor{burlywood4}{rgb}{0.55,0.45,0.33}
\definecolor{burlywood}{rgb}{0.87,0.72,0.53}
\definecolor{cadetblue}{rgb}{0.37,0.62,0.63}
\definecolor{chartreuse1}{rgb}{0.50,1.00,0.00}
\definecolor{chartreuse2}{rgb}{0.46,0.93,0.00}
\definecolor{chartreuse3}{rgb}{0.40,0.80,0.00}
\definecolor{chartreuse4}{rgb}{0.27,0.55,0.00}
\definecolor{chartreuse}{rgb}{0.50,1.00,0.00}
\definecolor{chocolate1}{rgb}{1.00,0.50,0.14}
\definecolor{chocolate2}{rgb}{0.93,0.46,0.13}
\definecolor{chocolate3}{rgb}{0.80,0.40,0.11}
\definecolor{chocolate4}{rgb}{0.55,0.27,0.07}
\definecolor{chocolate}{rgb}{0.82,0.41,0.12}
\definecolor{coral1}{rgb}{1.00,0.45,0.34}
\definecolor{coral2}{rgb}{0.93,0.42,0.31}
\definecolor{coral3}{rgb}{0.80,0.36,0.27}
\definecolor{coral4}{rgb}{0.55,0.24,0.18}
\definecolor{coral}{rgb}{1.00,0.50,0.31}
\definecolor{cornflowerblue}{rgb}{0.39,0.58,0.93}
\definecolor{cornsilk1}{rgb}{1.00,0.97,0.86}
\definecolor{cornsilk2}{rgb}{0.93,0.91,0.80}
\definecolor{cornsilk3}{rgb}{0.80,0.78,0.69}
\definecolor{cornsilk4}{rgb}{0.55,0.53,0.47}
\definecolor{cornsilk}{rgb}{1.00,0.97,0.86}
\definecolor{cyan1}{rgb}{0.00,1.00,1.00}
\definecolor{cyan2}{rgb}{0.00,0.93,0.93}
\definecolor{cyan3}{rgb}{0.00,0.80,0.80}
\definecolor{cyan4}{rgb}{0.00,0.55,0.55}
\definecolor{cyan}{rgb}{0.00,1.00,1.00}
\definecolor{darkblue}{rgb}{0.00,0.00,0.55}
\definecolor{darkcyan}{rgb}{0.00,0.55,0.55}
\definecolor{darkgoldenrod}{rgb}{0.72,0.53,0.04}
\definecolor{darkgray}{rgb}{0.66,0.66,0.66}
\definecolor{darkgreen}{rgb}{0.00,0.39,0.00}
\definecolor{darkgrey}{rgb}{0.66,0.66,0.66}
\definecolor{darkkhaki}{rgb}{0.74,0.72,0.42}
\definecolor{darkmagenta}{rgb}{0.55,0.00,0.55}
\definecolor{darkolive}{rgb}{0.33,0.42,0.18}
\definecolor{darkorange}{rgb}{1.00,0.55,0.00}
\definecolor{darkorchid}{rgb}{0.60,0.20,0.80}
\definecolor{darkred}{rgb}{0.55,0.00,0.00}
\definecolor{darksalmon}{rgb}{0.91,0.59,0.48}
\definecolor{darksea}{rgb}{0.56,0.74,0.56}
\definecolor{darkslate}{rgb}{0.18,0.31,0.31}
\definecolor{darkslate}{rgb}{0.18,0.31,0.31}
\definecolor{darkslate}{rgb}{0.28,0.24,0.55}
\definecolor{darkturquoise}{rgb}{0.00,0.81,0.82}
\definecolor{darkviolet}{rgb}{0.58,0.00,0.83}
\definecolor{deeppink}{rgb}{1.00,0.08,0.58}
\definecolor{deepsky}{rgb}{0.00,0.75,1.00}
\definecolor{dimgray}{rgb}{0.41,0.41,0.41}
\definecolor{dimgrey}{rgb}{0.41,0.41,0.41}
\definecolor{dodgerblue}{rgb}{0.12,0.56,1.00}
\definecolor{firebrick1}{rgb}{1.00,0.19,0.19}
\definecolor{firebrick2}{rgb}{0.93,0.17,0.17}
\definecolor{firebrick3}{rgb}{0.80,0.15,0.15}
\definecolor{firebrick4}{rgb}{0.55,0.10,0.10}
\definecolor{firebrick}{rgb}{0.70,0.13,0.13}
\definecolor{floralwhite}{rgb}{1.00,0.98,0.94}
\definecolor{forestgreen}{rgb}{0.13,0.55,0.13}
\definecolor{gainsboro}{rgb}{0.86,0.86,0.86}
\definecolor{ghostwhite}{rgb}{0.97,0.97,1.00}
\definecolor{gold1}{rgb}{1.00,0.84,0.00}
\definecolor{gold2}{rgb}{0.93,0.79,0.00}
\definecolor{gold3}{rgb}{0.80,0.68,0.00}
\definecolor{gold4}{rgb}{0.55,0.46,0.00}
\definecolor{goldenrod1}{rgb}{1.00,0.76,0.15}
\definecolor{goldenrod2}{rgb}{0.93,0.71,0.13}
\definecolor{goldenrod3}{rgb}{0.80,0.61,0.11}
\definecolor{goldenrod4}{rgb}{0.55,0.41,0.08}
\definecolor{goldenrod}{rgb}{0.85,0.65,0.13}
\definecolor{gold}{rgb}{1.00,0.84,0.00}
\definecolor{gray0}{rgb}{0.00,0.00,0.00}
\definecolor{gray100}{rgb}{1.00,1.00,1.00}
\definecolor{gray10}{rgb}{0.10,0.10,0.10}
\definecolor{gray11}{rgb}{0.11,0.11,0.11}
\definecolor{gray12}{rgb}{0.12,0.12,0.12}
\definecolor{gray13}{rgb}{0.13,0.13,0.13}
\definecolor{gray14}{rgb}{0.14,0.14,0.14}
\definecolor{gray15}{rgb}{0.15,0.15,0.15}
\definecolor{gray16}{rgb}{0.16,0.16,0.16}
\definecolor{gray17}{rgb}{0.17,0.17,0.17}
\definecolor{gray18}{rgb}{0.18,0.18,0.18}
\definecolor{gray19}{rgb}{0.19,0.19,0.19}
\definecolor{gray1}{rgb}{0.01,0.01,0.01}
\definecolor{gray20}{rgb}{0.20,0.20,0.20}
\definecolor{gray21}{rgb}{0.21,0.21,0.21}
\definecolor{gray22}{rgb}{0.22,0.22,0.22}
\definecolor{gray23}{rgb}{0.23,0.23,0.23}
\definecolor{gray24}{rgb}{0.24,0.24,0.24}
\definecolor{gray25}{rgb}{0.25,0.25,0.25}
\definecolor{gray26}{rgb}{0.26,0.26,0.26}
\definecolor{gray27}{rgb}{0.27,0.27,0.27}
\definecolor{gray28}{rgb}{0.28,0.28,0.28}
\definecolor{gray29}{rgb}{0.29,0.29,0.29}
\definecolor{gray2}{rgb}{0.02,0.02,0.02}
\definecolor{gray30}{rgb}{0.30,0.30,0.30}
\definecolor{gray31}{rgb}{0.31,0.31,0.31}
\definecolor{gray32}{rgb}{0.32,0.32,0.32}
\definecolor{gray33}{rgb}{0.33,0.33,0.33}
\definecolor{gray34}{rgb}{0.34,0.34,0.34}
\definecolor{gray35}{rgb}{0.35,0.35,0.35}
\definecolor{gray36}{rgb}{0.36,0.36,0.36}
\definecolor{gray37}{rgb}{0.37,0.37,0.37}
\definecolor{gray38}{rgb}{0.38,0.38,0.38}
\definecolor{gray39}{rgb}{0.39,0.39,0.39}
\definecolor{gray3}{rgb}{0.03,0.03,0.03}
\definecolor{gray40}{rgb}{0.40,0.40,0.40}
\definecolor{gray41}{rgb}{0.41,0.41,0.41}
\definecolor{gray42}{rgb}{0.42,0.42,0.42}
\definecolor{gray43}{rgb}{0.43,0.43,0.43}
\definecolor{gray44}{rgb}{0.44,0.44,0.44}
\definecolor{gray45}{rgb}{0.45,0.45,0.45}
\definecolor{gray46}{rgb}{0.46,0.46,0.46}
\definecolor{gray47}{rgb}{0.47,0.47,0.47}
\definecolor{gray48}{rgb}{0.48,0.48,0.48}
\definecolor{gray49}{rgb}{0.49,0.49,0.49}
\definecolor{gray4}{rgb}{0.04,0.04,0.04}
\definecolor{gray50}{rgb}{0.50,0.50,0.50}
\definecolor{gray51}{rgb}{0.51,0.51,0.51}
\definecolor{gray52}{rgb}{0.52,0.52,0.52}
\definecolor{gray53}{rgb}{0.53,0.53,0.53}
\definecolor{gray54}{rgb}{0.54,0.54,0.54}
\definecolor{gray55}{rgb}{0.55,0.55,0.55}
\definecolor{gray56}{rgb}{0.56,0.56,0.56}
\definecolor{gray57}{rgb}{0.57,0.57,0.57}
\definecolor{gray58}{rgb}{0.58,0.58,0.58}
\definecolor{gray59}{rgb}{0.59,0.59,0.59}
\definecolor{gray5}{rgb}{0.05,0.05,0.05}
\definecolor{gray60}{rgb}{0.60,0.60,0.60}
\definecolor{gray61}{rgb}{0.61,0.61,0.61}
\definecolor{gray62}{rgb}{0.62,0.62,0.62}
\definecolor{gray63}{rgb}{0.63,0.63,0.63}
\definecolor{gray64}{rgb}{0.64,0.64,0.64}
\definecolor{gray65}{rgb}{0.65,0.65,0.65}
\definecolor{gray66}{rgb}{0.66,0.66,0.66}
\definecolor{gray67}{rgb}{0.67,0.67,0.67}
\definecolor{gray68}{rgb}{0.68,0.68,0.68}
\definecolor{gray69}{rgb}{0.69,0.69,0.69}
\definecolor{gray6}{rgb}{0.06,0.06,0.06}
\definecolor{gray70}{rgb}{0.70,0.70,0.70}
\definecolor{gray71}{rgb}{0.71,0.71,0.71}
\definecolor{gray72}{rgb}{0.72,0.72,0.72}
\definecolor{gray73}{rgb}{0.73,0.73,0.73}
\definecolor{gray74}{rgb}{0.74,0.74,0.74}
\definecolor{gray75}{rgb}{0.75,0.75,0.75}
\definecolor{gray76}{rgb}{0.76,0.76,0.76}
\definecolor{gray77}{rgb}{0.77,0.77,0.77}
\definecolor{gray78}{rgb}{0.78,0.78,0.78}
\definecolor{gray79}{rgb}{0.79,0.79,0.79}
\definecolor{gray7}{rgb}{0.07,0.07,0.07}
\definecolor{gray80}{rgb}{0.80,0.80,0.80}
\definecolor{gray81}{rgb}{0.81,0.81,0.81}
\definecolor{gray82}{rgb}{0.82,0.82,0.82}
\definecolor{gray83}{rgb}{0.83,0.83,0.83}
\definecolor{gray84}{rgb}{0.84,0.84,0.84}
\definecolor{gray85}{rgb}{0.85,0.85,0.85}
\definecolor{gray86}{rgb}{0.86,0.86,0.86}
\definecolor{gray87}{rgb}{0.87,0.87,0.87}
\definecolor{gray88}{rgb}{0.88,0.88,0.88}
\definecolor{gray89}{rgb}{0.89,0.89,0.89}
\definecolor{gray8}{rgb}{0.08,0.08,0.08}
\definecolor{gray90}{rgb}{0.90,0.90,0.90}
\definecolor{gray91}{rgb}{0.91,0.91,0.91}
\definecolor{gray92}{rgb}{0.92,0.92,0.92}
\definecolor{gray93}{rgb}{0.93,0.93,0.93}
\definecolor{gray94}{rgb}{0.94,0.94,0.94}
\definecolor{gray95}{rgb}{0.95,0.95,0.95}
\definecolor{gray96}{rgb}{0.96,0.96,0.96}
\definecolor{gray97}{rgb}{0.97,0.97,0.97}
\definecolor{gray98}{rgb}{0.98,0.98,0.98}
\definecolor{gray99}{rgb}{0.99,0.99,0.99}
\definecolor{gray9}{rgb}{0.09,0.09,0.09}
\definecolor{gray}{rgb}{0.75,0.75,0.75}
\definecolor{green1}{rgb}{0.00,1.00,0.00}
\definecolor{green2}{rgb}{0.00,0.93,0.00}
\definecolor{green3}{rgb}{0.00,0.80,0.00}
\definecolor{green4}{rgb}{0.00,0.55,0.00}
\definecolor{greenyellow}{rgb}{0.68,1.00,0.18}
\definecolor{green}{rgb}{0.00,1.00,0.00}
\definecolor{grey0}{rgb}{0.00,0.00,0.00}
\definecolor{grey100}{rgb}{1.00,1.00,1.00}
\definecolor{grey10}{rgb}{0.10,0.10,0.10}
\definecolor{grey11}{rgb}{0.11,0.11,0.11}
\definecolor{grey12}{rgb}{0.12,0.12,0.12}
\definecolor{grey13}{rgb}{0.13,0.13,0.13}
\definecolor{grey14}{rgb}{0.14,0.14,0.14}
\definecolor{grey15}{rgb}{0.15,0.15,0.15}
\definecolor{grey16}{rgb}{0.16,0.16,0.16}
\definecolor{grey17}{rgb}{0.17,0.17,0.17}
\definecolor{grey18}{rgb}{0.18,0.18,0.18}
\definecolor{grey19}{rgb}{0.19,0.19,0.19}
\definecolor{grey1}{rgb}{0.01,0.01,0.01}
\definecolor{grey20}{rgb}{0.20,0.20,0.20}
\definecolor{grey21}{rgb}{0.21,0.21,0.21}
\definecolor{grey22}{rgb}{0.22,0.22,0.22}
\definecolor{grey23}{rgb}{0.23,0.23,0.23}
\definecolor{grey24}{rgb}{0.24,0.24,0.24}
\definecolor{grey25}{rgb}{0.25,0.25,0.25}
\definecolor{grey26}{rgb}{0.26,0.26,0.26}
\definecolor{grey27}{rgb}{0.27,0.27,0.27}
\definecolor{grey28}{rgb}{0.28,0.28,0.28}
\definecolor{grey29}{rgb}{0.29,0.29,0.29}
\definecolor{grey2}{rgb}{0.02,0.02,0.02}
\definecolor{grey30}{rgb}{0.30,0.30,0.30}
\definecolor{grey31}{rgb}{0.31,0.31,0.31}
\definecolor{grey32}{rgb}{0.32,0.32,0.32}
\definecolor{grey33}{rgb}{0.33,0.33,0.33}
\definecolor{grey34}{rgb}{0.34,0.34,0.34}
\definecolor{grey35}{rgb}{0.35,0.35,0.35}
\definecolor{grey36}{rgb}{0.36,0.36,0.36}
\definecolor{grey37}{rgb}{0.37,0.37,0.37}
\definecolor{grey38}{rgb}{0.38,0.38,0.38}
\definecolor{grey39}{rgb}{0.39,0.39,0.39}
\definecolor{grey3}{rgb}{0.03,0.03,0.03}
\definecolor{grey40}{rgb}{0.40,0.40,0.40}
\definecolor{grey41}{rgb}{0.41,0.41,0.41}
\definecolor{grey42}{rgb}{0.42,0.42,0.42}
\definecolor{grey43}{rgb}{0.43,0.43,0.43}
\definecolor{grey44}{rgb}{0.44,0.44,0.44}
\definecolor{grey45}{rgb}{0.45,0.45,0.45}
\definecolor{grey46}{rgb}{0.46,0.46,0.46}
\definecolor{grey47}{rgb}{0.47,0.47,0.47}
\definecolor{grey48}{rgb}{0.48,0.48,0.48}
\definecolor{grey49}{rgb}{0.49,0.49,0.49}
\definecolor{grey4}{rgb}{0.04,0.04,0.04}
\definecolor{grey50}{rgb}{0.50,0.50,0.50}
\definecolor{grey51}{rgb}{0.51,0.51,0.51}
\definecolor{grey52}{rgb}{0.52,0.52,0.52}
\definecolor{grey53}{rgb}{0.53,0.53,0.53}
\definecolor{grey54}{rgb}{0.54,0.54,0.54}
\definecolor{grey55}{rgb}{0.55,0.55,0.55}
\definecolor{grey56}{rgb}{0.56,0.56,0.56}
\definecolor{grey57}{rgb}{0.57,0.57,0.57}
\definecolor{grey58}{rgb}{0.58,0.58,0.58}
\definecolor{grey59}{rgb}{0.59,0.59,0.59}
\definecolor{grey5}{rgb}{0.05,0.05,0.05}
\definecolor{grey60}{rgb}{0.60,0.60,0.60}
\definecolor{grey61}{rgb}{0.61,0.61,0.61}
\definecolor{grey62}{rgb}{0.62,0.62,0.62}
\definecolor{grey63}{rgb}{0.63,0.63,0.63}
\definecolor{grey64}{rgb}{0.64,0.64,0.64}
\definecolor{grey65}{rgb}{0.65,0.65,0.65}
\definecolor{grey66}{rgb}{0.66,0.66,0.66}
\definecolor{grey67}{rgb}{0.67,0.67,0.67}
\definecolor{grey68}{rgb}{0.68,0.68,0.68}
\definecolor{grey69}{rgb}{0.69,0.69,0.69}
\definecolor{grey6}{rgb}{0.06,0.06,0.06}
\definecolor{grey70}{rgb}{0.70,0.70,0.70}
\definecolor{grey71}{rgb}{0.71,0.71,0.71}
\definecolor{grey72}{rgb}{0.72,0.72,0.72}
\definecolor{grey73}{rgb}{0.73,0.73,0.73}
\definecolor{grey74}{rgb}{0.74,0.74,0.74}
\definecolor{grey75}{rgb}{0.75,0.75,0.75}
\definecolor{grey76}{rgb}{0.76,0.76,0.76}
\definecolor{grey77}{rgb}{0.77,0.77,0.77}
\definecolor{grey78}{rgb}{0.78,0.78,0.78}
\definecolor{grey79}{rgb}{0.79,0.79,0.79}
\definecolor{grey7}{rgb}{0.07,0.07,0.07}
\definecolor{grey80}{rgb}{0.80,0.80,0.80}
\definecolor{grey81}{rgb}{0.81,0.81,0.81}
\definecolor{grey82}{rgb}{0.82,0.82,0.82}
\definecolor{grey83}{rgb}{0.83,0.83,0.83}
\definecolor{grey84}{rgb}{0.84,0.84,0.84}
\definecolor{grey85}{rgb}{0.85,0.85,0.85}
\definecolor{grey86}{rgb}{0.86,0.86,0.86}
\definecolor{grey87}{rgb}{0.87,0.87,0.87}
\definecolor{grey88}{rgb}{0.88,0.88,0.88}
\definecolor{grey89}{rgb}{0.89,0.89,0.89}
\definecolor{grey8}{rgb}{0.08,0.08,0.08}
\definecolor{grey90}{rgb}{0.90,0.90,0.90}
\definecolor{grey91}{rgb}{0.91,0.91,0.91}
\definecolor{grey92}{rgb}{0.92,0.92,0.92}
\definecolor{grey93}{rgb}{0.93,0.93,0.93}
\definecolor{grey94}{rgb}{0.94,0.94,0.94}
\definecolor{grey95}{rgb}{0.95,0.95,0.95}
\definecolor{grey96}{rgb}{0.96,0.96,0.96}
\definecolor{grey97}{rgb}{0.97,0.97,0.97}
\definecolor{grey98}{rgb}{0.98,0.98,0.98}
\definecolor{grey99}{rgb}{0.99,0.99,0.99}
\definecolor{grey9}{rgb}{0.09,0.09,0.09}
\definecolor{grey}{rgb}{0.75,0.75,0.75}
\definecolor{honeydew1}{rgb}{0.94,1.00,0.94}
\definecolor{honeydew2}{rgb}{0.88,0.93,0.88}
\definecolor{honeydew3}{rgb}{0.76,0.80,0.76}
\definecolor{honeydew4}{rgb}{0.51,0.55,0.51}
\definecolor{honeydew}{rgb}{0.94,1.00,0.94}
\definecolor{hotpink}{rgb}{1.00,0.41,0.71}
\definecolor{indianred}{rgb}{0.80,0.36,0.36}
\definecolor{ivory1}{rgb}{1.00,1.00,0.94}
\definecolor{ivory2}{rgb}{0.93,0.93,0.88}
\definecolor{ivory3}{rgb}{0.80,0.80,0.76}
\definecolor{ivory4}{rgb}{0.55,0.55,0.51}
\definecolor{ivory}{rgb}{1.00,1.00,0.94}
\definecolor{khaki1}{rgb}{1.00,0.96,0.56}
\definecolor{khaki2}{rgb}{0.93,0.90,0.52}
\definecolor{khaki3}{rgb}{0.80,0.78,0.45}
\definecolor{khaki4}{rgb}{0.55,0.53,0.31}
\definecolor{khaki}{rgb}{0.94,0.90,0.55}
\definecolor{lavenderblush}{rgb}{1.00,0.94,0.96}
\definecolor{lavender}{rgb}{0.90,0.90,0.98}
\definecolor{lawngreen}{rgb}{0.49,0.99,0.00}
\definecolor{lemonchiffon}{rgb}{1.00,0.98,0.80}
\definecolor{lightblue}{rgb}{0.68,0.85,0.90}
\definecolor{lightcoral}{rgb}{0.94,0.50,0.50}
\definecolor{lightcyan}{rgb}{0.88,1.00,1.00}
\definecolor{lightgoldenrod}{rgb}{0.93,0.87,0.51}
\definecolor{lightgoldenrod}{rgb}{0.98,0.98,0.82}
\definecolor{lightgray}{rgb}{0.83,0.83,0.83}
\definecolor{lightgreen}{rgb}{0.56,0.93,0.56}
\definecolor{lightgrey}{rgb}{0.83,0.83,0.83}
\definecolor{lightpink}{rgb}{1.00,0.71,0.76}
\definecolor{lightsalmon}{rgb}{1.00,0.63,0.48}
\definecolor{lightsea}{rgb}{0.13,0.70,0.67}
\definecolor{lightsky}{rgb}{0.53,0.81,0.98}
\definecolor{lightslate}{rgb}{0.47,0.53,0.60}
\definecolor{lightslate}{rgb}{0.47,0.53,0.60}
\definecolor{lightslate}{rgb}{0.52,0.44,1.00}
\definecolor{lightsteel}{rgb}{0.69,0.77,0.87}
\definecolor{lightyellow}{rgb}{1.00,1.00,0.88}
\definecolor{limegreen}{rgb}{0.20,0.80,0.20}
\definecolor{linen}{rgb}{0.98,0.94,0.90}
\definecolor{magenta1}{rgb}{1.00,0.00,1.00}
\definecolor{magenta2}{rgb}{0.93,0.00,0.93}
\definecolor{magenta3}{rgb}{0.80,0.00,0.80}
\definecolor{magenta4}{rgb}{0.55,0.00,0.55}
\definecolor{magenta}{rgb}{1.00,0.00,1.00}
\definecolor{maroon1}{rgb}{1.00,0.20,0.70}
\definecolor{maroon2}{rgb}{0.93,0.19,0.65}
\definecolor{maroon3}{rgb}{0.80,0.16,0.56}
\definecolor{maroon4}{rgb}{0.55,0.11,0.38}
\definecolor{maroon}{rgb}{0.69,0.19,0.38}
\definecolor{mediumaquamarine}{rgb}{0.40,0.80,0.67}
\definecolor{mediumblue}{rgb}{0.00,0.00,0.80}
\definecolor{mediumorchid}{rgb}{0.73,0.33,0.83}
\definecolor{mediumpurple}{rgb}{0.58,0.44,0.86}
\definecolor{mediumsea}{rgb}{0.24,0.70,0.44}
\definecolor{mediumslate}{rgb}{0.48,0.41,0.93}
\definecolor{mediumspring}{rgb}{0.00,0.98,0.60}
\definecolor{mediumturquoise}{rgb}{0.28,0.82,0.80}
\definecolor{mediumviolet}{rgb}{0.78,0.08,0.52}
\definecolor{midnightblue}{rgb}{0.10,0.10,0.44}
\definecolor{mintcream}{rgb}{0.96,1.00,0.98}
\definecolor{mistyrose}{rgb}{1.00,0.89,0.88}
\definecolor{moccasin}{rgb}{1.00,0.89,0.71}
\definecolor{navajowhite}{rgb}{1.00,0.87,0.68}
\definecolor{navyblue}{rgb}{0.00,0.00,0.50}
\definecolor{navy}{rgb}{0.00,0.00,0.50}
\definecolor{oldlace}{rgb}{0.99,0.96,0.90}
\definecolor{olivedrab}{rgb}{0.42,0.56,0.14}
\definecolor{orange1}{rgb}{1.00,0.65,0.00}
\definecolor{orange2}{rgb}{0.93,0.60,0.00}
\definecolor{orange3}{rgb}{0.80,0.52,0.00}
\definecolor{orange4}{rgb}{0.55,0.35,0.00}
\definecolor{orangered}{rgb}{1.00,0.27,0.00}
\definecolor{orange}{rgb}{1.00,0.65,0.00}
\definecolor{orchid1}{rgb}{1.00,0.51,0.98}
\definecolor{orchid2}{rgb}{0.93,0.48,0.91}
\definecolor{orchid3}{rgb}{0.80,0.41,0.79}
\definecolor{orchid4}{rgb}{0.55,0.28,0.54}
\definecolor{orchid}{rgb}{0.85,0.44,0.84}
\definecolor{palegoldenrod}{rgb}{0.93,0.91,0.67}
\definecolor{palegreen}{rgb}{0.60,0.98,0.60}
\definecolor{paleturquoise}{rgb}{0.69,0.93,0.93}
\definecolor{paleviolet}{rgb}{0.86,0.44,0.58}
\definecolor{papayawhip}{rgb}{1.00,0.94,0.84}
\definecolor{peachpuff}{rgb}{1.00,0.85,0.73}
\definecolor{peru}{rgb}{0.80,0.52,0.25}
\definecolor{pink1}{rgb}{1.00,0.71,0.77}
\definecolor{pink2}{rgb}{0.93,0.66,0.72}
\definecolor{pink3}{rgb}{0.80,0.57,0.62}
\definecolor{pink4}{rgb}{0.55,0.39,0.42}
\definecolor{pink}{rgb}{1.00,0.75,0.80}
\definecolor{plum1}{rgb}{1.00,0.73,1.00}
\definecolor{plum2}{rgb}{0.93,0.68,0.93}
\definecolor{plum3}{rgb}{0.80,0.59,0.80}
\definecolor{plum4}{rgb}{0.55,0.40,0.55}
\definecolor{plum}{rgb}{0.87,0.63,0.87}
\definecolor{powderblue}{rgb}{0.69,0.88,0.90}
\definecolor{purple1}{rgb}{0.61,0.19,1.00}
\definecolor{purple2}{rgb}{0.57,0.17,0.93}
\definecolor{purple3}{rgb}{0.49,0.15,0.80}
\definecolor{purple4}{rgb}{0.33,0.10,0.55}
\definecolor{purple}{rgb}{0.63,0.13,0.94}
\definecolor{red1}{rgb}{1.00,0.00,0.00}
\definecolor{red2}{rgb}{0.93,0.00,0.00}
\definecolor{red3}{rgb}{0.80,0.00,0.00}
\definecolor{red4}{rgb}{0.55,0.00,0.00}
\definecolor{red}{rgb}{1.00,0.00,0.00}
\definecolor{rosybrown}{rgb}{0.74,0.56,0.56}
\definecolor{royalblue}{rgb}{0.25,0.41,0.88}
\definecolor{saddlebrown}{rgb}{0.55,0.27,0.07}
\definecolor{salmon1}{rgb}{1.00,0.55,0.41}
\definecolor{salmon2}{rgb}{0.93,0.51,0.38}
\definecolor{salmon3}{rgb}{0.80,0.44,0.33}
\definecolor{salmon4}{rgb}{0.55,0.30,0.22}
\definecolor{salmon}{rgb}{0.98,0.50,0.45}
\definecolor{sandybrown}{rgb}{0.96,0.64,0.38}
\definecolor{seagreen}{rgb}{0.18,0.55,0.34}
\definecolor{seashell1}{rgb}{1.00,0.96,0.93}
\definecolor{seashell2}{rgb}{0.93,0.90,0.87}
\definecolor{seashell3}{rgb}{0.80,0.77,0.75}
\definecolor{seashell4}{rgb}{0.55,0.53,0.51}
\definecolor{seashell}{rgb}{1.00,0.96,0.93}
\definecolor{sienna1}{rgb}{1.00,0.51,0.28}
\definecolor{sienna2}{rgb}{0.93,0.47,0.26}
\definecolor{sienna3}{rgb}{0.80,0.41,0.22}
\definecolor{sienna4}{rgb}{0.55,0.28,0.15}
\definecolor{sienna}{rgb}{0.63,0.32,0.18}
\definecolor{skyblue}{rgb}{0.53,0.81,0.92}
\definecolor{slateblue}{rgb}{0.42,0.35,0.80}
\definecolor{slategray}{rgb}{0.44,0.50,0.56}
\definecolor{slategrey}{rgb}{0.44,0.50,0.56}
\definecolor{snow1}{rgb}{1.00,0.98,0.98}
\definecolor{snow2}{rgb}{0.93,0.91,0.91}
\definecolor{snow3}{rgb}{0.80,0.79,0.79}
\definecolor{snow4}{rgb}{0.55,0.54,0.54}
\definecolor{snow}{rgb}{1.00,0.98,0.98}
\definecolor{springgreen}{rgb}{0.00,1.00,0.50}
\definecolor{steelblue}{rgb}{0.27,0.51,0.71}
\definecolor{tan1}{rgb}{1.00,0.65,0.31}
\definecolor{tan2}{rgb}{0.93,0.60,0.29}
\definecolor{tan3}{rgb}{0.80,0.52,0.25}
\definecolor{tan4}{rgb}{0.55,0.35,0.17}
\definecolor{tan}{rgb}{0.82,0.71,0.55}
\definecolor{thistle1}{rgb}{1.00,0.88,1.00}
\definecolor{thistle2}{rgb}{0.93,0.82,0.93}
\definecolor{thistle3}{rgb}{0.80,0.71,0.80}
\definecolor{thistle4}{rgb}{0.55,0.48,0.55}
\definecolor{thistle}{rgb}{0.85,0.75,0.85}
\definecolor{tomato1}{rgb}{1.00,0.39,0.28}
\definecolor{tomato2}{rgb}{0.93,0.36,0.26}
\definecolor{tomato3}{rgb}{0.80,0.31,0.22}
\definecolor{tomato4}{rgb}{0.55,0.21,0.15}
\definecolor{tomato}{rgb}{1.00,0.39,0.28}
\definecolor{turquoise1}{rgb}{0.00,0.96,1.00}
\definecolor{turquoise2}{rgb}{0.00,0.90,0.93}
\definecolor{turquoise3}{rgb}{0.00,0.77,0.80}
\definecolor{turquoise4}{rgb}{0.00,0.53,0.55}
\definecolor{turquoise}{rgb}{0.25,0.88,0.82}
\definecolor{violetred}{rgb}{0.82,0.13,0.56}
\definecolor{violet}{rgb}{0.93,0.51,0.93}
\definecolor{wheat1}{rgb}{1.00,0.91,0.73}
\definecolor{wheat2}{rgb}{0.93,0.85,0.68}
\definecolor{wheat3}{rgb}{0.80,0.73,0.59}
\definecolor{wheat4}{rgb}{0.55,0.49,0.40}
\definecolor{wheat}{rgb}{0.96,0.87,0.70}
\definecolor{whitesmoke}{rgb}{0.96,0.96,0.96}
\definecolor{white}{rgb}{1.00,1.00,1.00}
\definecolor{yellow1}{rgb}{1.00,1.00,0.00}
\definecolor{yellow2}{rgb}{0.93,0.93,0.00}
\definecolor{yellow3}{rgb}{0.80,0.80,0.00}
\definecolor{yellow4}{rgb}{0.55,0.55,0.00}
\definecolor{yellowgreen}{rgb}{0.60,0.80,0.20}
\definecolor{yellow}{rgb}{1.00,1.00,0.00}
 
\usepackage[usenames,dvipsnames]{xcolor}
\usepackage[psamsfonts]{amsfonts}
\usepackage{graphicx}
\usepackage{amssymb}

\usepackage{hyperref} 
\hypersetup{pdfpagemode=FullScreen,   
backref,
colorlinks=true, 
citecolor=Bittersweet,
linkcolor=Bittersweet, 
urlcolor=Bittersweet
}
 
\setattribute{journal}{name}{}

\startlocaldefs

\theoremstyle{plain}
\newtheorem{Lem}{Lemma}[section]
\newtheorem{Prop}{Proposition}[section]
\newtheorem{Theor}{Theorem}[section]

\numberwithin{equation}{section}

\newcommand{\n}{^{(n)}}

\newcommand{\cqfd}{\hfill $\square$}

\newcommand{\R}{\mathbb R}

\newcommand{\Xb}{\mathbf{X}}
\newcommand{\Sb}{\mathbf{S}}

\newcommand{\ub}{\ensuremath{\mathbf{u}}}

\newcommand{\vb}{\ensuremath{\mathbf{v}}}
\newcommand{\xb}{\ensuremath{\mathbf{x}}}

\newcommand{\Ab}{\ensuremath{\mathbf{A}}}

\newcommand{\Hb}{\ensuremath{\mathbf{H}}}
\newcommand{\Ob}{\ensuremath{\mathbf{O}}}

\newcommand{\Tb}{\ensuremath{\mathbf{T}}}

\newcommand{\Wb}{\ensuremath{\mathbf{W}}}
\newcommand{\Yb}{\ensuremath{\mathbf{Y}}}

\newcommand{\thetab}{{\pmb \theta}}

\newcommand{\Deltab}{{\pmb \Delta}}
\newcommand{\taub}{{\pmb \tau}}

\newcommand{\Gamb}{{\pmb \Gamma}}

\newcommand{\pr}{^{\prime}}

\newcommand{\ny}{n\rightarrow\infty}

\endlocaldefs

\begin{document}

\begin{frontmatter}
\title{Inference for spherical location\\ under high concentration}
\runtitle{Spherical inference under high concentration}

\begin{aug}
\author{\fnms{Davy} \snm{Paindaveine}\thanksref{t1}\ead[label=e1]{dpaindav@ulb.ac.be}
\ead[label=u1,url]{http://homepages.ulb.ac.be/dpaindav}}
\,\! \and\,\!
\author{\fnms{ Thomas} \snm{Verdebout}\thanksref{t2}
\ead[label=e2]{tverdebo@ulb.ac.be}
\ead[label=u2,url]{http://tverdebo.ulb.ac.be}}

\thankstext{t1}{Corresponding author. Davy Paindaveine's research is supported by a research fellowship from the Francqui Foundation and by the Program of Concerted Research Actions (ARC) of the Universit\'{e} libre de Bruxelles.}
\thankstext{t2}{Thomas Verdebout's research is supported by the ARC Program of the Universit\'{e} libre de Bruxelles and by the Cr\'{e}dit de Recherche  J.0134.18 of the FNRS (Fonds National pour la Recherche Scientifique), Communaut\'{e} Fran\c{c}aise de Belgique. }
\runauthor{D. Paindaveine and Th. Verdebout}

\affiliation{Universit\'{e} libre de Bruxelles}

\address{Universit\'{e} libre de Bruxelles\\
ECARES and D\'{e}partement de Math\'{e}matique\\
Avenue F.D. Roosevelt, 50\\ 
ECARES, CP114/04\\
B-1050, Brussels\\ 
Belgium\\
\printead{e1}\\  
\printead{u1}\\
}

\address{Universit\'{e} libre de Bruxelles\\
ECARES and D\'{e}partement de Math\'{e}matique\\
Boulevard du Triomphe, CP210\\
B-1050, Brussels\\
Belgium\\ 
\printead{e2}\\ 
\printead{u2}\\
}

\end{aug}
\vspace{3mm}

\begin{abstract}
Motivated by the fact that circular or spherical data are often much concentrated around a location $\pmb\theta$, we consider inference about~$\pmb\theta$ under \emph{high concentration} asymptotic scenarios for which the probability of any fixed spherical cap centered at $\pmb\theta$ converges to one as the sample size~$n$ diverges to infinity. Rather than restricting to Fisher--von Mises--Langevin distributions, we consider a much broader, semiparametric, class of rotationally symmetric distributions indexed by the location parameter $\pmb\theta$, a scalar concentration parameter $\kappa$ and a functional nuisance $f$. We determine the class of distributions for which high concentration is obtained as $\kappa$ diverges to infinity. For such distributions, we then consider inference (point estimation, confidence zone estimation, hypothesis testing) on $\pmb\theta$ in asymptotic scenarios where $\kappa_n$ diverges to infinity at an arbitrary rate with the sample size $n$. Our asymptotic investigation reveals that, interestingly, optimal inference procedures on $\pmb\theta$ show consistency rates that depend on $f$. Using asymptotics ``\`{a} la Le Cam", we show that the spherical mean is, at any $f$, a parametrically super-efficient estimator of~$\thetab$ and that the Watson and Wald tests for $\mathcal{H}_0:{\pmb\theta}={\pmb\theta}_0$ enjoy similar, non-standard, optimality properties. We illustrate our results through simulations and treat a real data example. On a technical point of view, our asymptotic derivations require challenging expansions of rotationally symmetric functionals for large arguments of~$f$.
\end{abstract}

\begin{keyword}[class=MSC]
\kwd[Primary ]{62E20, 62F30}
\kwd[; secondary ]{62F05, 62F12}
\end{keyword}

\begin{keyword}
\kwd{Concentrated distributions}
\kwd{Directional statistics}
\kwd{Le Cam's asymptotic theory of statistical experiments}
\kwd{Local asymptotic normality}
\kwd{Super-efficiency}
\end{keyword}

\end{frontmatter}


\section{Introduction} 
\label{sec:intro}

Directional statistics is concerned with data on the unit sphere~$\mathcal{S}^{p-1}=\{ \xb \in\R^p : \|\xb\|^2=\xb'\xb=1\}$ of~$\R^p$ 
or more generally on Riemannian manifolds such as a torus or an infinite cylinder. Directional data are present in many fields and have attracted a lot of attention in the last decade. Recent applications include analysis of magnetic remanence through copulae on product manifolds in \cite{Ju15}, analysis of animal movement using angular regression in \cite{Rietal16}, or analysis of flight trajectories through principal component analysis for functional data on $\mathcal{S}^{p-1}$ in \cite{DaiMu18}, to cite only a few. For an overview of the topic, we refer to \cite{MarJup2000} and \cite{LVbook2017}.

In this paper, we consider a class of distributions on~$\mathcal{S}^{p-1}$ admitting a density at~$\xb$ that is proportional to~$f(\kappa \xb'\thetab)$, where~$\thetab\in\mathcal{S}^{p-1}$, $\kappa>0$ and~$f$ is a monotone increasing function from~$\R$ to~$\R^+$ (throughout, densities on~$\mathcal{S}^{p-1}$ will be with respect to the surface area measure). The resulting distribution on the sphere will be denoted as~${\rm Rot}_p(\thetab,\kappa,f)$ to stress its \emph{rotational symmetry}: if~$\Xb\sim {\rm Rot}_p(\thetab,\kappa,f)$, then~$\Ob\Xb$ and~$\Xb$ are equal in distribution for any $p\times p$ orthogonal matrix~$\Ob$ such that~$\Ob\thetab=\thetab$. Clearly, $\thetab$ is the modal location on the sphere, hence plays the role of a location parameter. In contrast, $\kappa$ is a scale or \emph{concentration} parameter. This terminology is justified by the fact that, for many functions~$f$, the distribution~${\rm Rot}_p(\thetab,\kappa,f)$ becomes arbitrarily concentrated around~$\thetab$ as~$\kappa$ diverges to infinity; it is in particular so for the celebrated \emph{Fisher--von Mises--Langevin (FvML)} distributions, that are obtained with~$f=\exp$. FvML distributions play a central role in directional statistics, a role that can be compared to the one played by Gaussian distributions in classical multivariate setups. For instance, the responses of the circular/spherical regression models in \cite{Riv1986}, \cite{DoMa02}, \cite{Ashis13} and \cite{Ro14} are FvML with a location parameter that depends on the predictors. 

In most applications, the location parameter~$\thetab$ is the parameter of interest, whereas the concentration parameter~$\kappa$ and the infinite-dimensional parameter~$f$ are unspecified nuisances. The most classical estimator of~$\thetab$ is the \emph{spherical mean}, whereas the most celebrated test for~$\mathcal{H}_0:\thetab=\thetab_0$, where~$\thetab_0\in\mathcal{S}^{p-1}$ is fixed, is the Watson test (see Sections~\ref{sec:estimation} and~\ref{sec:testing}, respectively).
In the standard asymptotic scenario under which~$n$ diverges to infinity with~$\kappa$ fixed, the asymptotic properties of these procedures are well-known; see, e.g., \cite{MarJup2000}. In particular, the spherical mean is \mbox{root-$n$} consistent, whereas the Watson test shows non-trivial asymptotic powers under sequences of local alternatives of the form~$\mathcal{H}_1\n:\thetab=\thetab_n$ with~$\sqrt{n}\|\thetab_n-\thetab_0\|\to c>0$. 

In practice, the asymptotic results above are relevant in cases where the underlying concentration~$\kappa$ is neither too small nor too large. For small values of~$\kappa$, the \mbox{fixed-$\kappa$} asymptotic distribution of the spherical mean and the corresponding asymptotic null distribution of~$W_n$ only poorly approximate the exact distribution of these statistics, unless the sample size~$n$ at hand is extremely large. This motivates considering a double asymptotic scenario where~$\kappa=\kappa_n$ goes to zero as~$n$ diverges to infinity. The observations~$\Xb_{n1},\ldots,\Xb_{nn}$ are then assumed to form a random sample from the distribution~${\rm Rot}_p(\thetab,\kappa_n,f)$, with~$\kappa_n=o(1)$, which makes it here strictly necessary to consider triangular arrays of observations. Such a ``low-concentration double asymptotic scenario" was considered in \cite{PaiVer17b}, where it was proved that the faster~$\kappa_n$ goes to zero, the poorer the consistency rates of the aforementioned inference procedures. More precisely, (i) if~$\kappa_n=o(1)$ with~$\kappa_n\sqrt{n}\to\infty$, then~$\kappa_n\sqrt{n}(\hat{\thetab}_n-\thetab)$ is asymptotically normal, so that the consistency rate of the spherical mean deteriorates from~$\sqrt{n}$ (in the standard fixed-$\kappa$ case) to~$\kappa_n\sqrt{n}$ (in the present case); (ii) if~$\kappa_n=O(1/\sqrt{n})$, then the spherical mean is not consistent anymore. Similarly, in situation~(i), the Watson test shows non-trivial asymptotic powers under sequences of local alternatives of the form~$\mathcal{H}_1\n\!:\thetab=\thetab_n$ with~$\kappa_n\sqrt{n}\|\thetab_n-\thetab_0\|\to c>0$, and, in situation~(ii), there is no sequence of alternatives under which this test would be consistent. These behaviors of the spherical mean and of the Watson test are non-standard yet expected: as the concentration~$\kappa_n$ gets smaller, the distribution~${\rm Rot}_p(\thetab,\kappa_n,f)$ becomes increasingly closer to the uniform distribution on~$\mathcal{S}^{p-1}$ for which the parameter of interest~$\thetab$ is not identifiable. In other words, inference on~$\thetab$ is increasingly challenging as~$\kappa$ decreases to zero, which reflects in the deterioration of the consistency rates above. 

The situation for large concentrations~$\kappa$ is similar yet different. On the one hand, it is still so that a standard, fixed-$\kappa$, asymptotic analysis could in principle fail describing in a suitable way the finite-sample behaviors of the spherical mean and of the Watson test statistic under high concentration. On the other hand, inference about~$\thetab$ intuitively becomes increasingly easy as the distribution gets more and more concentrated around~$\thetab$, which should make it possible to define ``super-efficient" estimators and tests on~$\thetab$. Inference for ``concentrated" FvML distributions actually has already been quite much considered in the literature. One of the first papers tackling inference problems for the location parameter of FvML distributions under large values of~$\kappa$ is  \cite{Wat1984}, where asymptotic results as~$\kappa\to\infty$ with~$n$ fixed were derived. In the same asymptotic scenario, \cite{Riv1986} investigated the null limiting behavior of a goodness-of-fit test for FvML distributions, whereas \cite{Riv1989}, \cite{DoMa02} and \cite{Do03} considered spherical regression in a concentrated FvML setup. \cite{Ro14} analyzed concentrated data using a regression model with an FvML noise. \cite{FujWat1992} obtained the asymptotic null distributions of various test statistics for~$\mathcal{H}_0\!:\thetab=\thetab_0$ again as~$\kappa\to\infty$ with~$n$ fixed, and derived the asymptotic powers of the corresponding tests under appropriate sequences of local alternatives. Still in the framework of FvML distributions, \cite{Wat1996} reviewed point estimation and (one-sample and multi-sample) hypothesis testing in the standard asymptotic scenario where~$n\to\infty$ with~$\kappa$ fixed and in the concentrated scenario where~$\kappa\to\infty$ with~$n$ fixed. \cite{ArJu13} and \cite{ArJu18} considered estimation of ``highly concentrated rotations". Finally, \cite{Chi2003} considered inference for concentrated matrix FvML distributions, still in a setup where~$\kappa\to\infty$ with~$n$ fixed; see also \cite{Chi2003book}. Monographs covering inference for concentrated FvML distributions include \cite{Wat1983b} and \cite{MarJup2000}.

This review of the literature shows that inference on~$\thetab$ under high concentration is a classical topic in directional statistics. Yet this review also reveals some important limitations in previous studies: (i) all asymptotic results available are as~$\kappa\to\infty$ with~$n$ fixed, while, parallel to the low-concentration case above, a double asymptotic scenario where~$\kappa=\kappa_n$ would go to infinity with~$n$ would be at least as natural (particularly so if~$\kappa_n$ would be allowed to diverge to infinity at an arbitrary rate as a function of~$n$); 
\linebreak
(ii) all results are limited to the parametric case of FvML distributions, so that the asymptotic properties of the spherical mean and of the Watson test remain unknown in the broader semiparametric class of~${\rm Rot}_p(\thetab,\kappa,f)$ distributions; (iii) for hypothesis testing, most works focused on the null hypothesis: very few results try and describe asymptotic powers under sequences of local alternatives, and, more importantly, not a single optimality result, to the best of our knowledge, was obtained in the literature. In this paper, we therefore fill an important gap by deriving results that are getting rid of the limitations~(i)--(iii).  

The outline of the paper is as follows. In Section~\ref{sec:notationassumption}, we fix the notation, introduce the assumptions that will be used throughout and characterize the rotationally symmetric distributions that provide high concentration for arbitrarily large values of~$\kappa$. In Section~\ref{sec:estimation}, we derive the asymptotic distribution of the spherical mean in a double asymptotic scenario where~$\kappa_n$ diverges to infinity at an arbitrary rate with~$n$. Interestingly, in contrast with what happens for low concentrations, the consistency rate here depends on the nuisance function~$f$. We also provide confidence zones for~$\thetab$ that quite naturally take the form of spherical caps centered at the spherical mean. In Section~\ref{sec:testing}, we study the asymptotic behaviour of the Watson and Wald tests. In Section~\ref{sec:LAN}, we turn to optimality issues and show that, under mild assumptions on~$f$, the sequence of statistical experiments considered is locally asymptotically normal. We establish the Le Cam optimality of the spherical mean estimator and of the Watson and Wald tests under high concentration. 
Finally, a real data application is conducted in Section \ref{sec:realdata} and a wrap up is provided in Section~\ref{sec:wrapup}. Proofs are collected in the appendix.


\section{High concentration} 
\label{sec:notationassumption}

Throughout, we will denote as~${\rm P}\n_{\thetab_n,\kappa_n,f}$ the hypothesis under which the observations $\Xb_{n1},\ldots,\Xb_{nn}$ form a random sample from the distribution~${\rm Rot}_p(\thetab_n,\kappa_n,f)$ described in the introduction, that is, the hypothesis under which these observations are mutually independent and share the common density
\begin{equation}
\label{rotsympdf}
	\xb \mapsto \frac{c_{p,\kappa_n,f}\Gamma(\frac{p-1}{2})}{2\pi^{(p-1)/2}} f(\kappa_n \xb'\thetab_n)
,
\end{equation}
where $\Gamma(\cdot)$ is the Euler Gamma function and the constant~$c_{p,\kappa,f}$ is given by 
\begin{equation}
\label{valuec}	
c_{p,\kappa,f}
:=
1\, \Big/ \int_{-1}^1 (1-s^2)^{(p-3)/2} f(\kappa s)\,ds
.
\end{equation} 
In the sequel,~$f:\R\to\R^+$ is assumed to be monotone non-decreasing on~$(-\infty,0]$ and monotone increasing on~$[0,\infty)$. Under this assumption, the location parameter~$\thetab_n$ is properly identified as the modal location on the sphere. One way to also make~$\kappa_n$ and~$f$ identifiable would be to further impose~$f(0)=f'(0)=1$. We will not impose these conditions since we also want to consider functions~$f$ that are not differentiable at zero. The resulting lack of identifiability will not be an issue in the sequel since~$\kappa_n$ and~$f$ play the role of nuisance parameters when conducting inference on~$\thetab_n$. 

We will often make use of the \emph{tangent-normal decomposition of~$\Xb_{ni}$} with respect to~$\thetab_n$, which reads~$\Xb_{ni}=u_{ni} \thetab_n+v_{ni}\Sb_{ni}$, with
$$
u_{ni}
=
\Xb_{ni}'\thetab_n
,
\quad
v_{ni}:=\sqrt{1-u_{ni}^2}
,
$$
and 
$$
\Sb_{ni}
:=
\frac{(\mathbf{I}_p-\thetab_n\thetab_n')\Xb_{ni}}{\|(\mathbf{I}_p-\thetab_n\thetab_n')\Xb_{ni}\|}
=
\frac{1}{v_{ni}}
(\mathbf{I}_p-\thetab_n\thetab_n')\Xb_{ni}
.
$$ 
The cosine~$u_{ni}$ is associated with the latitude of~$\Xb_{ni}$ with respect to the ``north pole"~$\thetab_n$, whereas~$\Sb_{ni}$ determines the corresponding hyper-longitude. Under~${\rm P}\n_{\thetab_n,\kappa_n,f}$,~$u_{n1}$ and~$\Sb_{n1}$ are mutually independent, $\Sb_{n1}$ is uniformly distributed on~$\mathcal{S}^\perp_{\thetab_n}:=\{\xb\in\mathcal{S}^{p-1}:\xb'\thetab_n=0\}$, and~$u_{n1}$ admits the density
\begin{equation} 
	\label{pdfu}
s
\mapsto 
c_{p,\kappa_n,f} (1-s^2)^{(p-3)/2} f(\kappa_n s)\,\mathbb{I}[s\in[-1,1]]
,
\end{equation}
where~$\mathbb{I}[A]$ stands for the indicator function of the set~$A$. The moments of~$u_{n1}$ under~${\rm P}\n_{\thetab_n,\kappa_n,f}$ will play an important role below and will be denoted as~$e_{n\ell}:={\rm E}[u_{n1}^\ell]$, $\ell=1,2,\ldots$ We will also write~$\tilde{e}_{n2}=e_{n2}-e_{n1}^2$ for the corresponding variance. The function~$f$ governs (jointly with~$\kappa_n$) the distribution of the angle~$\arccos (u_{n1})$ between~$\Xb_{n1}$ and~$\thetab_n$, hence is sometimes referred to as an \emph{angular function}.

The present paper is concerned with sequences of rotationally symmetric distributions that are asymptotically highly concentrated, meaning that the probability mass of any fixed spherical cap centered at~$\thetab_n$ converges to one as~$n$ diverges to infinity. More precisely, we will say that the sequence of hypotheses~${\rm P}\n_{\thetab_n,\kappa_n,f}$ is asymptotically highly concentrated if and only if, for any sequence~($\kappa_n$) diverging to infinity and any~$\varepsilon\in(0,2)$, we have
\begin{equation}
\label{HCdef}	
{\rm P}\n_{\thetab_n,\kappa_n,f}\big[\Xb_{n1}'\thetab_n>1-\varepsilon\big]
=
c_{p,\kappa_n,f} \int_{1-\varepsilon}^1 (1-s^2)^{(p-3)/2} f(\kappa_n s)\,ds \to 1 
,
\end{equation}   
that is, if and only if~$u_{n1}$ converges in probability to one as soon as~$(\kappa_n)$ diverges to infinity. Since this is clearly a property that depends on~$f$ only, we will say that~$f$ provides high concentration if and only if~(\ref{HCdef}) holds. Not all functions~$f$ provide high concentration. The polynomial functions~$z\mapsto f(z)=t^b\mathbb{I}[t\geq 0]$ are examples since, for any~$\varepsilon\in(0,1)$, they yield
$$
c_{p,\kappa_n,f} \int_{1-\varepsilon}^1 (1-s^2)^{(p-3)/2} f(\kappa_n s)\,ds
=
\frac{\int_{1-\varepsilon}^1 (1-s^2)^{(p-3)/2} s^b\, ds}{\int_{0}^1 (1-s^2)^{(p-3)/2} s^b\, ds}
=:
C
< 
1 
,
$$
where~$C$ does not depend on~$n$. It is easy to check that~$z\mapsto f(z)=\frac{\pi}{2}+\arctan(z)$ does not provide high concentration either, but that the  angular FvML function~$z\mapsto f(z)=\exp(z)$ does. It is therefore desirable to characterize the functions~$f$ providing high concentration, which is the aim of the following result. 

\begin{Theor}
\label{Theordumal}
Let~$f:\R\to\R^+$ be monotone non-decreasing on~$(-\infty,0]$ and monotone increasing on~$[0,\infty)$. Assume that~$f$ is differentiable in a neighborhood of~$\infty$ $($in the sense that there exists~$M$ such that~$f$ is differentiable over~$(M,\infty))$ and put~$\varphi_f:=f'/f$, where~$f'$ is the derivative of~$f$. Then we have the following:
\begin{itemize}
\item[(i)] If~$\kappa \varphi_f(\kappa)\nearrow\infty$ as~$\kappa\to\infty$, then~$f$ provides high concentration. 
\item[(ii)] If~$\kappa \varphi_f(\kappa)\to c(>0)$ as~$\kappa\to\infty$, then~$f$ does not provide high concentration. 
\item[(iii)] If~$\kappa \varphi_f(\kappa)\searrow 0$ as~$\kappa\to\infty$, then~$f$ does not provide high concentration. 
\end{itemize}
\end{Theor}

In this result,~$g(\kappa)\nearrow\infty$ (resp., $g(\kappa)\searrow 0$) as~$\kappa\to\infty$ means that (a)~$g(\kappa)$ diverges to infinity (resp., converges to zero) as~$\kappa$ diverges to infinity and that (b) there exists~$M$ such that~$\kappa\mapsto g(\kappa)$ is monotone non-decreasing (resp., monotone non-increasing) over~$(M,\infty)$. Essentially, Theorem~\ref{Theordumal} states that high concentration is obtained if~$f(z)$ diverges to infinity at least exponentially fast as~$z$ diverges to infinity. In particular, this result confirms that the polynomial and arctan functions~$f$ above do not provide high concentration, but that the FvML one does. Writing throughout~$z^b:={\rm sgn}(z)|z|^b$, it also shows that all functions~$z\mapsto f_b(z):=\exp(z^b)$, with $b>0$, do provide high concentration. These functions~$f$, which include the FvML one, will be our main running examples below. 

In the rest of the paper, $\mathcal{F}$ will stand for the collection of functions~$f:\R\to\R^+$ that (i) are monotone non-decreasing on~$(-\infty,0]$ and monotone increasing on~$[0,\infty)$, (ii) are differentiable in a neighborhood of~$\infty$, (iii) are such that~$\kappa \varphi_f(\kappa)\nearrow\infty$ as~$\kappa\to\infty$ and (iv) satisfy, for any~$\xi,\zeta>-1$,
\begin{equation}
	\label{conds}
\int_{-1}^1
g_{\xi, \zeta}(s)
\big|
e^{\log f(\kappa s)- \log f(\kappa)}
-
e^{(s-1) \kappa\varphi_f(\kappa)}
\big|
 \,ds
=
o\bigg( \frac{1}{(\kappa \varphi_f(\kappa))^{\xi+1}} \bigg)
\end{equation}
as~$\kappa\to\infty$, with $g_{\xi, \zeta}(s):=(1-s)^\xi
(1+s)^\zeta$. As the following result shows, our prototypical examples of angular functions~$f$ providing high concentration meet these properties.

\begin{Prop}
\label{PropCondexpb}
For any~$b>0$, the function~$z\mapsto f_b(z)=\exp(z^b)$ belongs to~$\mathcal{F}$. 
\end{Prop}

As already mentioned, the moments of~$u_{n1}=\Xb_{n1}'\thetab_n$ under~${\rm P}\n_{\thetab_n,\kappa_n,f}$ will play a key role in the sequel. It will actually be important to understand the  asymptotic behavior of these moments under high concentration. This is the role of the following result.

\begin{Theor}
\label{Theorexpansions}
Fix an integer~$p\geq 2$ and~$f\in\mathcal{F}$. Let~$(\kappa_n)$ be a positive real sequence that diverges to infinity. Then, 
%
$$
\hspace{-8mm}  
(i)
\qquad
1-e_{n2}
=
\frac{p-1}{\kappa_n \varphi_f(\kappa_n)}
+
o\bigg(\frac{1}{\kappa_n \varphi_f(\kappa_n)}\bigg)
,
$$
%
$$
\hspace{-5mm} 
(ii)
\qquad
\tilde{e}_{n2}
=
\frac{p-1}{2(\kappa_n \varphi_f(\kappa_n))^2}
+
o\bigg(\frac{1}{(\kappa_n \varphi_f(\kappa_n))^2}\bigg)
$$
and
%
$$
\hspace{-1mm} 
(iii)
\qquad
{\rm E}\big[ v_{n1}^4\big]
=
\frac{p^2-1}{(\kappa_n \varphi_f(\kappa_n))^2}
+
o\bigg(\frac{1}{(\kappa_n \varphi_f(\kappa_n))^2}\bigg)
$$
as~$n\to\infty$. 
\end{Theor}
\vspace{1mm}

As a corollary, we have
\begin{equation}
	\label{sjffl}
\frac{(1-e_{n2})^2}{\tilde{e}_{n2}}
=
2(p-1)
+
o(1)
\quad
\textrm{ and }
\quad
\frac{{\rm E}\big[ v_{n1}^4\big]}{\tilde{e}_{n2}}
=
2(p+1)+o(1)
\end{equation}
as~$n\to\infty$. Also, Vitali's Theorem (see, e.g., Theorem~5.5 in \citealp{Sho2000}) readily implies that, under the conditions of Theorem~\ref{Theorexpansions}, $e_{n1}=1+o(1)$ as~$n\to\infty$. One could obtain an expansion of~$1-e_{n1}$ that is similar to the one in Theorem~\ref{Theorexpansions}(i), but we will not do so since this is not needed for our purposes.


\section{Point estimation} 
\label{sec:estimation}

As mentioned in the introduction, the most classical estimator of location under rotational symmetry is the spherical mean, which is given by
$$
\hat{\thetab}_n
:=
\frac{\bar{\Xb}_n}{\|\bar{\Xb}_n\|}
,
$$
with~$\bar{\Xb}_n:=\frac{1}{n} \sum_{i=1}^n \Xb_{ni}$. Under~${\rm P}\n_{\thetab,\kappa_n,f}$, ${\rm E}[\Xb_{n1}]=\lambda_{\kappa_n,f}\thetab$ for some positive scalar factor~$\lambda_{\kappa_n,f}$, so that the spherical mean is a moment-type estimator of~$\thetab$. It is easy to check that it is also the maximum likelihood estimator of~$\thetab$ in the class of FvML distributions. This makes it desirable to investigate the asymptotic behavior of this estimator under high concentration. We have the following result.

\begin{Theor}
\label{Theorsphericalmean}
Fix an integer~$p\geq 2$, $\thetab\in\mathcal{S}^{p-1}$ and~$f\in\mathcal{F}$. Let~$(\kappa_n)$ be a positive real sequence that diverges to infinity. Then, under~${\rm P}\n_{\thetab,\kappa_n,f}$, 
\begin{equation}
	\label{consist}
\sqrt{n\kappa_n \varphi_f(\kappa_n)} \, (\hat{\thetab}_n - \thetab)
\stackrel{\mathcal{D}}{\to}
\mathcal{N}
\big(
{\bf 0}
,
{\bf I}_p- \thetab \thetab\pr 
\big)
\end{equation}
as~$n\to\infty$, so that, still under~${\rm P}\n_{\thetab,\kappa_n,f}$, 
\begin{equation}
	\label{ICpre}
n\kappa_n\varphi_f(\kappa_n)
\big(1-(\thetab\pr\hat{\thetab}_n)^2\big) 
\stackrel{\mathcal{D}}{\to}
\chi^2_{p-1}
\end{equation}
as~$n\to\infty$ (throughout, $\stackrel{\mathcal{D}}{\to}$ denotes convergence in distribution).
\end{Theor}
\vspace{3mm}

Since the sequence~$(\kappa_n\varphi_f(\kappa_n))$ diverges to infinity under high concentration, Theorem~\ref{Theorsphericalmean} shows that the consistency rate of the spherical mean is faster than the usual parametric root-$n$ rate. Interestingly, this consistency rate depends on the angular function~$f$. For instance, for~$f(z)=\exp(z^b)$ with~$b>0$, the rate is~$n^{(b+1)/2}$, hence can be arbitrary close to the standard root-$n$ rate for small~$b$, but can also provide arbitrary fast polynomial convergence. Clearly, even faster rates can be achieved by considering more extreme high concentration patterns. 

The asymptotic result~(\ref{ICpre}) in principle allows constructing confidence zones for~$\thetab$. More precisely, it follows from this result that a confidence zone for~$\thetab$ at asymptotic confidence level~$1-\alpha$ is given by
$$
\Bigg\{
\thetab\in\mathcal{S}^{p-1}\!: |\thetab\pr\hat{\thetab}_n| 
\geq
\sqrt{1-\frac{\chi^2_{p-1,1-\alpha}}{n\kappa_n\varphi_f(\kappa_n)}}
\
\Bigg\}
,
$$
where~$\chi^2_{p-1,1-\alpha}$ denotes the upper $\alpha$-quantile of the~$\chi^2_{p-1}$ distribution. This confidence zone, however, is problematic in two respects. First, it is not connected, as it takes the form of two antipodal spherical caps centered at~$\pm\hat\thetab_n$, which is not natural. Second, while the $f$-dependent consistency rate in Theorem~\ref{Theorsphericalmean} is interesting, it also leads to confidence zones that cannot be used in practice since~$f$ is usually an unspecified nuisance. The first problem can be dealt with by deriving a weak limiting result for~$\thetab'\hat\thetab_n$ obtained from a second-order delta method (while Theorem~\ref{Theorsphericalmean} results from a classical, first-order, delta method). We have the following result.

\begin{Theor}
\label{TheorsphericalCap}
Fix an integer~$p\geq 2$, $\thetab\in\mathcal{S}^{p-1}$ and~$f\in\mathcal{F}$. Let~$(\kappa_n)$ be a positive real sequence that diverges to infinity. Then, under~${\rm P}\n_{\thetab,\kappa_n,f}$, 
$2n\kappa_n\varphi_f(\kappa_n)
(1-\thetab\pr\hat{\thetab}_n) 
\stackrel{\mathcal{D}}{\to}
\chi^2_{p-1}$
as~$n\to\infty$.
\end{Theor}

This second-order result provides confidence zones at asymptotic confidence level~$1-\alpha$ that are given by 
\begin{equation}
\label{sphercapICpresque}	
\Bigg\{
\thetab\in\mathcal{S}^{p-1}\!: \thetab\pr\hat{\thetab}_n 
\geq
1
- 
\frac{\chi^2_{p-1,1-\alpha}}{2n\kappa_n\varphi_f(\kappa_n)}
\Bigg\}
,
\end{equation}
hence take, quite naturally, the form of (connected) spherical caps centered at~$\hat\thetab_n$. Of course, these confidence zones still cannot be used in practice since~$f$ is unspecified. Fortunately, Theorem~\ref{Theorexpansions}(i) allows replacing the unknown quantity~$\kappa_n\varphi_f(\kappa_n)$ by the quantity~$(p-1)/(1-e_{n2})=(p-1)/(1-{\rm E}[(\Xb_{n1}\pr \thetab)^2])$, which can be naturally estimated by~$(p-1)/(1-\hat{e}_{n2})$, where we let~$\hat{e}_{n2}:=\frac{1}{n} \sum_{i=1}^n(\Xb_{ni}\pr \hat{\thetab}_n)^2$. The following result, that guarantees that this replacement has no asymptotic impact, opens the door to the construction of feasible confidence zones. 

\begin{Theor}
\label{Theorsphericalmean2}
Fix an integer~$p\geq 2$, $\thetab\in\mathcal{S}^{p-1}$ and~$f\in\mathcal{F}$. Let~$(\kappa_n)$ be a positive real sequence that diverges to infinity. Then, under~${\rm P}\n_{\thetab,\kappa_n,f}$, 
$$
\frac{\sqrt{n(p-1)}(\hat{\thetab}_n- \thetab)}{\sqrt{1-\hat{e}_{n2}}}
\stackrel{\mathcal{D}}{\to}
\mathcal{N}
\big(
{\bf 0}
,
{\bf I}_p- \thetab \thetab\pr 
\big)
$$
as~$n\to\infty$, and, still under~${\rm P}\n_{\thetab,\kappa_n,f}$, 
$$
\frac{n(p-1)\big(1-(\thetab\pr\hat{\thetab}_n)^2\big)}{1-\hat{e}_{n2}}  
\stackrel{\mathcal{D}}{\to}
\chi^2_{p-1}
\quad \textrm{ and } \quad
\frac{2n(p-1)(1-\thetab\pr\hat{\thetab}_n)}{1-\hat{e}_{n2}}  
\stackrel{\mathcal{D}}{\to}
\chi^2_{p-1}
$$
as~$n\to\infty$, where, in all cases,~$\hat{e}_{n2}=\frac{1}{n} \sum_{i=1}^n(\Xb_{ni}\pr \hat{\thetab}_n)^2$.
\end{Theor}

\begin{figure}[htbp!] 
\begin{center}
\includegraphics[width=\textwidth]{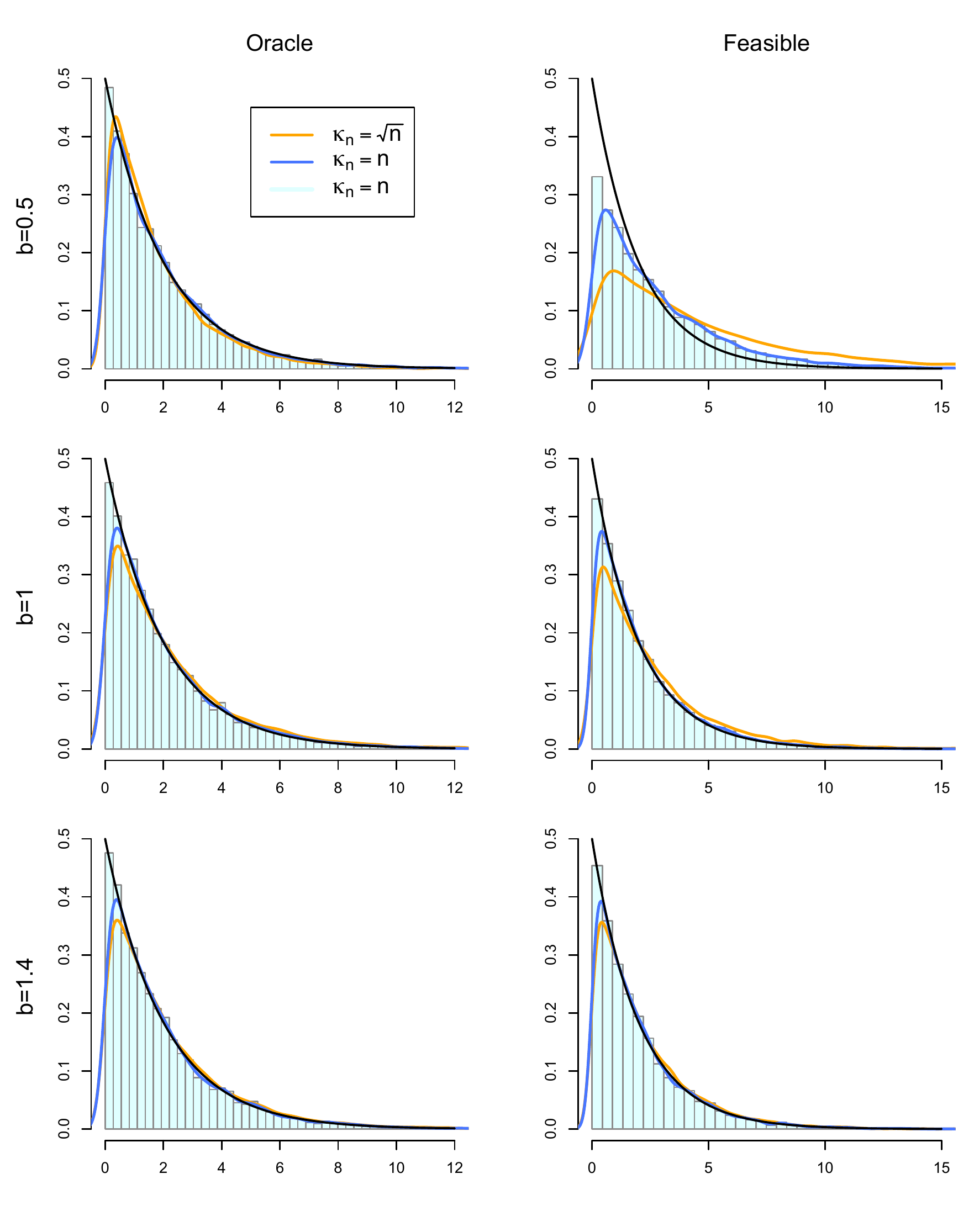} 
\end{center}
\vspace{-4mm}
\caption{Histograms of~$2n\kappa_n\varphi_f(\kappa_n)(1-\thetab\pr\hat{\thetab}_n)$ (left) and~$2n(p-1)(1-\thetab\pr\hat{\thetab}_n)/(1-\hat{e}_{n2})$ (right) computed from $10,\!000$ random samples from~${\rm P}\n_{\thetab,%
\kappa_n,f_b}$, 
\vspace{-.6mm}
with~$n=100$, $\thetab=(1,0,0)\pr$, $\kappa_n=n$, and~$f_b(z)=\exp(z^b)$, for~$b=0.5$ (top), $b=1$
 (middle) and~$b=1.4$ (bottom). The blue curve is the kernel density estimate resulting from the R command \texttt{density} with default parameter values. The orange curve is the corresponding kernel density estimate for random samples generated with~$\kappa_n=\sqrt{n}$. The theoretical limiting density, namely the density of the~$\chi^2_2$ distribution, is plotted in black.}
\label{Fig1}
\end{figure}

\begin{figure}[htbp!] 
\begin{center}
\includegraphics[width=\textwidth]{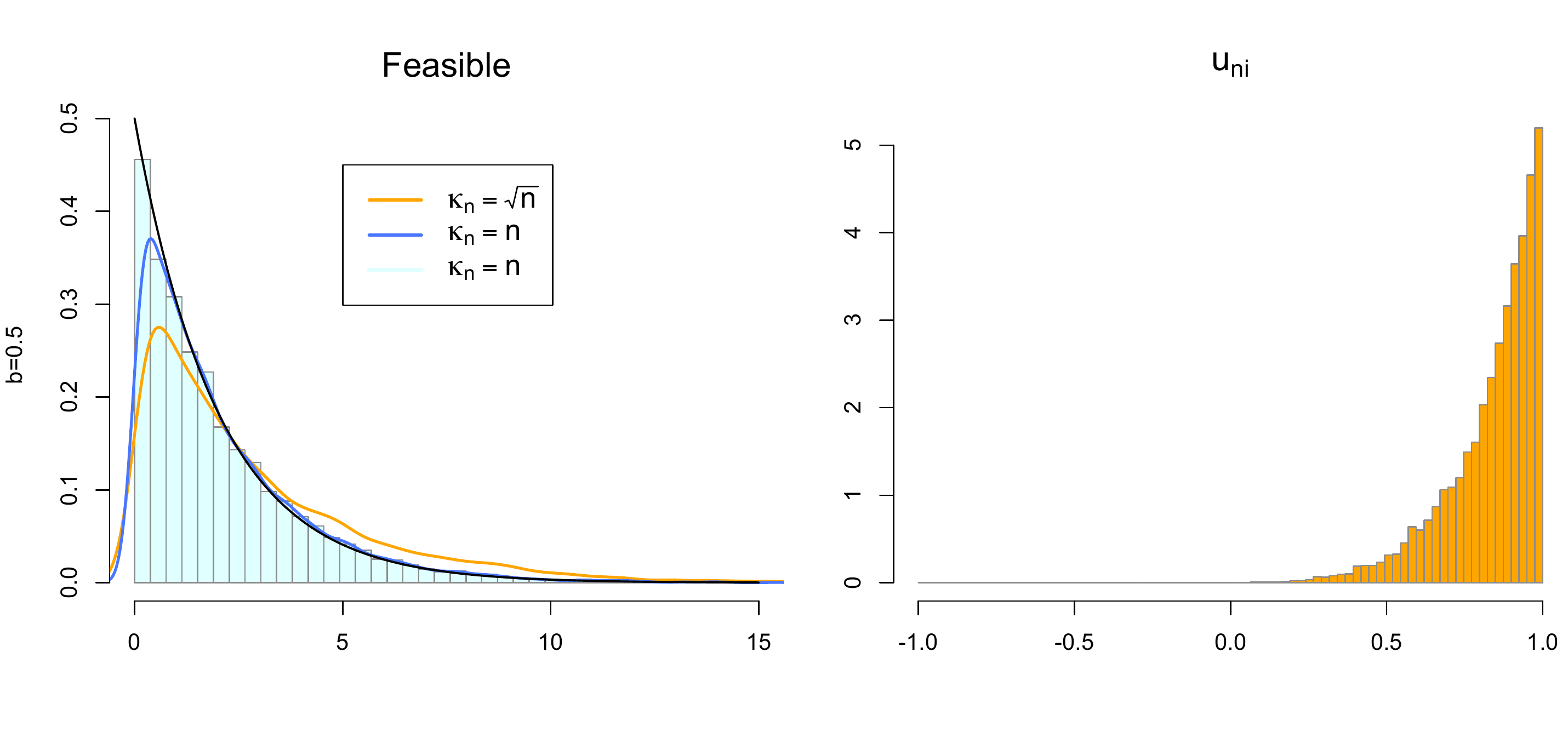}
\end{center}
\vspace{-5mm}
\caption{(Left:) Histogram 
\vspace{-.7mm}
of~$2n(p-1)(1-\thetab\pr\hat{\thetab}_n)/(1-\hat{e}_{n2})$ computed from $2,\!500$ random samples from~${\rm P}\n_{\thetab,\kappa_n,f_b}$, with~$n=10,\!000$, $\thetab=(1,0,0)\pr$, $\kappa_n=n$, and~$f_b(z)=\exp(z^b)$, for~$b=0.5$. The blue curve is the kernel density estimate resulting from the R command \texttt{density} with default parameter values. The orange curve is the corresponding kernel density estimate for random samples generated with~$\kappa_n=\sqrt{n}$. The theoretical limiting density is still plotted in black. 
\vspace{-.8mm}
(Right:) Histogram of~$u_{ni}=\Xb_{ni}\pr\thetab$, $i=1,\ldots,n$, where the~$\Xb_{ni}$'s form a random sample from~${\rm P}\n_{\thetab,\kappa_n,f_b}$, 
\vspace{-.8mm}
with~$n=10,\!000$, $\thetab=(1,0,0)\pr$, $\kappa_n=\sqrt{n}$, and~$f_b(z)=\exp(z^b)$, for~$b=0.5$.}
\label{Fig2}
\end{figure}

As a direct corollary, a feasible version of the spherical cap confidence zone in~(\ref{sphercapICpresque}) is  
\begin{equation}
\label{seccapbis}
\Bigg\{
\thetab\in\mathcal{S}^{p-1}\!: \thetab\pr\hat{\thetab}_n 
\geq
1 
- 
\frac{1-\hat{e}_{n2}}{2n(p-1)}
\chi^2_{p-1,1-\alpha} 
\Bigg\} 
.
\end{equation}
We conducted the following Monte Carlo exercises to check the validity of Theorems~\ref{TheorsphericalCap}--\ref{Theorsphericalmean2}. For each combination of~$a\in\{0.5,1\}$ and~$b\in\{0.5,1,1.4\}$, we generated $M=10,\!000$ random samples of size~$n=100$ from the rotationally symmetric distribution with location~$\thetab=(1,0,0)'\in\mathcal{S}^2$, concentration~$\kappa_n=n^a$, and angular function~$z\mapsto f_b(z)=\exp(z^b)$ (numerical overflows prevented us from considering larger values of~$b$). For each~$a$ and~$b$, Figure~\ref{Fig1} plots kernel density estimates of the resulting $M$ values of~$T_n^{\rm Oracle}:=2n\kappa_n\varphi_f(\kappa_n)(1-\thetab\pr\hat{\thetab}_n)$ and~$T^{\rm Feasible}_n:=2n(p-1)(1-\thetab\pr\hat{\thetab}_n)/(1-\hat{e}_{n2})$ (for~$a=1$, raw histograms are also provided). Clearly, Figure~\ref{Fig1} supports the theoretical results above, with possibly one exception only, namely the case of~$T^{\rm Feasible}_n$ with~$b=0.5$. We therefore focused on this case and repeated the same Monte Carlo exercise with~$n=10,\!000$. The results, that are shown in Figure~\ref{Fig2}, are now in perfect agreement with the theory for~$a=1$, whereas the fit still is not excellent for~$a=0.5$. A closer inspection provides the explanation: despite the large sample size~$n$ considered in Figure~\ref{Fig2}, the distribution associated with~$a=b=0.5$ is far for being highly concentrated; see the right panel of this figure. The fit observed for~$a=0.5$ in the left panel of  Figure~\ref{Fig2} therefore does not contradict our theoretical results, which would materialize for higher concentrations.


\section{Hypothesis testing} 
\label{sec:testing}

We now turn to hypothesis testing and, more specifically, to the generic problem of testing the null hypothesis~$\mathcal{H}_0:\thetab=\thetab_0$ against the alternative~$\mathcal{H}_1:\thetab\neq \thetab_0$, where~$\thetab_0$ is a fixed unit $p$-vector. In this section, we consider the Watson test (\citealp[p.\!~140]{Wat1983b}) and the Wald test (\citealp{Hay1990,HayPur1985}), that respectively reject the null hypothesis at asymptotic level~$\alpha$ whenever 
\begin{equation}
\label{WatsonTestStatistics}	
W_n 
:=
\frac{n (p-1)\bar{\Xb}_n\pr ({\bf I}_p- \thetab_{0} \thetab_{0}\pr) \bar{\Xb}_n}{1-\frac{1}{n} \sum_{i=1}^n(\Xb_{ni}\pr\thetab_{0})^2} 
\end{equation}
and
\begin{equation}
\label{WaldTestStatistics}
S_n
=
\frac{n(p-1) (\bar{\Xb}_{n}\pr\thetab_{0})^2
\,
 \hat{\thetab}_n\pr({\bf I}_p- \thetab_{0}\thetab_{0}\pr)\hat{\thetab}_n}{1-\frac{1}{n} \sum_{i=1}^n(\Xb_{ni}\pr\thetab_{0})^2}
\end{equation}
exceed the critical value~$\chi^2_{p-1,1-\alpha}$. In standard asymptotic scenarios where the sample size~$n$ diverges to infinity with~$\kappa$ fixed, the Watson and Wald test statistics are asymptotically equivalent in probability under the null hypothesis, hence also under sequences of contiguous alternatives, so that these tests may be considered asymptotically equivalent. As shown in \cite{PaiVer17b}, however, this asymptotic equivalence does not survive asymptotic scenarios for which~$\kappa_n=O(1/\sqrt{n})$ as~$n$ diverges to infinity. This suggests investigating the asymptotic behavior of these tests under the high concentration scenarios considered in the previous sections. 
 
To do so, let
$$
R_n
=
\frac{1-\frac{1}{n} \sum_{i=1}^n(\Xb_{ni}\pr\thetab_{0})^2}{\sqrt{2(p-1)} \tilde{e}_{n2}^{1/2}}
$$
and decompose the Watson and Wald test statistics into
$$
W_n 
=: 
\frac{\tilde{W}_n}{R_n}
\quad\textrm{and}\quad
S_n
=:
\frac{(\bar{\Xb}_{n}\pr\thetab_{0})^2 \tilde{S}_n}{R_n}
\cdot
$$
We then have the following lemma.

\begin{Lem}
\label{Lemtest}
Fix an integer~$p\geq 2$, $\thetab_0\in\mathcal{S}^{p-1}$ and~$f\in\mathcal{F}$. Let~$(\kappa_n)$ be a positive real sequence that diverges to infinity. Let~$(\taub_n)$ be a bounded sequence in~$\R^p$ such that~$\thetab_n=\thetab_0+\nu_n\taub_n\in\mathcal{S}^{p-1}$ for all~$n$, with~$\nu_n:=1/\sqrt{n\kappa_n\varphi_f(\kappa_n)}$. Then, under~${\rm P}\n_{\thetab_n,\kappa_n,f}$, we have~$R_n=1+o_{\rm P}(1)$ and~$\bar{\Xb}_{n}\pr\thetab_{0}=1+o_{\rm P}(1)$ as~$n\to\infty$, so that~$W_n=\tilde{W}_n+o_{\rm P}(1)$ and~$S_n=\tilde{S}_n+o_{\rm P}(1)$ as~$n\to\infty$.  	
\end{Lem}

This lemma ensures that, both under the sequence of null hypotheses \linebreak ${\rm P}\n_{\thetab_0,\kappa_n,f}$ (taking $\taub_n \equiv {\bf 0}$) and under sequences of local alternatives of the form~${\rm P}\n_{\thetab_n,\kappa_n,f}$, one may focus on~$\tilde{W}_n$ and~$\tilde{S}_n$ when studying the asymptotic behaviors of the Watson and Wald test statistics in~(\ref{WatsonTestStatistics})--\eqref{WaldTestStatistics}. These asymptotic behaviors are provided in the following result.

\begin{Theor}
\label{Theortest}
Fix an integer~$p\geq 2$, $\thetab_0\in\mathcal{S}^{p-1}$ and~$f\in\mathcal{F}$. Let~$(\kappa_n)$ be a positive real sequence that diverges to infinity. Let~$(\taub_n)$ be a sequence in~$\R^p$ converging to~$\taub$ and such that~$\thetab_n=\thetab_0+\nu_n\taub_n\in\mathcal{S}^{p-1}$ for all~$n$, with~$\nu_n:=1/\sqrt{n\kappa_n\varphi_f(\kappa_n)}$. Then, (i) under~${\rm P}\n_{\thetab_0,\kappa_n,f}$, 
$$
W_n=S_n+o_{\rm P}(1)
\stackrel{\mathcal{D}}{\to}
\chi^2_{p-1}
$$
as~$n\to\infty$; (ii) under~${\rm P}\n_{\thetab_n,\kappa_n,f}$, 
$$
W_n=S_n+o_{\rm P}(1)
\stackrel{\mathcal{D}}{\to}
\chi^2_{p-1}
\big(
\|\taub\|^2
\big)
$$
as~$n\to\infty$, where~$\chi^2_{p-1}(c)$ denotes the non-central chi-square distribution with~$p-1$ degrees of freedom and non-centrality parameter~$c$. 
\end{Theor}

This result shows that, under high concentration, the Watson and Wald test statistics remain asymptotically equivalent in probability both under the null hypothesis and under the considered sequences of local alternatives. Both tests show asymptotic size~$\alpha$ under the null hypothesis, irrespective of the angular function~$f$ and of the rate at which the concentration~$\kappa_n$ diverges to infinity. Theorem~\ref{Theortest} also reveals that~$\nu_n$ describes the consistency rate of these tests, in the sense that the Watson and Wald tests show non-trivial asymptotic powers (that is, asymptotic powers in~$(\alpha,1)$) under sequences of local alternatives of the form~${\rm P}\n_{\thetab_n,\kappa_n,f}$, with~$\nu_n^{-1}\|\thetab_n-\thetab_0\|\to c>0$. Like in point estimation, this rate depends on~$f$ and is faster than the standard parametric root-$n$ rate that is obtained for fixed~$\kappa$; that is, compared to the alternatives that can be detected in the standard fixed-$\kappa$ situation, less severe---hence, more challenging---alternatives can be detected under high concentration. 

\begin{figure}[htbp!]  
\begin{center}
\includegraphics[width=\textwidth]{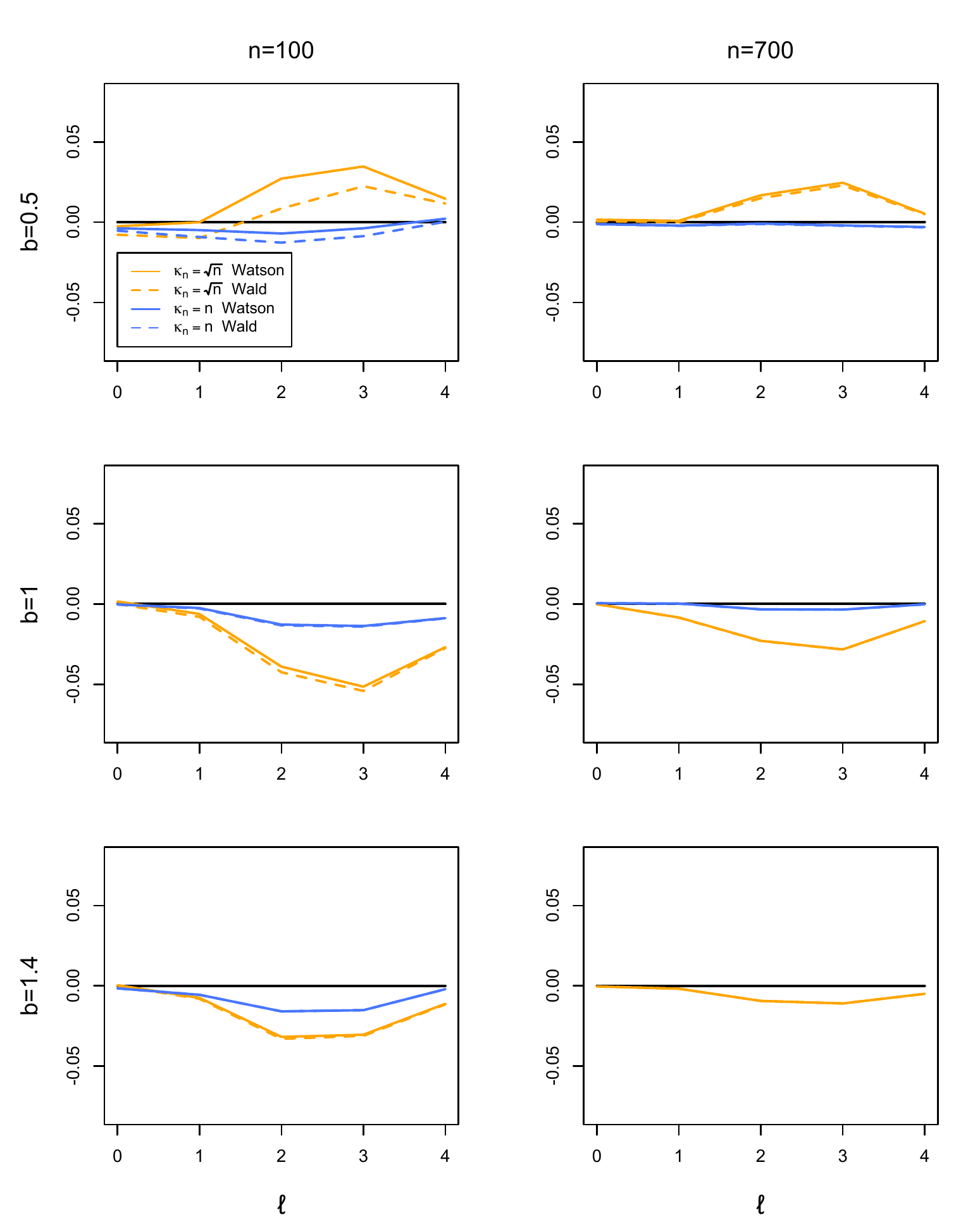}
\end{center}  
\vspace{-2mm}
\caption{Plots of the differences between (i) the rejection frequencies of the Watson and Wald tests and (ii) their theoretical limiting powers in~(\ref{theorpowersimu}). Rejection frequencies are obtained from a collection of~$M=10,\!000$ random samples of size~$n=100$ (left) or~$n=700$ (right) from the rotationally symmetric distribution with location~$\thetab_{n\ell}$ in~(\ref{alternativesimu}), concentration~$\kappa_n=\sqrt{n}$ or~$n$, and angular function~$z\mapsto f_b=\exp(z^b)$, with~$b=0.5$ (top), $1$ (middle) or $1.4$ (bottom); see Section~\ref{sec:testing} for details.}    
\label{Fig3} 
\end{figure}

We performed the following Monte Carlo exercise to illustrate the results in   Theorem~\ref{Theortest}. For each combination of~$a\in\{0.5,1\}$,~$b\in\{0.5,1,1.4\}$ and~$\ell\in\{0,1,2,3,4\}$, we generated $M=10,\!000$ random samples of size~$n=100$ from the rotationally symmetric distribution with concentration~$\kappa_n=n^a$, angular function~$z\mapsto f_b(z)=\exp(z^b)$, and location
\begin{equation}
	\label{alternativesimu}
\thetab_{n\ell}
=
\Bigg(
\!
\begin{array}{ccc}
\cos \alpha_{n\ell} & -\sin \alpha_{n\ell} & 0\\[-1.5mm]
\sin \alpha_{n\ell} & \hspace{3mm} \cos \alpha_{n\ell}  & 0\\[-1.5mm]
0  & 0 & 1
\end{array}
\!
\Bigg)
\thetab_0
,
\end{equation}
where we let~$\thetab_0=(1,0,0)'$ and~$\alpha_{n\ell} := 2\arcsin(\ell/(2\nu_n))$, with 
$\nu_n=
\linebreak 
1/\sqrt{n\kappa_n\varphi_{f_b}(\kappa_n)}$. The alternative locations~$\thetab_{n\ell}$ rewrite~$\thetab_0+\nu_n\taub_{n\ell}$ for some $p$-vector~$\taub_{n\ell}$ with norm~$\ell$. Clearly, $\ell=0$ refers to the null hypothesis~$\mathcal{H}_0:\thetab=\thetab_0$ and~$\ell=1,2,3,4$ correspond to increasingly severe alternatives. In each sample, we performed the Watson and Wald tests at asymptotic level~$\alpha=5\%$. Figure~\ref{Fig3} plots, as a function of~$\ell$, the resulting rejection frequencies, or more precisely, the difference between these rejection frequencies and the corresponding theoretical limiting powers
\begin{equation}
	\label{theorpowersimu}
{\rm P}[Y_\ell>\chi^2_{p-1,1-\alpha}],
\quad
\textrm{ with } Y_\ell\sim\chi^2_{p-1}
\big(
\ell^2
\big);
\end{equation}
see Theorem~\ref{Theortest}(ii). The figure also reports the results for sample size~$n=700$, but for the case with highest concentration (i.e., the case~$(a,b)=(1,1.4)$) for which data generation led to numerical overflow. Rejection frequencies agree well with the limiting powers (note the scale of the vertical axes), particularly for~$\kappa_n=n$ which provides a higher concentration than~$\kappa_n=\sqrt{n}$. The agreement improves as the sample size increases. In all cases but the one with lowest concentration (i.e., the case~$(a,b)=(0.5,0.5)$), the asymptotic equivalence between the Watson and Wald tests materializes already for~$n=100$.


\section{Local asymptotic normality} 
\label{sec:LAN}

The Watson test was shown to enjoy strong optimality properties, both in the standard asymptotic scenario where the concentration~$\kappa_n$ is fixed and in the non-standard one where the concentration goes to zero; see \cite{PaiVer17b}. In the latter scenario, the Wald test, on the contrary, fails to be optimal. In this section, we investigate the optimality properties of the Watson and Wald tests and of the spherical mean estimator under high concentration. Optimality will be in the Le Cam sense, which requires studying the \emph{Local Asymptotic Normality (LAN)} of the sequence of fixed-$f$ parametric submodels at hand. 

To do so, we will need to reinforce our assumptions on~$f$. Let~$p(\geq 2)$ be an integer, $(\kappa_n)$ be a positive sequence diverging to infinity, and~$(t_n)$ be a bounded positive sequence. In the sequel, we will denote as~$\mathcal{F}_{\rm LAN}(p,\kappa_n,t_n)$ the collection of angular functions~$f\in\mathcal{F}$ such that, as~$\kappa\to\infty$, 
$$
\frac{1}{f(\kappa)}
\int_{-1}^1 \! \!
\big(\varphi_f(\kappa s)-\varphi_f(\kappa)\big)^2
(1-s^2)^{(p-3)/2}
f(\kappa s)
\, \!
ds
\! =
o
\bigg(
\frac{1}{\kappa^{(p+1)/2} (\varphi_f(\kappa))^{(p-3)/2}}
\bigg)
$$
and such that, letting~$h^\pm_n(s,w):=-{\textstyle{\frac{1}{2}}}t_n^2 \kappa_n\nu_n^2 s \pm c_n t_n \kappa_n\nu_n (1-s^2)^{1/2} w^{1/2}$, with~$\nu_n:=1/\sqrt{n\kappa_n\varphi_f(\kappa_n)}$ and~$c_n:=( 1-\frac{1}{4}\nu_n^2 t_n^2)^{1/2}$,  
\begin{eqnarray}
	\lefteqn{
\hspace{-13.5mm}
\frac{1}{f(\kappa_n)}
\int_{-1}^1
\int_{0}^1 \! \!
\big|
\log f(\kappa_n s+h^\pm_n(s,w))
-
\log f(\kappa_n s)
-
h^\pm_n(s,w) \varphi_f(\kappa_n s)
\big| f(\kappa_n s)
}
\nonumber
\\[2mm]
& & 
\hspace{3mm} 
\times (1-s^2)^{(p-3)/2} 
\,
dG_p(w)
ds
=
o\bigg(
\frac{1}{n(\kappa_n \varphi_f(\kappa_n))^{(p-1)/2}}
\bigg)
\label{assudeath}
\end{eqnarray}
as~$n\to\infty$, where, for~$p\geq 3$,~$G_p$ \label{pageGp} is the cumulative distribution function of the ${\rm Beta}(\frac{1}{2},\frac{p-2}{2})$ distribution, whereas, for~$p=2$, $G_p$ is the cumulative distribution function of the Dirac distribution in~$1$. As shown in the next result, most angular functions~$f_b$ do satisfy these extra assumptions, sometimes under an extremely mild restriction on the rate at which the sequence~$(\kappa_n)$ diverges to infinity with~$n$.

\begin{Prop}
\label{PropCondLANexpb}
Let~$p(\geq 2)$ be an integer, $(\kappa_n)$ be a positive sequence diverging to infinity, and~$(t_n)$ be a bounded positive sequence. Then, for any~$b\geq 1$, the function~$z\mapsto f_b(z)=\exp(z^b)$ belongs to~$\mathcal{F}_{\rm LAN}(p,\kappa_n,t_n)$. Provided that there exists~$\varepsilon\in(0,2)$ such that~$\kappa_n^b/(\log n)\geq (1-b)/(2-\varepsilon)$ for~$n$ large enough, the same holds for~$f_b$, with~$b\in(\frac{1}{2},1)$. 
\end{Prop}

In other words,  $f_b$, with~$b\geq 1$, belongs to~$\mathcal{F}_{\rm LAN}(p,\kappa_n,t_n)$ irrespective of the sequences~($\kappa_n$) and~$(t_n)$, whereas all angular functions~$f_b$, with~$b\in(\frac{1}{2},1)$, belong to~$\mathcal{F}_{\rm LAN}(p,\kappa_n,t_n)$ in particular when~$(\kappa_n)$ diverges to infinity at least as fast as~$(\log n)^2$, hence \mbox{e.g.} when~$\kappa_n=n^a$, with~$a>0$. We then have the following LAN result.

\begin{Theor}
\label{TheorLAN}
Fix an integer~$p\geq 2$ and $\thetab\in\mathcal{S}^{p-1}$. Let~$(\kappa_n)$ be a positive real sequence that diverges to infinity. Let~$(\taub_n)$ be a bounded sequence in~$\R^p$ such that~$\thetab_n=\thetab+\nu_n\taub_n\in\mathcal{S}^{p-1}$ for all~$n$, with~$\nu_n:=1/\sqrt{n\kappa_n\varphi_f(\kappa_n)}$. Assume that~$f$ belongs to~$\mathcal{F}_{\rm LAN}(p,\kappa_n,\|\taub_n\|)$.
Then, as~$n\to\infty$ under~${\rm P}\n_{\thetab,\kappa_n,f}$, 
$$
\Lambda_{\thetab+\nu_n\taub_n/\thetab,\kappa_n,f}
:=
\log \frac{d{\rm P}\n_{\thetab+\nu_n\taub_n,\kappa_n,f}}{d{\rm P}\n_{\thetab,\kappa_n,f}} 
=
\taub_n' 
\Deltab\n_{\thetab,f}
-
\frac{1}{2} 
\taub_n\pr 
\Gamb_{\thetab}
\taub_n
+
o_{\rm P}(1)
,
$$
where the \emph{central sequence}
$
\Deltab\n_{\thetab,f}
:=
\nu_n^{-1}
(\mathbf{I}_p-\thetab\thetab')
\bar{\Xb}_n
$,
still under~${\rm P}\n_{\thetab,\kappa_n,f}$, is asymptotically normal with mean zero and covariance matrix~$\Gamb_\thetab:=\mathbf{I}_p-\thetab\thetab'\!$.  
\end{Theor}

This result shows that the rate~$\nu_n$ identified in the previous sections is actually the contiguity rate associated with the sequence of statistical experiments at hand. Remarkably, this provides one of the few semiparametric examples (if any) where the contiguity rate depends on the fixed value of the functional nuisance~$f$. Since the contiguity rate coincides with the rate of convergence of the spherical mean (see Theorem~\ref{Theorsphericalmean}), we conclude that the spherical mean is rate-consistent. Better: since the proof of Theorem~\ref{Theorsphericalmean} establishes that 
$$
\sqrt{n\kappa_n \varphi_f(\kappa_n)} \, (\hat{\thetab}_n - \thetab)=\Deltab\n_{\thetab,f}+o_{\rm P}(1)
$$
as $\ny$ under~${\rm P}\n_{\thetab,\kappa_n,f}$,
it actually follows from Theorems~\ref{Theorsphericalmean} and~\ref{TheorLAN} that the spherical mean is an asymptotically optimal estimator in the sense of the convolution theorem; see, e.g., Theorem~8.8 from \cite{van1998}. Turning to hypothesis testing, it also follows from the LAN result above that the Watson and Wald tests from the previous section are rate-consistent, since Theorem~\ref{Theortest}(ii) indicates that these tests show non-trivial asymptotic powers under the sequence of contiguous alternatives involved in Theorem~\ref{TheorLAN}. Actually, in the present LAN framework, an application of the Le Cam third lemma confirms these asymptotic local powers. 

To show this, fix a positive real sequence~$(\kappa_n)$ that diverges to infinity and local alternatives as in Theorem~\ref{TheorLAN}. Then, under the sequence of null hypotheses~${\rm P}\n_{\thetab_0,\kappa_n,f}$, 
$$
\Tb_n^W
:=
\frac{\sqrt{n}(p-1)^{1/4}(\bar{\Xb}_n-e_{n1}\thetab_0)}{2^{1/4}{\tilde{e}_{n2}^{1/4}}}
$$
is asymptotically normal with mean zero and covariance matrix~$\Gamb_{\thetab_0}$; this follows from~(\ref{sdf}) in the proof of Theorem~\ref{Theorsphericalmean}. Now, by using  Theorem~\ref{Theorexpansions}(ii), we obtain that, under the same sequence of hypotheses,
\begin{eqnarray*}
	\lefteqn{
{\rm Cov}\big[\Tb_n^W\!,\Lambda_{\thetab_0+\nu_n\taub_n/\thetab_0,\kappa_n,f}
\big]
=
{\rm Cov}\big[\Tb_n^W,\Deltab\n_{\thetab_0,f}\big]\taub_n
+
o(1)
}
\\[2mm]
& & 
\hspace{3mm} 
=
\frac{n\sqrt{\kappa_n\varphi_f(\kappa_n)}(p-1)^{1/4}}{2^{1/4}{\tilde{e}_{n2}^{1/4}}}
{\rm E}
\big[
(\bar{\Xb}_n-e_{n1}\thetab_0)
(\bar{\Xb}_n-e_{n1}\thetab_0)'
\big]
\Gamb_{\thetab_0}
\taub_n
+
o(1)
\\[2mm]
& & 
\hspace{3mm} 
=
{\rm E}
\big[
\Tb_n^W
(\Tb_n^W)'
\big]
\Gamb_{\thetab_0}
\taub_n
+
o(1)
=
\taub_n
+o(1)
.
\end{eqnarray*}
Thus, Le Cam's third lemma entails that, under the sequence of contiguous alternatives~${\rm P}\n_{\thetab_n,\kappa_n,f}$, with~$\thetab_n\!=\thetab_0+\nu_n\taub_n$, $\nu_n\!=1/\sqrt{n\kappa_n\varphi_f(\kappa_n)}$ \mbox{and~$(\taub_n)\!\to\taub$}, $\Tb_n^W$ is asymptotically normal with mean~$\taub$ and covariance matrix~$\Gamb_{\thetab_0}$, so that, under this sequence of hypotheses, 
$
\tilde{W}_n
=
(\Tb^W_n)'\Gamb_{\thetab_0}^- \Tb^W_n
\stackrel{\mathcal{D}}{\to}
\chi^2_{p-1}
\big(
\|\taub\|^2
\big)
,
$
where~$\Ab^-$ stands for the Moore-Penrose inverse of~$\Ab$. From contiguity, we thus obtain that 
$
W_n=\tilde{W}_n+o_{\rm P}(1)\stackrel{\mathcal{D}}{\to}
\chi^2_{p-1}
\big(
\|\taub\|^2
\big)
$ 
under the  alternatives considered, which, as announced, is in agreement with Theorem~\ref{Theortest}(ii). As for the Wald test, the fact that~$S_n\stackrel{\mathcal{D}}{\to}\chi^2_{p-1}
\big(
\|\taub\|^2
\big)
$ under the same sequence of alternatives directly follows from the result for the Watson test and from the fact that the null asymptotic equivalence~$W_n=S_n+o_{\rm P}(1)$ in Theorem~\ref{Theortest}(i) extends, from contiguity, to the present contiguous alternatives.  

Beyond this, one of the main interests of the LAN result in Theorem~\ref{TheorLAN} is to pave the way to the construction of Le Cam optimal tests for the problem of testing~$\mathcal{H}_0:\thetab=\thetab_0$ versus~$\mathcal{H}_1:\thetab\neq \thetab_0$ under angular function~$f$. It directly follows from this result that, for this problem, the test rejecting the null hypothesis at asymptotic level~$\alpha$ whenever
$$
Q_n
:=
(\Deltab\n_{\thetab_0,f})' \Gamb_{\thetab_0}^- \Deltab\n_{\thetab_0,f}
> 
\chi^2_{p-1,1-\alpha}
$$
is Le Cam optimal (more precisely, locally asymptotically maximin) at asym\-ptotic level~$\alpha$. Since Theorem~\ref{Theorexpansions}(ii) ensures that, under the null hypothesis, 
\begin{eqnarray*}
\Deltab\n_{\thetab_0,f}
&\!\!\!=\!\!\!&
\sqrt{n\kappa_n\varphi_f(\kappa_n)} 
(\mathbf{I}_p-\thetab_0\thetab_0')
(\bar{\Xb}_n-e_{n1}\thetab_0)
\\[2mm]
&\!\!\!=\!\!\!&
\frac{\sqrt{n}(p-1)^{1/4}}{2^{1/4}\tilde{e}_{n2}^{1/4}}
(\mathbf{I}_p-\thetab_0\thetab_0')
(\bar{\Xb}_n-e_{n1}\thetab_0)
+
o_{\rm P}(1)
,
\end{eqnarray*}
Lemma~\ref{Lemtest} readily entails that
$
Q_n
=
\tilde{W}_n
+
o_{\rm P}(1)
=
W_n
+
o_{\rm P}(1)
$
under the null hypothesis, hence, from contiguity, also under the sequences of local alternatives above. It follows that, under the assumptions of Theorem~\ref{TheorLAN}, the Watson test is optimal in the Le Cam sense. Since the Watson test does not depend on~$f$, this optimality holds at any~$f$ meeting the assumptions of Theorem~\ref{TheorLAN}. From the asymptotic equivalence result in Theorem~\ref{Theortest}(i) and from contiguity, this extends to the Wald test. 

In the high concentration framework considered, it may be intuitively appealing to linearize the problem and apply a standard Euclidean test to the data projected onto the tangent space to~$\mathcal{S}^{p-1}$ at the null location~$\thetab_0$---or equivalently, to the data~$\Yb_{ni}:={\bf P}_{\thetab_0}'\Xb_{ni}$, $i=1,\ldots,n$, where~${\bf P}_{\thetab_0}$ is an arbitrary $p\times (p-1)$ matrix whose columns form an orthornormal basis of the orthogonal complement of~$\thetab_0$ in~$\R^p$. The null hypothesis~$\mathcal{H}_0:\thetab=\thetab_0$ translates into testing that the mean of the common (under rotational symmetry about~$\thetab_0$, spherically symmetric) distribution of the~$\Yb_{ni}$'s is the zero vector. The Watson test can actually be seen as the (spherical) Hotelling test rejecting the null hypothesis at asymptotic level~$\alpha$ whenever $n \bar{\bf Y}_n'\Sb_n^{-1}\bar{\bf Y}_n>\chi^2_{p-1,1-\alpha}$, with~$\bar{\bf Y}_n:= n^{-1} \sum_{i=1}^n \Yb_{ni}$ and with a standardization matrix~$\Sb_n$ that, in line with the underlying spherical symmetry, is a multiple of the identity matrix. Quite nicely, Theorem~\ref{TheorLAN} formally proves that this linearization provides a test that is Le Cam optimal at any~$f$. We insist, however, that it was unclear that such a linearization would provide a test that achieves optimality in the original sequence of curved statistical experiments. Not only because the impact of linearization is difficult to control, but also because it is unknown whether or not  the spherical Hotelling test is optimal in any sense under the, highly concentrated and skewed, alternatives obtained in the tangent space (to the best of our knowledge, the only optimality results for the spherical Hotelling test relate to shifted spherical Gaussian distributions; see, e.g., \citealp{HP2002}).


\section{Real data illustration}
\label{sec:realdata}

The real dataset we analyze here consists in measurements of magnetic remanence directions in $n=62$ rock specimens. The objective of Remanent magnetism or equivalently Paleomagnetism is to study the strength and the direction of the Earth's magnetic field over time. The orientation and intensity of the Earth's magnetic field can be obtained through the record of remanent magnetism preserved in rocks. The directions of remanent magnetization allow scientists to determine the position of the Earth's magnetic pole with respect to the study location at the time when the magnetization was acquired.  

\begin{figure}[htbp!] 
\begin{center}
\vspace{-0mm}
\includegraphics[width=\textwidth]{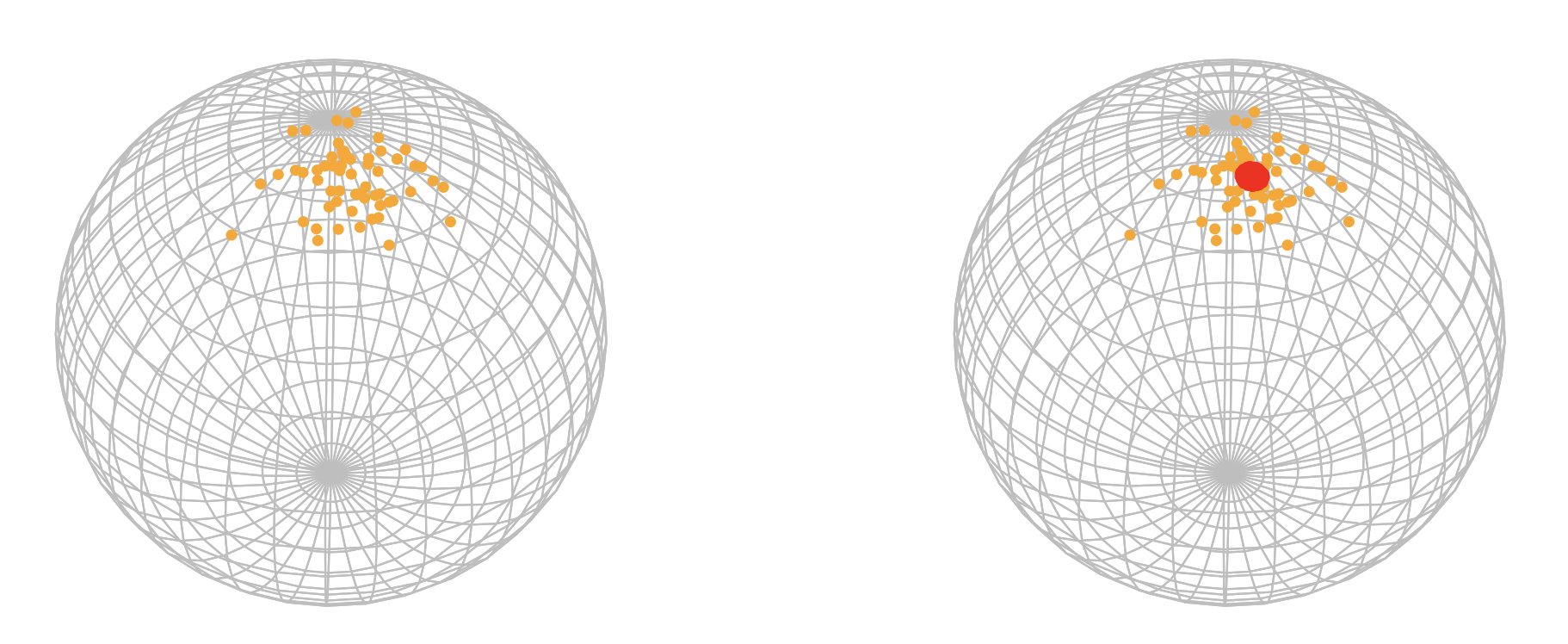} 
\end{center}
\vspace{-2mm}
\caption{(Left:) plot of the $n=62$ directions of magnetic remanence associated with the real dataset cosidered in Section~\ref{sec:realdata}. (Right:) plot of the same dataset with, in red, the corresponding 95\% confidence cap in \eqref{seccapbis}.}
\label{Fig4}
\end{figure}

\begin{figure}[htbp!] 
\begin{center}
\vspace{-10mm}
\includegraphics[height=10cm, width=.8 \textwidth]{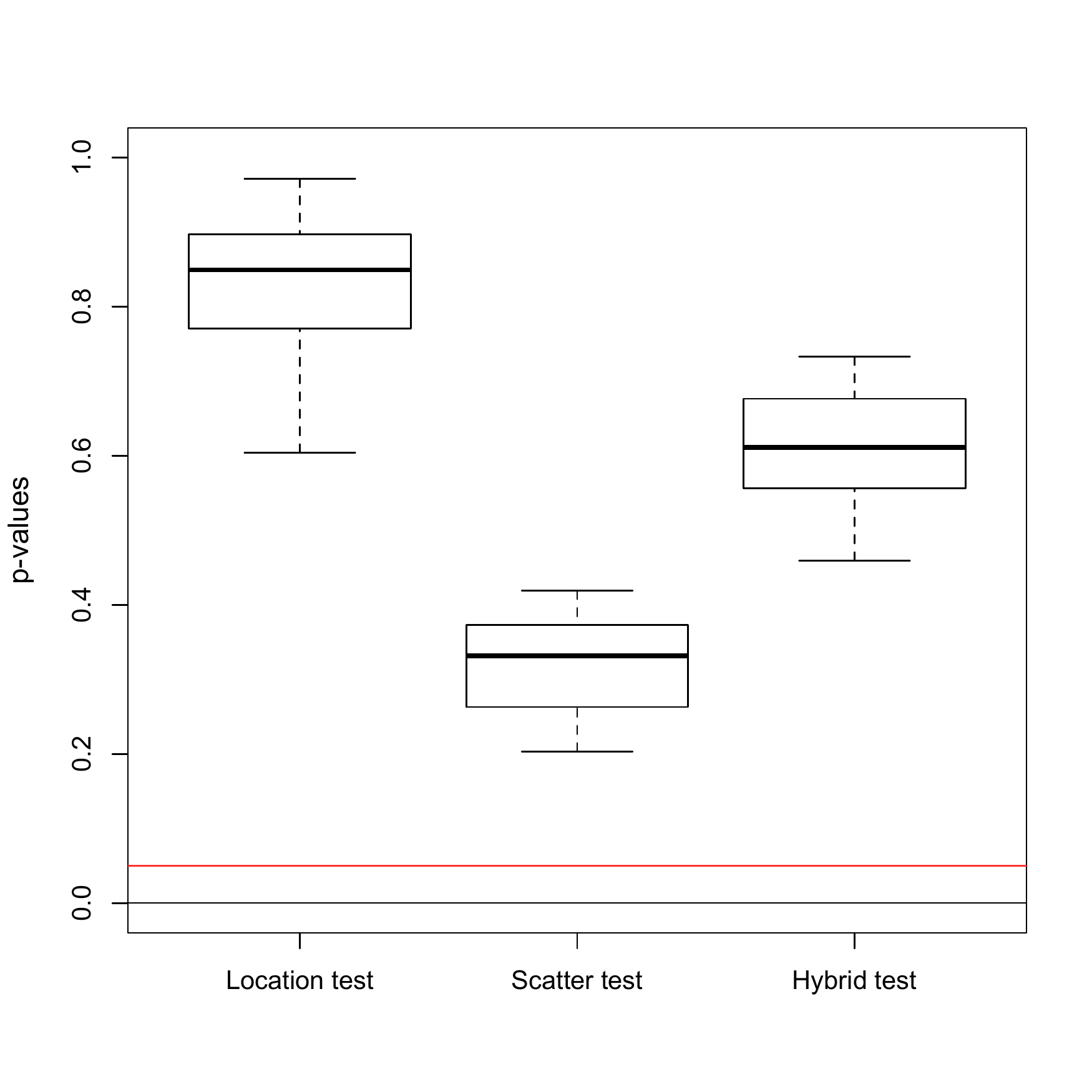} 
\end{center}
\vspace{-9mm}
\caption{Boxplots of the $p$-values of three tests of rotational symmetry over~$\mathcal{S}^2$ obtained from the 62 leave-one-out samples obtained for the magnetic remanence data set; see Section~\ref{sec:realdata} for details. The tests considered are (the unspecified-$\thetab$ versions of) the location test, scatter test and hybrid test of rotational symmetry from \cite{GPPV19}.}
\label{Fig5}
\end{figure}

We consider here a well-known dataset on~$\mathcal{S}^2$ that has already been used for inference on spherical location in \cite{Fish87}. The dataset, which is provided as Dataset~A in Appendix~B8 of this monograph, is showed in the left panel of Figure~\ref{Fig4}. Clearly, the data is highly concentrated. In line with this, the FvML maximum likelihood estimator of the concentration parameter~$\kappa$ takes value~$\hat{\kappa}=76.12$, which is of the same order of magnitude as the sample size~$n=62$. Figure~\ref{Fig4} also suggests that rotational symmetry is a plausible assumption. To assess this, we performed the three tests of rotational symmetry on~$\mathcal{S}^2$  that were recently proposed in \cite{GPPV19}: a location test and a scatter test, that respectively show power against location-type alternatives and scatter-type alternatives to rotational symmetry (we refer to \citealp{GPPV19} for details), as well as a hybrid test that shows power against both types of alternatives. These three tests, that are meant to test the null hypothesis of rotational symmetry about an unspecified location~$\thetab$, provided the $p$-values~$.844$, $.305$ and $.607$, respectively, hence did not lead to rejection at any usual nominal level. To somewhat assess the robustness of this result, we performed the following analysis: on the~$62$ samples of size~$61$ obtained by leaving one of the original observations out, we performed the same three tests of rotational symmetry and provided in Figure~\ref{Fig5} the boxplots of the 62 $p$-values obtained for each of the three tests. Again, at any usual nominal level, none of these subsamples led any of the three tests to reject the null hypothesis of rotational symmetry. 

The various statistical methods studied in this paper are therefore perfectly suitable for the present dataset. To illustrate one of these methods, we computed the 95$\%$ confidence cap for the spherical location defined in~\eqref{seccapbis}. The resulting confidence cap is showed in the right panel of Figure~\ref{Fig4}. This confidence zone is centered at the spherical mean~$\hat{\thetab}=(.210, .104, .972)\pr$ and, as expected in the present high concentration setup, has a very small size.


\section{Wrap up} 
\label{sec:wrapup} 

We discussed inference on the location parameter of rotationally symmetric distributions under high concentration. We did so by considering double asymptotic scenarios where the underlying concentration parameter~$\kappa_n$ diverges to infinity at an arbitrary rate with the sample size~$n$. This significantly improves over the state of the art for directional inference under high concentration, since previous works not only focused on a parametric class of distributions (namely, the FvML one) but also restricted to asymptotics as~$\kappa$ diverges to infinity with~$n$ fixed.  Our asymptotic results indicate that standard fixed-$\kappa$ methods are robust to high concentration, in the sense that they will remain valid in the aforementioned double asymptotic scenarios: the spherical mean remains consistent and asymptotically normal, whereas the Watson and Wald tests still asymptotically meet the level constraint. Under high concentration, however, these statistical procedures enjoy faster consistency rates than in the standard fixed-$\kappa$ asymptotic scenario. Remarkably, these consistency rates depend on the type of rotationally symmetric distributions considered, that is, they depend on the underlying angular function~$f$; this dependence is such that the higher the concentration, the faster the consistency rates. In contrast with all previous works on high concentration, we also considered optimality issues. We showed that, under mild assumptions on~$f$, the aforementioned inference procedures enjoy strong, Le Cam-type, optimality properties. For some (not all) angular functions, optimality requires that~$\kappa_n$ diverges to infinity sufficiently fast as a function of~$n$; the corresponding restriction, as we have seen, is extremely mild for our running example associated with~$f_b(z)=\exp(z^b)$, as optimality, for~$b\in(\frac{1}{2},1)$ holds in particular when~$\kappa_n$ diverges to infinity at least as fast as~$(\log n)^2$, whereas no restriction of this sort is required for~$b\geq 1$, hence in particular for the usual FvML case.


%



\appendix
\section{Proofs}

\subsection{Proof of Theorem~\ref{Theordumal}}
\label{Supsec1}

 The proof requires the following preliminary result.

\begin{Lem}
\label{Lemdumal}	
If~$\kappa \varphi_f(\kappa)\nearrow\infty$ (resp., $\kappa \varphi_f(\kappa)\searrow 0$) as~$\kappa\to\infty$, then there exists~$z_0$ such that~$f$ is convex (resp., concave) in~$[z_0,\infty)$.   
\end{Lem}

{\sc Proof of Lemma~\ref{Lemdumal}.}
Assume that~$\kappa \varphi_f(\kappa)\nearrow\infty$. Pick~$z_0$ large enough so that, in~$[z_0,\infty)$, $z\mapsto z\varphi_f(z)$ is monotone non-decreasing and takes its values in~$[1,\infty)$. Then, letting~$g(z):=z/f(z)$, the mean value theorem implies that, for any~$a,b$ with~$z_0\leq a<b$, 
\begin{eqnarray*}
0
\leq 
b\varphi_f(b)-a\varphi_f(a)
&=&
(g(b)-g(a))f'(b)+g(a)(f'(b)-f'(a))
\\[2mm]
&=&
\frac{1-c\varphi_f(c)}{f(c)} (b-a) f'(b)+g(a)(f'(b)-f'(a))
,
\end{eqnarray*}
for some~$c\in(a,b)$. Since~$c\varphi(c)\geq 1$, we must have~$f'(b)\geq f'(a)$. Therefore, $f'$ is monotone non-decreasing in~$[z_0,\infty)$, so that~$f$ is convex on the same set. The proof is entirely similar for the case~$\kappa \varphi_f(\kappa)\searrow 0$, where~$z_0$ is taken so that, in~$[z_0,\infty)$, $z\mapsto z\varphi_f(z)$ is monotone non-increasing and takes its values in~$[0,1]$. 
\cqfd
\vspace{3mm}

{\sc Proof of Theorem~\ref{Theordumal}.}
Writing 
$$
A_\kappa
:=
\int_{1-\varepsilon}^{1} (1-s^2)^{(p-3)/2} f(\kappa s)\,ds
\quad\textrm{ and }
\quad
B_\kappa
:=
\int_{-1}^{1-\varepsilon} (1-s^2)^{(p-3)/2} f(\kappa s)\,ds
,
$$
note that~$f$ provides high concentration if and only if~$A_\kappa/(A_\kappa+B_\kappa)\to 1$ as~$\kappa\to\infty$, or equivalently, if and only if~$A_\kappa/B_\kappa\to \infty$ as~$\kappa\to\infty$. In this proof, $C$ denotes a positive quantity that does not depend on~$\kappa$ and whose value may change from line to line. 

(i) Assume that~$\kappa \varphi_f(\kappa)\nearrow\infty$. Without loss of generality, restrict then to~$\kappa\geq \kappa_0$, where~$\kappa_0$ is such that~$f$ is convex in~$[\kappa_0(1-\varepsilon),\infty)$ (Lemma~\ref{Lemdumal}). 
%
Then, using the fact that $(1-s^2)^{(p-3)/2}  (s-(1-\varepsilon))$ is positive for $s \in (1- \varepsilon, 1)$, we have
\begin{eqnarray*}
	\lefteqn{
	A_\kappa
\geq
\int_{1-\varepsilon}^{1} 
(1-s^2)^{(p-3)/2} \{ f(\kappa (1-\varepsilon)) + \kappa(s-(1-\varepsilon)) f'(\kappa (1-\varepsilon) )\}
\,ds
}
\\[2mm]
& & 
\hspace{-7mm} 
\geq
\kappa f'(\kappa (1-\varepsilon) )
\int_{1-\varepsilon}^{1} 
(1-s^2)^{(p-3)/2}  (s-(1-\varepsilon))
\,ds
=
C
\kappa(1-\varepsilon) f'(\kappa (1-\varepsilon) ).
\end{eqnarray*}
Since
$$
B_\kappa
\leq
f(\kappa (1-\varepsilon))
\int_{-1}^{1-\varepsilon} 
(1-s^2)^{(p-3)/2} 
\,ds
=
C
f(\kappa (1-\varepsilon))
,
$$
we conclude that
$$
\frac{A_\kappa}{B_\kappa}
\geq
C
\kappa(1-\varepsilon) \varphi_f( \kappa(1-\varepsilon) )
\to \infty
$$
as~$\kappa$ diverges to infinity, so that~$f$ provides high concentration.
\vspace{3mm}
 
(ii) Assume that~$\kappa \varphi_f(\kappa)\to c$ for some~$c>0$, that is, $z \varphi_f(z)=c+o(1)$ as~$z\to\infty$. This means that~$\varphi_f(z)-c/z=(\log f(z)-c\log z)'=g(z)$ for a function~$g$ that satisfies~$g(z)=o(1/z)$ as~$z\to\infty$, hence that is integrable in a neighborhood of~$\infty$. For~$z_0$ large enough so that~$g(z)\leq 1$ for~$z\geq z_0$ and~$g$ is integrable in~$[z_0,\infty)$, we then have
$$
\log \Big(\frac{f(z)}{z^c}\Big)
= 
\log \Big(\frac{f(z_0)}{z_0^c}\Big)
+
\int_{z_0}^z
g(y)
\,dy
=
C+o(1)
$$
as~$z\to\infty$, which rewrites  
\begin{equation}
	\label{re}	
	f(z)
=
C z^c
+
o(z^c)
\end{equation}
for some constant~$C$ as~$z\to\infty$. This entails that, for any~$0< a < b \leq 1$, 
\begin{eqnarray*}
\lefteqn{
\Bigg|
\frac{1}{\kappa^{c}}
\int_{a}^{b} 
(1-s^2)^{(p-3)/2} 
(f(\kappa s)-C(\kappa s)^c) \, ds
\Bigg|
}
\\[2mm]
& & \hspace{23mm} 
\leq 
\frac{1}{\kappa^{c}}
\bigg(\sup_{z\in[\kappa a,\kappa b]} |f(z)-Cz^c|\bigg)
\int_{a}^{b}  
(1-s^2)^{(p-3)/2} 
\,ds
\\[2mm]
& & \hspace{23mm} 
\leq 
\frac{\pi}{2}
\sup_{z\in[\kappa a,\kappa b]} \bigg|\frac{f(z)-Cz^c}{z^c}\bigg|
\to 0
\end{eqnarray*}
as~$\kappa\to\infty$. Fixing~$\varepsilon\in(0,1/2)$, this implies that
\begin{eqnarray*}
\lefteqn{	
\frac{A_\kappa}{B_\kappa}
\leq
\frac{\int_{1-\varepsilon}^{1} (1-s^2)^{(p-3)/2} f(\kappa s)\,ds}
{\int_{1-2\varepsilon}^{1-\varepsilon} (1-s^2)^{(p-3)/2} f(\kappa s)\,ds}
=
\frac{\kappa^{-c}\int_{1-\varepsilon}^{1} (1-s^2)^{(p-3)/2} f(\kappa s)\,ds}
{\kappa^{-c}\int_{1-2\varepsilon}^{1-\varepsilon} (1-s^2)^{(p-3)/2} f(\kappa s)\,ds}
}
\\[2mm]
& & \hspace{-3mm} 
=
\frac{
\kappa^{-c}\int_{1-\varepsilon}^{1} (1-s^2)^{(p-3)/2} C(\kappa s)^c\,ds+o(1)
}
{
\kappa^{-c}\int_{1-2\varepsilon}^{1-\varepsilon} (1-s^2)^{(p-3)/2} C(\kappa s)^c\,ds+o(1)
}
\to
\frac{\int_{1-\varepsilon}^{1} (1-s^2)^{(p-3)/2} s^c\,ds}
{\int_{1-2\varepsilon}^{1-\varepsilon} (1-s^2)^{(p-3)/2} s^c\,ds}
,
\end{eqnarray*}
so that~$A_\kappa/B_\kappa=O(1)$ as~$\kappa\to\infty$, which shows that~$f$ does not provide high concentration. 
\vspace{3mm}

(iii)
Assume that~$\kappa \varphi_f(\kappa) \searrow 0$. Fix~$\tilde{\varepsilon}>\varepsilon$ and restrict, without loss of generality, to~$\kappa\geq \kappa_0$, where~$\kappa_0$ is such that~$f$ is concave in~$[\kappa_0(1-\tilde{\varepsilon}),\infty)$ (Lemma~\ref{Lemdumal}). Concavity ensures that
$$
A_\kappa
\leq
f(\kappa)
\int_{1-\varepsilon}^{1} (1-s^2)^{(p-3)/2}\,ds
=
C
f(\kappa)
\leq 
C
\big(
f( \kappa(1-\tilde{\varepsilon}) )
+
\kappa \tilde{\varepsilon} 
f'( \kappa(1-\tilde{\varepsilon}) )
\big)
.
$$
Since 
$$
B_\kappa
\geq
\int_{1-\tilde{\varepsilon}}^{1-\varepsilon} 
(1-s^2)^{(p-3)/2} f(\kappa s)
\,ds
\geq
C f(\kappa(1-\tilde{\varepsilon}))
,
$$
we obtain
$$
\frac{A_\kappa}{B_\kappa}
\leq
C
+
C
\kappa(1-\varepsilon) \varphi_f( \kappa(1-\varepsilon) )
\to C<\infty
$$
as~$\kappa$ diverges to infinity, so that~$f$ does not provide high concentration.
\cqfd


\subsection{Proof of Proposition~\ref{PropCondexpb}}
\label{Supsec2}

The proof of Proposition~\ref{PropCondexpb} requires both following preliminary results.

\begin{Lem}
\label{lemDCT}
For any~$\xi,\zeta>-1$, 
$$
\int_0^{2c} 
z^{\xi}
\big(2-{\textstyle{\frac{z}{c}}}\big)^\zeta 
e^{-z}
 \,dz
 \to
2^\zeta \Gamma(\xi+1)
$$
as~$c\to\infty$.
\end{Lem}

{\sc Proof of Lemma~\ref{lemDCT}.}
Letting $y=2-z/c$ (i.e., $z=c(2-y)$), we have
$$
\int_0^{2c} 
z^{\xi}
\big\{
\big(2-{\textstyle{\frac{z}{c}}}\big)^\zeta 
-
2^\zeta
\big\}
e^{-z}
 \,dz
=
\int_0^{2} 
c^{\xi+1}
(2-y)^\xi
\big\{
y^\zeta 
-
2^\zeta
\big\}
 e^{-c(2-y)}
 \,dy
 .
$$
For any~$y\in(0,2)$, we have that  
$$
c^{\xi+1}
(2-y)^\xi
\big|
y^\zeta 
-
2^\zeta
\big|
 e^{-c(2-y)}
\leq
(\xi+1)^{\xi+1}
 e^{-(\xi+1)}
 \frac{
\big|
y^\zeta
-
2^\zeta
\big|
}{2-y}
$$
for any~$c>0$. The result then follows from the Lebesgue Dominated Convergence Theorem. 
\cqfd

\begin{Lem}
\label{lemtec}
(i) For~$b\in(0,1)$, $(1-r)^b -1+br\leq 0$ for any~$r\in[0,2]$ $($recall that~$z^b:={\rm sgn}(z)|z|^b)$.  (ii) For~$b\geq 1$, there exists~$c\in(0,1)$ 
such that
$
0\leq (1-r)^b -1+br \leq cbr
$
for any~$r\in[0,2]$. (iii) For~$b>0$, there exists~$C>0$ such that
$|(1-r)^b-1+br| \leq Cr^2$ for any~$r\in[0,2]$. 
\end{Lem}

{\sc Proof of Lemma~\ref{lemtec}.} 
(i)
Fix~$b\in(0,1)$ and put~$g(r)=(1-r)^b -1+br$. For~$r\in(0,1)$, 
$
g'(r)
=-b(1-r)^{b-1}+b
\leq 0
$
and for~$r\in(1,2)$, 
$
g'(r)
=(-(r-1)^b -1+br)'
=-b(r-1)^{b-1}+b
\leq 0
.
$
Since~$g$ is continuous over~$[0,2]$, this implies that~$g$ is monotone non-increasing over~$[0,2]$. The result thus follows from the fact that~$g(0)=0$. 
\vspace{3mm}
 
(ii)
Fix~$b\geq 1$. Then~$g'(r)\geq 0$ for any~$r\in(0,1)\cup(1,2)$. The continuity of~$g$ over~$[0,2]$ and the fact that~$g(0)=0$ thus imply that~$g(r)\geq 0$ for any~$r\in(0,2)$. It remains to show that there exists~$c\in(0,1)$ such that
$
(1-r)^b +br-1
\leq 
cbr
$ for any~$r\in[0,2]$, or equivalently, that there exists a positive integer~$k$ for which
\begin{equation}
	\label{toshowk}
h_k(r)
:=
(1-r)^b -1+br - \frac{k-1}{k} br\leq 0
\end{equation}
for any~$r\in[0,2]$. Clearly, $h_k(r)\to h(r):=(1-r)^b-1$ as~$k\to\infty$ and the convergence is uniform in~$r\in[0,2]$. Since~$h_2'$ is right-continuous at~$0$ and satisfies
$$
h_2'(0)
=
-\frac{b}{2}
<
0
,
$$
there exists~$\eta>0$ such that~$h'_2(r)<0$ for all~$r\in(0,\eta]$, which (since~$h_2(0)=0$) yields~$h_2(r)<0$ for all~$r\in(0,\eta]$. Since, for any~$r\in(0,2]$,~$h_k(r)$ is monotone decreasing in~$k$, we deduce that~$h_k(r)<0$ for all~$r\in(0,\eta]$ and all~$k\geq 2$. Now, put~$\varepsilon:=-h_2(\eta)>0$. The uniform convergence of~$(h_k)$ to~$h$ ensures that there exists~$k_0$ such that~$|h_{k_0}(r)-h(r)|<\varepsilon/2$ for any~$r\in[0,2]$. This and the fact that~$h$ is monotone decreasing in~$[0,2]$ implies that, for any~$r\in[\eta,2]$,
$$
h_{k_0}(r) 
< h(r) + \varepsilon/2 
\leq h(\eta) + \varepsilon/2 
< h_{k_0}(\eta) + \varepsilon 
\leq h_{2}(\eta) + \varepsilon 
= 0
.
$$
We conclude that~$h_{k_0}(0)=0$ and~$h_{k_0}(r)<0$ for any~$r\in(0,2]$, so that~(\ref{toshowk}) holds for~$k=k_0$.
\vspace{3mm}
 
(iii) The Cauchy formula for the remainder of Taylor expansions yields that, for any~$r\in[0,\frac{1}{2}]$, we have
$
(1-r)^b-1+br
=
\frac{1}{2}b(b-1)(1-\eta_{b,r} r)^{b-2} r^2
$ 
for some~$\eta_{b,r}\in(0,1)$. This implies that there exists~$c_1>0$ such that
$$
|(1-r)^b-1+br|
\leq c_1 r^2
$$ 
for any~$r\in[0,\frac{1}{2}]$. Now, the mapping~$r\mapsto (1-r)^b-1+br$ is continuous over~$r\in[\frac{1}{2},2]$, so that, for any~$r\in[\frac{1}{2},2]$, we have 
$$
|(1-r)^b-1+br|\leq c_2\leq 4c_2 r^2
.
$$ 
The claim therefore holds with~$C:=\max(c_1,4c_2)$. 
\cqfd    
\vspace{3mm}


{\sc Proof of Proposition~\ref{PropCondexpb}.}
We only need to prove that Condition~(\ref{conds}) holds for any~$b>0$ (the other conditions are indeed trivially fulfilled). To do so, fix~$b>0$ and note that, for~$f_b(z)=\exp(z^b)$,~(\ref{conds}) rewrites 
 $$
\int_{-1}^1
(1-s)^\xi
(1+s)^\zeta
\big|
e^{(\kappa s)^b- \kappa^b}
-
e^{(s-1) b\kappa^b}
\big|
 \,ds
=
o\bigg( \frac{1}{\kappa^{b(\xi+1)}} \bigg)
$$
(in this proof, all convergences are as~$\kappa\to\infty$), that is, letting~$s=1-r$, 
\begin{equation}
	\label{toprzz}
\int_{0}^2
r^\xi
(2-r)^\zeta
\big|
e^{\kappa^b\{(1-r)^b-1\}}
-
e^{-b\kappa^b r}
\big|
 \,dr
=
o\bigg( \frac{1}{\kappa^{b(\xi+1)}} \bigg)
.
\end{equation}
If~$b\in(0,1)$, then Parts~(i) and~(iii) of Lemma~\ref{lemtec} and the mean value theorem yield (below, $\eta_{b,r}\in(0,1)$) 
\begin{eqnarray*}
\big|
e^{\kappa^b\{(1-r)^b-1\}}
-
e^{-b\kappa^b r}
\big|
&\!\!\!=\!\!\!&
\big(
1
-
e^{\kappa^b\{(1-r)^b-1+br\}}
\big)
e^{-b\kappa^b r}
\\[2mm]
&\!\!\!=\!\!\!&
\kappa^b
| (1-r)^b-1+br |
e^{\eta_{b,r}\kappa^b\{(1-r)^b-1+br\}}
e^{-b\kappa^b r}
\\[2mm]
&\!\!\!=\!\!\!&
C
\kappa^b
r^2
e^{-b\kappa^b r}
.
\end{eqnarray*}
Now, if~$b\geq 1$, then Lemma~\ref{lemtec}(ii)--(iii) and the mean value theorem yield
\begin{eqnarray*}
\big|
e^{\kappa^b\{(1-r)^b-1\}}
-
e^{-b\kappa^b r}
\big|
&\!\!\!=\!\!\!&
\big(
e^{\kappa^b\{(1-r)^b-1+br\}}
-
1
\big)
e^{-b\kappa^b r}
\\[2mm]
&\!\!\!=\!\!\!&
\kappa^b
|(1-r)^b-1+br|
e^{\eta_{b,r}\kappa^b\{(1-r)^b-1+br\}}
e^{-b\kappa^b r}
\\[2mm]
&\!\!\!=\!\!\!&
C \kappa^b r^2
e^{(c-1)b \kappa^b r}
.
\end{eqnarray*}
We therefore showed that, for any~$b>0$, there exists~$K>0$ such that
$$
\big|
e^{\kappa^b\{(1-r)^b-1\}}
-
e^{-b\kappa^b r}
\big|
=
C \kappa^b r^2
e^{-K \kappa^b r}
.
$$ 
By letting~$z=K\kappa^b r$, this yields
 \begin{eqnarray*}
\lefteqn{
\hspace{-12mm} 
\int_{0}^2
r^\xi
(2-r)^\zeta
\big|
e^{\kappa^b\{(1-r)^b-1\}}
-
e^{-b\kappa^b r}
\big|
 \,dr
\leq
C
\kappa^b
\int_{0}^2
r^{\xi+2}
(2-r)^\zeta
e^{-K \kappa^b r}
 \,dr
}
\\[2mm]
& & 
\hspace{-3mm} 
\leq
\frac{C}{K^{\xi+3}\kappa^{b(\xi+2)}} 
\int_{0}^{2K\kappa^b}
z^{\xi+2}
(2-{\textstyle{\frac{z}{K\kappa^b}}})^\zeta
e^{-z}
 \,dz
=
O\bigg(\frac{1}{\kappa^{b(\xi+2)}}\bigg) 
,
\end{eqnarray*}
where we used Lemma~\ref{lemDCT}. This proves~(\ref{toprzz}), hence establishes the result. 
\cqfd
\vspace{3mm}


\subsection{Proof of Theorem~\ref{Theorexpansions}}
\label{SupSec3}

The proof crucially relies on the following lemma.

\begin{Lem}
\label{LemExpansionConstant}
Fix an integer~$p\geq 2$ and~$f\in\mathcal{F}$. Let~$(\kappa_n)$ be a positive real sequence that diverges to~$\infty$. Then,
$$
(i)
\qquad
\frac{1}{f(\kappa_n)c_{p,\kappa_n,f}}
=
\frac{
2^{(p-3)/2}  
\Gamma\big( {\textstyle{\frac{p-1}{2}}} \big)
}
{
(\kappa_n \varphi_f(\kappa_n))^{(p-1)/2}
}
+
o\bigg(\frac{1}{(\kappa_n \varphi_f(\kappa_n))^{(p-1)/2}}\bigg)
$$
and
$$
(ii)
\qquad
\frac{\tilde{e}_{n2}}{(f(\kappa_n))^2 c^2_{p,\kappa_n,f}}
=
\frac{2^{p-4} (p-1) \Gamma^2\big(\textstyle{\frac{p-1}{2}}\big)
}{(\kappa_n \varphi_f(\kappa_n))^{p+1}}
+ 
o\bigg( 
\frac{1}{(\kappa_n \varphi_f(\kappa_n))^{p+1}}
 \bigg)
$$
as~$n\to\infty$.
\end{Lem}
\vspace{1mm}

{\sc Proof of Lemma~\ref{LemExpansionConstant}}.
(i) Write
$$
\frac{1}{f(\kappa_n)c_{p,\kappa_n,f}}
=
\frac{1}{f(\kappa_n)}
\int_{-1}^1  (1-s^2)^{(p-3)/2} f(\kappa_n s)\,ds
=
T_{n1}
+
T_{n2}
,
$$
with
$$
T_{n1}
:=
\int_{-1}^{1} 
(1-s^2)^{(p-3)/2} 
e^{(s-1) \kappa_n\varphi_f(\kappa_n)}
\,ds
$$
and
$$
T_{n2}
:=
\int_{-1}^{1} 
(1-s^2)^{(p-3)/2} 
\big(
e^{\log f(\kappa_n s)- \log f(\kappa_n)}
-
e^{(s-1) \kappa_n\varphi_f(\kappa_n)}
\big)
\,ds
.
$$
Letting~$z=(1-s) \kappa_n \varphi_f(\kappa_n)$, Lemma~\ref{lemDCT} readily yields
\begin{eqnarray*}
T_{n1}
&\!\!\!=\!\!\!&
\frac{1}{(\kappa_n\varphi_f(\kappa_n))^{(p-1)/2}}
\int_{0}^{2\kappa_n\varphi_f(\kappa_n)} 
z^{(p-3)/2} 
(2-{\textstyle{\frac{z}{\kappa_n\varphi_f(\kappa_n)}}})^{(p-3)/2} 
e^{-z}
\,dz
\\[2mm]
&\!\!\!=\!\!\!&
\frac{
2^{(p-3)/2}  
\Gamma\big( {\textstyle{\frac{p-1}{2}}} \big)
}
{
(\kappa_n \varphi_f(\kappa_n))^{(p-1)/2}
}
+
o\bigg(\frac{1}{(\kappa_n \varphi_f(\kappa_n))^{(p-1)/2}}\bigg)
,
\end{eqnarray*}
Since~(\ref{conds}) ensures that
$$
T_{n2}
=
o\bigg( \frac{1}{(\kappa_n \varphi_f(\kappa_n))^{(p-1)/2}} \bigg)
,
$$
the result follows. 
\vspace{1mm}

(ii) Using the U-statistic formulation of the variance, we have
\begin{eqnarray*}
\lefteqn{
\hspace{-6mm}
	\frac{\tilde{e}_{n2}}{(f(\kappa_n))^2 c^2_{p,\kappa_n,f}}
=
\frac{1}{2(f(\kappa_n))^2}
\int_{-1}^1 \int_{-1}^1 
(s-\tilde{s})^2 
(1-s^2)^{(p-3)/2} 
(1-\tilde{s}^2)^{(p-3)/2}
}
\\[2mm]
& & 
\hspace{33mm} 
\times
  f(\kappa_n s) f(\kappa_n \tilde{s})
\,ds d\tilde{s}
=
S_{n1}+2S_{n2}+S_{n3}
,
\end{eqnarray*}
where
\begin{eqnarray*}
\lefteqn{
\hspace{-0mm} 
S_{n1}
:=
\frac{1}{2}
\int_{-1}^1 \int_{-1}^1 
(s-\tilde{s})^2 
}
\\[2mm]
& & 
\hspace{13mm} 
\times
(1-s^2)^{(p-3)/2} 
(1-\tilde{s}^2)^{(p-3)/2} 
e^{(s-1) \kappa_n\varphi_f(\kappa_n)}
e^{(\tilde{s}-1) \kappa_n \varphi_f(\kappa_n)}
\,ds d\tilde{s}
,
\end{eqnarray*}
\begin{eqnarray*}
\lefteqn{
\hspace{-2mm} 
S_{n2}
:=
\frac{1}{2}
\int_{-1}^1 \int_{-1}^1 
(s-\tilde{s})^2 
(1-s^2)^{(p-3)/2} (1-\tilde{s}^2)^{(p-3)/2} 
e^{(s-1) \kappa_n\varphi_f(\kappa_n)}
}
\\[2mm]
& & 
\hspace{43mm} 
\times
\big(
e^{\log f(\kappa_n \tilde{s})- \log f(\kappa_n)}
-
e^{(\tilde{s}-1) \kappa_n\varphi_f(\kappa_n)}
\big)
\,ds d\tilde{s}
\end{eqnarray*}
and
\begin{eqnarray*}
\lefteqn{
\hspace{-14mm} 
S_{n3}
:=
\frac{1}{2}
\int_{-1}^1 \int_{-1}^1 
(s-\tilde{s})^2 
(1-s^2)^{(p-3)/2} (1-\tilde{s}^2)^{(p-3)/2} 
}
\\[2mm]
& & 
\hspace{10mm} 
\times
\big(
e^{\log f(\kappa_n s)- \log f(\kappa_n)}
-
e^{(s-1) \kappa_n\varphi_f(\kappa_n)}
\big)
\\[2mm]
& & 
\hspace{18mm} 
\times
\big(
e^{\log f(\kappa_n \tilde{s})- \log f(\kappa_n)}
-
e^{(\tilde{s}-1) \kappa_n\varphi_f(\kappa_n)}
\big)
\,ds d\tilde{s}
.
\end{eqnarray*}

We start with~$S_{n1}$. Letting~$z=(1-s) \kappa_n \varphi_f(\kappa_n)$ and~$\tilde{z}=(1-\tilde{s}) \kappa_n \varphi_f(\kappa_n)$, we obtain
\begin{eqnarray*}
\lefteqn{
S_{n1}
=
\frac{1}{2(\kappa_n \varphi_f(\kappa_n))^{p+1}}
\int_{0}^{2\kappa_n \varphi_f(\kappa_n)} 
\int_{0}^{2\kappa_n \varphi_f(\kappa_n)} 
(z-\tilde{z})^2 
(2-{\textstyle{\frac{z}{\kappa_n\varphi_f(\kappa_n)}}})^{(p-3)/2} 
}
\\[2mm]
& & 
\hspace{23mm} 
\times
(2-{\textstyle{\frac{\tilde{z}}{\kappa_n\varphi_f(\kappa_n)}}})^{(p-3)/2} 
z^{(p-3)/2} 
\tilde{z}^{(p-3)/2} 
e^{-z}
e^{-\tilde{z}}
\,dz d\tilde{z}
\\[2mm]
& & 
\hspace{0mm} 
=
\frac{1}{(\kappa_n \varphi_f(\kappa_n))^{p+1}}
\bigg(
\int_{0}^{2\kappa_n \varphi_f(\kappa_n)} 
z^{(p+1)/2} 
(2-{\textstyle{\frac{z}{\kappa_n\varphi_f(\kappa_n)}}})^{(p-3)/2}e^{-z}
\,dz
\bigg)
\\[2mm]
& & 
\hspace{23mm} 
\times
\bigg(
\int_{0}^{2\kappa_n \varphi_f(\kappa_n)} 
z^{(p-3)/2}
(2-{\textstyle{\frac{z}{\kappa_n\varphi_f(\kappa_n)}}})^{(p-3)/2}e^{-z}
\,dz
\bigg)
\\[2mm]
& & 
\hspace{3mm} 
-
\frac{1}{(\kappa_n \varphi_f(\kappa_n))^{p+1}}
\bigg(
\int_{0}^{2\kappa_n \varphi_f(\kappa_n)} 
z^{(p-1)/2} 
(2-{\textstyle{\frac{z}{\kappa_n\varphi_f(\kappa_n)}}})^{(p-3)/2}e^{-z}
\,dz
\bigg)^2
,
\end{eqnarray*}
so that Lemma~\ref{lemDCT} provides
\begin{eqnarray*}
S_{n1}
&\!\!\!=\!\!\!&
\frac{2^{p-3}}{(\kappa_n\varphi_f(\kappa_n))^{p+1}}
\big(
\Gamma\big({\textstyle{\frac{p+3}{2}}}\big)
\Gamma\big({\textstyle{\frac{p-1}{2}}}\big)
-
\Gamma^2\big({\textstyle{\frac{p+1}{2}}}\big)
\big)
+
o\bigg(\frac{1}{(\kappa_n \varphi_f(\kappa_n))^{p+1}}\bigg)
\\[2mm]
&\!\!\!=\!\!\!&
\frac{2^{p-3}}{(\kappa_n\varphi_f(\kappa_n))^{p+1}}
\big(
\big({\textstyle{\frac{p+1}{2}}}\big)
\big({\textstyle{\frac{p-1}{2}}}\big)
-
\big({\textstyle{\frac{p-1}{2}}}\big)^2
\big)
\Gamma^2\big({\textstyle{\frac{p-1}{2}}}\big)
+
o\bigg(\frac{1}{(\kappa_n\varphi_f(\kappa_n))^{p+1}}\bigg)
\\[3mm]
&\!\!\!=\!\!\!&
\frac{2^{p-4}(p-1) \Gamma^2\big(\textstyle{\frac{p-1}{2}}\big)
}{(\kappa_n\varphi_f(\kappa_n))^{p+1}}
+
o\bigg(\frac{1}{(\kappa_n\varphi_f(\kappa_n))^{p+1}}\bigg)
.
\end{eqnarray*}

We turn to~$S_{n2}$. Upper-bounding~$(s-\tilde{s})^2=((1-s)-(1-\tilde{s}))^2$ by~$2(1-s)^2+2(1-\tilde{s})^2$, we obtain
\begin{eqnarray*}
\lefteqn{
S_{n2}
\leq
\bigg(
\int_{-1}^1 
(1-s)^2 
(1-s^2)^{(p-3)/2} 
e^{(s-1) \kappa_n\varphi_f(\kappa_n)}
\,ds
\bigg)
}
\\[2mm]
& & 
\hspace{3mm} 
\times
\bigg(
\int_{-1}^1 
(1-s^2)^{(p-3)/2} 
\big(
e^{\log f(\kappa_n s)- \log f(\kappa_n)}
-
e^{(s-1) \kappa_n\varphi_f(\kappa_n)}
\big)
\,ds
\bigg)
\\[2mm]
& & 
\hspace{1mm} 
+
\bigg(
\int_{-1}^1 
(1-s^2)^{(p-3)/2} 
e^{(s-1) \kappa_n\varphi_f(\kappa_n)}
\,ds
\bigg)
\\[2mm]
& & 
\hspace{3mm} 
\times
\bigg(
\!
\int_{-1}^1 
(1-s)^2 
(1-s^2)^{(p-3)/2} 
\big(
e^{\log f(\kappa_n s)- \log f(\kappa_n)}
-
e^{(s-1) \kappa_n\varphi_f(\kappa_n)}
\big)
\,ds
\!
\bigg)
.
\end{eqnarray*}
Letting~$z=(1-s) \kappa_n \varphi_f(\kappa_n)$ in two of the four integrals above, \eqref{conds} yields
\begin{eqnarray*}
\lefteqn{
S_{n2}
\leq
\frac{1}{(\kappa_n \varphi_f(\kappa_n))^{(p+3)/2}}
\bigg(
\int_{0}^{2\kappa_n \varphi_f(\kappa_n)} 
z^{(p+1)/2} 
(2-{\textstyle{\frac{z}{\kappa_n\varphi_f(\kappa_n)}}})^{(p-3)/2}e^{-z}
\,dz
\bigg)
}
\\[2mm]
& & 
\hspace{3mm} 
\times
\bigg(
\int_{-1}^1 
(1-s^2)^{(p-3)/2} 
\big(
e^{\log f(\kappa_n s)- \log f(\kappa_n)}
-
e^{(s-1) \kappa_n\varphi_f(\kappa_n)}
\big)
\,ds
\bigg)
\\[2mm]
& & 
\hspace{2mm} 
+
\frac{1}{(\kappa_n \varphi_f(\kappa_n))^{(p-1)/2}}
\bigg(
\int_{0}^{2\kappa_n \varphi_f(\kappa_n)} 
z^{(p-3)/2} 
(2-{\textstyle{\frac{z}{\kappa_n\varphi_f(\kappa_n)}}})^{(p-3)/2}e^{-z}
\,dz
\bigg)
\\[2mm]
& & 
\hspace{3mm} 
\times
\bigg(
\int_{-1}^1 
(1-s)^2 
(1-s^2)^{(p-3)/2} 
\big(
e^{\log f(\kappa_n s)- \log f(\kappa_n)}
-
e^{(s-1) \kappa_n\varphi_f(\kappa_n)}
\big)
\,ds
\bigg)
\\[2mm]
& & 
\hspace{-3mm} 
=
O\bigg(\frac{1}{(\kappa_n \varphi_f(\kappa_n))^{(p+3)/2}}\bigg)
o\bigg(\frac{1}{(\kappa_n \varphi_f(\kappa_n))^{(p-1)/2}}\bigg)
\\[2mm]
& & 
\hspace{1mm} 
+
O\bigg(\frac{1}{(\kappa_n \varphi_f(\kappa_n))^{(p-1)/2}}\bigg)
o\bigg(\frac{1}{(\kappa_n \varphi_f(\kappa_n))^{(p+3)/2}}\bigg)
=
o\bigg(\frac{1}{(\kappa_n \varphi_f(\kappa_n))^{p+1}}\bigg)
.
\end{eqnarray*}

We treat~$S_{n3}$ by upper-bounding again~$(s-\tilde{s})^2$ by~$2(1-s)^2+2(1-\tilde{s})^2$, which yields
\begin{eqnarray*}
\lefteqn{
S_{n3}
\leq
2
\bigg(
\int_{-1}^1 
(1-s)^2
(1-s^2)^{(p-3)/2} 
\big(
e^{\log f(\kappa_n s)- \log f(\kappa_n)}
-
e^{(s-1) \kappa_n\varphi_f(\kappa_n)}
\big)
\,ds
\bigg)
}
\\[2mm]
& & 
\hspace{13mm} 
\times
\bigg(
\int_{-1}^1 
(1-s^2)^{(p-3)/2} 
\big(
e^{\log f(\kappa_n s)- \log f(\kappa_n)}
-
e^{(s-1) \kappa_n\varphi_f(\kappa_n)}
\big)
\,ds
\bigg)
\\[2mm]
& & 
\hspace{0mm} 
=
o\bigg(\frac{1}{(\kappa_n \varphi_f(\kappa_n))^{(p+3)/2}}\bigg)
o\bigg(\frac{1}{(\kappa_n \varphi_f(\kappa_n))^{(p-1)/2}}\bigg)
=
o\bigg(\frac{1}{(\kappa_n \varphi_f(\kappa_n))^{p+1}}\bigg)
.
\end{eqnarray*}
This completes the proof. 
\cqfd
\vspace{3mm}

We can now prove Theorem~\ref{Theorexpansions}. 
\vspace{3mm}

{\sc Proof of Theorem~\ref{Theorexpansions}.}
First note that Lemma~\ref{LemExpansionConstant} readily yields
\begin{eqnarray}
\lefteqn{
	\frac{1-e_{n2}}{f(\kappa_n)c_{p,\kappa_n,f}}
=
\frac{1}{f(\kappa_n)}
\int_{-1}^1  (1-s^2) (1-s^2)^{(p-3)/2} f(\kappa_n s)\,ds
}
\nonumber
\\[2mm]
& & 
\hspace{3mm} 
=
\frac{1}{f(\kappa_n)c_{p+2,\kappa_n,f}}
=
\frac{
2^{(p-1)/2}  
\Gamma\big( {\textstyle{\frac{p+1}{2}}} \big)
}
{
(\kappa_n \varphi_f(\kappa_n))^{(p+1)/2}
}
+
o\bigg(\frac{1}{(\kappa_n \varphi_f(\kappa_n))^{(p+1)/2}}\bigg)
\label{oupti}
\end{eqnarray}
and
\begin{eqnarray}
\lefteqn{
	\frac{{\rm E}\big[ (1-u_{n1}^2)^2\big]}{f(\kappa_n)c_{p,\kappa_n,f}}
=
\frac{1}{f(\kappa_n)}
\int_{-1}^1  (1-s^2)^2 (1-s^2)^{(p-3)/2} f(\kappa_n s)\,ds
}
\nonumber
\\[2mm]
& & 
\hspace{3mm} 
=
\frac{1}{f(\kappa_n)c_{p+4,\kappa_n,f}}
=
\frac{
2^{(p+1)/2}  
\Gamma\big( {\textstyle{\frac{p+3}{2}}} \big)
}
{
(\kappa_n \varphi_f(\kappa_n))^{(p+3)/2}
}
+
o\bigg(\frac{1}{(\kappa_n \varphi_f(\kappa_n))^{(p+3)/2}}\bigg)
.
\label{oupti2}
\end{eqnarray}
The result then follows by writing
$$
1-e_{n2}
=
\bigg(\frac{1-e_{n2}}{f(\kappa_n)c_{p,\kappa_n,f}}\bigg)
\Big/
\bigg(\frac{1}{f(\kappa_n)c_{p,\kappa_n,f}}\bigg)
,
$$
%
$$
\tilde{e}_{n2}
=
\bigg(\frac{\tilde{e}_{n2}}{(f(\kappa_n))^2 c^2_{p,\kappa_n,f}}\bigg)
/
\bigg(\frac{1}{f(\kappa_n)c_{p,\kappa_n,f}} \bigg)^2
$$
and
$$
{\rm E}[v_{n1}^4]={\rm E}\big[ (1-u_{n1}^2)^2\big]
=
\bigg(\frac{{\rm E}\big[ (1-u_{n1}^2)^2\big]}{f(\kappa_n)c_{p,\kappa_n,f}}\bigg)
\Big/
\bigg(\frac{1}{f(\kappa_n)c_{p,\kappa_n,f}}\bigg)
,
$$
and by using Lemma~\ref{LemExpansionConstant} along with~(\ref{oupti})--(\ref{oupti2}). 
\cqfd
\vspace{3mm}


\subsection{Proofs of Theorems~\ref{Theorsphericalmean}, \ref{TheorsphericalCap} and~\ref{Theorsphericalmean2}}
\label{SupSec4}

Several proofs of this section rely on the following uniform second-order delta method (the proof is a trivial extension of the proof of Theorem~3.8 in \citealp{van1998}).

\begin{Lem} 
\label{Lemseconddelta}	
Let~$\phi:\R^p\to\R$ be twice continuously differentiable in a neighborhood of~$\vb$. Let~$(\vb_n)$ be a sequence in~$\R^p$ converging to~$\vb$. Let~$(\Tb_n)$ be a sequence of random vectors taking their values in the domain of~$\phi$ and such that~$r_n(\Tb_n-\vb_n)$ is~$O_{\rm P}(1)$ for a sequence~$(r_n)$ that diverges to infinity. Then,
$$
r_n^2 
\Big\{
\phi(\Tb_n)-\phi(\vb_n)
- (\Tb_n-\vb_n)' \nabla\phi(\vb)
- \frac{1}{2}(\Tb_n-\vb_n)' \Hb\phi(\vb) (\Tb_n-\vb_n)
\Big\}
=
o_{\rm P}(1)
,
$$
where~$\nabla\phi(\vb)$ and~$\Hb\phi(\vb)$ denote the gradient and Hessian matrix of~$\phi$ at~$\vb$, respectively. 
\end{Lem}

Assuming that~$\sqrt{n} \tilde{e}_{n2}^{-1/4} (\bar{\Xb}_n- e_{n1} \thetab)$ is~$O_{\rm P}(1)$ (this will be proved later in this section), this lemma entails that
\begin{eqnarray}
\label{2ndorderDelta}
\lefteqn{	
\hspace{3mm} 
\frac{n}{{\tilde{e}_{n2}^{1/2}}}
\left\{
\hat{\thetab}_n- \thetab
- ({\bf I}_p - \thetab \thetab\pr)
(\bar{\Xb}_n- e_{n1} \thetab)
\right.
}
\\[2mm]
& & 
\hspace{13mm} 
\left. 
- \frac{1}{2} 
\left(
\begin{array}{c}
(\bar{\Xb}_n- e_{n1} \thetab)\pr \Hb g_1(\thetab) (\bar{\Xb}_n- e_{n1} \thetab)\\[0mm]
\vdots\\[0mm]
(\bar{\Xb}_n- e_{n1} \thetab)\pr \Hb g_p(\thetab) (\bar{\Xb}_n- e_{n1} \thetab)
\end{array}
\right)
\right\}
= 
o_{\rm P}(1)
,
\nonumber
\end{eqnarray}
where~$g_j:\R^p\setminus\{{\bf 0}\}\to\R$ is the mapping defined through~$g_j(\xb)=x_j/\|\xb\|$. Note that this in particular yields 
\begin{equation}
\label{1storderDelta}
\frac{\sqrt{n}}{{\tilde{e}_{n2}^{1/4}}}
(\hat{\thetab}_n- \thetab)
=
\frac{\sqrt{n}}{{\tilde{e}_{n2}^{1/4}}}
 ({\bf I}_p - \thetab \thetab\pr)
(\bar{\Xb}_n- e_{n1} \thetab)
+
O_{\rm P}\bigg(\frac{{\tilde{e}_{n2}^{1/4}}}{\sqrt{n}}\bigg)
.
\end{equation}

We can now prove Theorem~\ref{Theorsphericalmean}. 
\vspace{1mm}

{\sc Proof of Theorem~\ref{Theorsphericalmean}.}
In this proof, all expectations and variances are under~${\rm P}\n_{\thetab,\kappa_n,f}$ and all stochastic convergences are as~$n\to\infty$ under~${\rm P}\n_{\thetab,\kappa_n,f}$. Using the tangent-normal decomposition of~$\Xb_{ni}$ with respect to~$\thetab$, write
$$
\frac{{\sqrt{n}}(\bar{\Xb}_n- e_{n1} \thetab)}{{\tilde{e}_{n2}^{1/4}}} 
=
\frac{1}{\sqrt{n}}\sum_{i=1}^n \frac{(u_{ni}- e_{n1})}{\tilde{e}_{n2}^{1/4}} \thetab+ \frac{1}{\sqrt{n}} \sum_{i=1}^n \frac{v_{ni}}{{{\tilde{e}_{n2}^{1/4}}}} \Sb_{ni} 
 =: 
V_{n}\thetab + \Wb_n
.
$$
Since Theorem~\ref{Theorexpansions}(ii) implies that
$$
{\rm E}\bigg[
\bigg(
\frac{1}{\sqrt{n}}\sum_{i=1}^n \frac{(u_{ni}- e_{n1})}{\tilde{e}_{n2}^{1/4}}
\bigg)^2
\bigg]
=
\tilde{e}_{n2}^{1/2}
\,
{\rm Var}\bigg[
\frac{1}{\sqrt{n}}\sum_{i=1}^n \frac{(u_{ni}- e_{n1})}{{\sqrt{\tilde{e}_{n2}}}}
\bigg]
=
\tilde{e}_{n2}^{1/2}
=
o(1)
,
$$
we obtain that 
$$
\frac{{\sqrt{n}}(\bar{\Xb}_n- e_{n1} \thetab)}{{\tilde{e}_{n2}^{1/4}}} 
 =
  \Wb_n+o_{\rm P}(1)
.
$$
For any unit $p$-vector~$\ub$, write 
$$
\ub'\Wb_n
=
\sum_{i=1}^n 
\frac{v_{ni}\ub'\Sb_{ni}}{\sqrt{n}\tilde{e}_{n2}^{1/4}} 
=:
\sum_{i=1}^n Z_{ni}
.
$$
For any~$n$, the~$Z_{ni}$'s are centered \mbox{i.i.d.} random variables such that
$$
s_n^2
:=
{\rm Var}
\bigg[
\sum_{i=1}^n Z_{ni}
\bigg]
=
\frac{(1- e_{n2})\ub'( {\bf I}_p- \thetab \thetab\pr)\ub}{(p-1)\tilde{e}_{n2}^{1/2}} 
=
\frac{\sqrt{2}\ub'( {\bf I}_p- \thetab \thetab\pr)\ub}{\sqrt{p-1}} + o(1)
,
$$ 
where we used the first result in~(\ref{sjffl}). Aiming at establishing the asymptotic normality of~$\ub'\Wb_n$, the Lindeberg condition reads
\begin{equation}
	\label{Lind}
\ell_n
:=
\frac{1}{s_n^2}
\sum_{i=1}^n 
{\rm E}
\Big[ 
Z_{ni}^2 \mathbb{I}[|Z_{ni}|>\varepsilon s_n]
\Big]
=
o(1)
.
\end{equation}
Applying the Cauchy-Schwarz and Chebyshev inequalities yields 
\begin{eqnarray*}
	\lefteqn{
\ell_n
=
\frac{n}{s_n^2}
{\rm E}
\Big[ 
Z_{n1}^2 \mathbb{I}[|Z_{n1}|>\varepsilon s_n]
\Big]
\leq
\frac{n}{s_n^2}
\sqrt{{\rm E}\big[ Z_{n1}^4 \big] P[|Z_{n1}|>\varepsilon s_n]}
}
\\[2mm]
& & 
\hspace{13mm} 
\leq
\frac{n}{s_n^2}
\sqrt{\frac{{\rm E}\big[ Z_{n1}^4 \big]{\rm Var}[Z_{n1}]}{\varepsilon^2 s_n^2}}
=
\frac{1}{\varepsilon s_n^2}
\sqrt{
n {\rm E}\big[ Z_{n1}^4 \big]
}
.
\end{eqnarray*}
Since the second convergence in~(\ref{sjffl}) provides
$$
n {\rm E}\big[ Z_{n1}^4 \big]
=
n {\rm E}
\bigg[ 
\bigg( 
\frac{v_{n1}\ub'\Sb_{n1}}{\sqrt{n}\tilde{e}_{n2}^{1/4}}
\bigg)^4 
\bigg] 
\leq
\frac{{\rm E}[ v_{n1}^4]}{n\tilde{e}_{n2}}
=
O\bigg(\frac{1}{n}\bigg)
,
$$
the Lindeberg condition in~(\ref{Lind}) is satisfied, so that $s_n^{-1} \ub'\Wb_n$ is asymptotically standard normal for any unit $p$-vector~$\ub$. Consequently, 
$$
\ub'\Wb_n 
\stackrel{\mathcal{D}}{\to}
\mathcal{N}
\bigg(
0
,
\frac{\sqrt{2}\ub'( {\bf I}_p- \thetab \thetab\pr)\ub}{\sqrt{p-1}} 
\bigg)
$$ 
 for any unit $p$-vector~$\ub$, which entails that
$$
\Wb_n 
\stackrel{\mathcal{D}}{\to}
\mathcal{N}
\bigg(
{\bf 0}
,
\frac{\sqrt{2}}{\sqrt{p-1}}( {\bf I}_p- \thetab \thetab\pr) 
\bigg)
.
$$ 
It follows that 
\begin{equation}
	\label{sdf}
\frac{{\sqrt{n}}(\bar{\Xb}_n- e_{n1} \thetab)}{{\tilde{e}_{n2}^{1/4}}} 
\stackrel{\mathcal{D}}{\to}
\mathcal{N}
\bigg(
{\bf 0}
,
\frac{\sqrt{2}}{\sqrt{p-1}}( {\bf I}_p- \thetab \thetab\pr) 
\bigg)
.
\end{equation}
Therefore, 
(\ref{1storderDelta}) holds and readily yields
$$
\frac{\sqrt{n}(\hat{\thetab}_n - \thetab)}{{\tilde{e}_{n2}^{1/4}}}
\stackrel{\mathcal{D}}{\to}
\mathcal{N}
\bigg( 
{\bf 0}
,
\frac{\sqrt{2}}{\sqrt{p-1}}( {\bf I}_p- \thetab \thetab\pr) 
\bigg)
,
$$  
which, by using Theorem~\ref{Theorexpansions}(ii), provides the weak limiting result in~(\ref{consist}). The one in~(\ref{ICpre}) then follows by noting that~$1-(\thetab\pr\hat{\thetab}_n)^2=(\hat{\thetab}_n - \thetab)'({\bf I}_p- \thetab \thetab\pr)^-(\hat{\thetab}_n - \thetab)$, where~$\Ab^-$ stands the Moore-Penrose inverse of~$\Ab$.
\cqfd
\vspace{3mm}
 

{\sc Proof of Theorem~\ref{TheorsphericalCap}.}
Direct computations allow checking that the function~$g_j:\R^p\setminus\{{\bf 0}\}\to\R$ defined through~$g_j(\xb)=x_j/\|\xb\|$ has the Hessian matrix
$$
Hg_j(\xb)
=
\frac{1}{\| \xb \|^3}
\bigg[
\frac{3x_j \xb\xb'}{\| \xb \|^2}
- {\bf e}_j \xb'
- \xb {\bf e}_j'
- x_j {\bf I}_p
\bigg]
,
$$
where~${\bf e}_j$ stands for the $j$th vector of the canonical basis of~$\R^p$. Therefore, premultiplying both sides of~(\ref{2ndorderDelta}) by~$\thetab'$ yields
$$
\frac{n(\thetab'\hat{\thetab}_n- 1)}{\tilde{e}_{n2}^{1/2}}
=
\frac{n}{2\tilde{e}_{n2}^{1/2}}
(\bar{\Xb}_n- e_{n1} \thetab)\pr 
\Hb
(\bar{\Xb}_n- e_{n1} \thetab)
+
o_{\rm P}(1)
,
$$
where (the $\theta_j$'s are the components of $\thetab$)
$$
	\Hb
:=
\sum_{j=1}^p \theta_j Hg_j(\thetab)
=
 3 \thetab\thetab'
-\thetab\thetab'
-\thetab\thetab'
- {\bf I}_p
=
-({\bf I}_p-\thetab\thetab')
. 
$$
Therefore, using~(\ref{sdf}), we obtain that
\begin{eqnarray*}
\lefteqn{
\frac{\sqrt{2(p-1)}n(1-\thetab'\hat{\thetab}_n)}{\tilde{e}_{n2}^{1/2}}
=
\frac{n}{\tilde{e}_{n2}^{1/2}}
(\bar{\Xb}_n- e_{n1} \thetab)\pr 
}
\\[2mm]
& & 
\hspace{23mm} 
\times \bigg(
\frac{\sqrt{p-1}}{\sqrt{2}}
({\bf I}_p-\thetab\thetab')
\bigg)
(\bar{\Xb}_n- e_{n1} \thetab)
+
o_{\rm P}(1)
\stackrel{\mathcal{D}}{\to}
\chi^2_{p-1}
.
\end{eqnarray*}
The result then follows from Theorem~\ref{Theorexpansions}(ii). 
\cqfd
\vspace{3mm}

 
The proof of Theorem~\ref{Theorsphericalmean2} requires the following preliminary result.

\begin{Lem} 
\label{LemS}	
Fix an integer~$p\geq 2$ and~$f\in\mathcal{F}$. Let~$(\thetab_n)$ be a sequence in~$\mathcal{S}^{p-1}$ and~$(\kappa_n)$ be a positive real sequence that diverges to infinity. Then, under~${\rm P}\n_{\thetab_n,\kappa_n,f}$, 
$$
\frac{1}{n\tilde{e}_{n2}^{1/4}} 
\sum_{i=1}^n
\Xb_{ni}\Xb_{ni}\pr
=
\frac{1}{n\tilde{e}_{n2}^{1/4}} 
\sum_{i=1}^n
u_{ni}^2
\thetab_n\thetab_n\pr
+o_{\rm P}(1)
$$
as~$n\to\infty$, where~$u_{ni}=\Xb_{ni}'\thetab_n$ refers to the tangent-normal decomposition of~$\Xb_{ni}$ with respect to~$\thetab_n$. 
\end{Lem}

{\sc Proof of Lemma~\ref{LemS}.}
Using the tangent-normal decomposition of~$\Xb_{ni}$ with respect to~$\thetab_n$, write
\begin{eqnarray*}
\frac{1}{n\tilde{e}_{n2}^{1/4}} 
\sum_{i=1}^n
\Xb_{ni}\Xb_{ni}\pr
&\!\!\!=\!\!\!&
\frac{1}{n\tilde{e}_{n2}^{1/4}} 
\sum_{i=1}^n
(u_{ni}\thetab_n+v_{ni}\Sb_{ni})
(u_{ni}\thetab_n+v_{ni}\Sb_{ni})\pr
\\[2mm]
&\!\!\!=\!\!\!&
\Bigg(
\frac{1}{n\tilde{e}_{n2}^{1/4}} 
\sum_{i=1}^n
u_{ni}^2
\thetab_n\thetab_n\pr
\Bigg)
+
\Tb_{n1}
+
\Tb_{n2}
+
\Tb_{n3}
,
\end{eqnarray*}
where
$$
\Tb_{n1}
:=
\frac{1}{n\tilde{e}_{n2}^{1/4}} 
\sum_{i=1}^n
u_{ni}v_{ni} 
(\thetab_n\Sb_{ni}'+\Sb_{ni}\thetab_n')
,
\qquad
\Tb_{n2}
:=
\frac{1-e_{n2}}{(p-1)\tilde{e}_{n2}^{1/4}} (\mathbf{I}_p-\thetab_n\thetab_n')
$$
and
$$
\Tb_{n3}
:=
\frac{1}{n\tilde{e}_{n2}^{1/4}} 
\sum_{i=1}^n
\bigg(
v_{ni}^2\Sb_{ni}\Sb_{ni}\pr
-
\frac{1-e_{n2}}{p-1} (\mathbf{I}_p-\thetab_n\thetab_n')
\bigg)
.
$$
Applying the Cauchy-Schwarz inequality and using~(\ref{sjffl}) yields
$$
{\rm E}
\bigg[
\bigg\|
\frac{1}{n\tilde{e}_{n2}^{1/4}} 
\sum_{i=1}^n
u_{ni} v_{ni} \Sb_{ni}
\bigg\|^2
\bigg]
=
\frac{{\rm E}[u_{n1}^2 v_{n1}^2]}{n\sqrt{\tilde{e}_{n2}}} 
\leq
\frac{\sqrt{{\rm E}[v_{n1}^4]}}{n \sqrt{\tilde{e}_{n2}}} 
=
o(1)
,
$$
which implies that~$\Tb_{n1}$ converges to zero in probability. Using~(\ref{sjffl}) along with the fact~$\tilde{e}_{n2}=o(1)$ (Theorem~\ref{Theorexpansions}), we obtain that~$\Tb_{n2}=o(1)$. Now, denoting as~${\rm vec}$ the operator that stacks the columns of a matrix on top of each other and using the identity~${\rm E}[\Sb_{n1}\Sb_{n1}\pr]=(\mathbf{I}_p-\thetab_n\thetab_n')/(p-1)$, we obtain
\begin{eqnarray*}
\lefteqn{
{\rm E}
\big[
\big\|
{\rm vec}\, \Tb_{n3}
\big\|^2
\big]
=
\frac{1}{n\sqrt{\tilde{e}_{n2}}} 
{\rm E}
\bigg[
\bigg\|
{\rm vec}
\bigg(
v_{n1}^2\Sb_{n1}\Sb_{n1}\pr
-
\frac{1-e_{n2}}{p-1} (\mathbf{I}_p-\thetab_n\thetab_n')
\bigg)
\bigg\|^2
\bigg]
}
\\[2mm]
& & 
\hspace{-5mm} 
=
\frac{1}{n\sqrt{\tilde{e}_{n2}}} 
{\rm Tr}
\bigg[
{\rm E}
\bigg[
\bigg(
v_{n1}^2\Sb_{n1}\Sb_{n1}\pr
-
\frac{1-e_{n2}}{p-1} (\mathbf{I}_p-\thetab_n\thetab_n')
\bigg)^2
\bigg]
\bigg]
\\[2mm]
& & 
\hspace{-5mm} 
=
\frac{1}{n\sqrt{\tilde{e}_{n2}}} 
{\rm Tr}
\bigg[
{\rm E}
\bigg[
v_{n1}^4\Sb_{n1}\Sb_{n1}\pr
\!
-
\!
\frac{2(1-e_{n2})}{p-1} v_{n1}^2
\Sb_{n1}\Sb_{n1}\pr
\!
+
\!
\frac{(1-e_{n2})^2}{(p-1)^2} (\mathbf{I}_p-\thetab_n\thetab_n')
\bigg]
\bigg]
\\[2mm]
& & 
\hspace{-5mm} 
=
\frac{1}{n\sqrt{\tilde{e}_{n2}}} 
{\rm Tr}
\bigg[
 \frac{{\rm E}[v_{n1}^4]}{p-1} 
 (\mathbf{I}_p-\thetab_n\thetab_n')
-
\frac{(1-e_{n2})^2}{(p-1)^2} (\mathbf{I}_p-\thetab_n\thetab_n')
\bigg]
\\[2mm]
& & 
\hspace{-5mm} 
=
\frac{1}{n\sqrt{\tilde{e}_{n2}}} 
\bigg(
{\rm E}[v_{n1}^4]
-
\frac{(1-e_{n2})^2}{p-1}
\bigg)
=
o(1)
.
\end{eqnarray*}
Using again~(\ref{sjffl}) along with the fact~$\tilde{e}_{n2}=o(1)$ thus  shows that~$\Tb_{n3}$ converges to zero in probability, which establishes the result. 
\cqfd
\vspace{3mm}


{\sc Proof of Theorem~\ref{Theorsphericalmean2}.}
By using Theorem~\ref{Theorexpansions}(i), it follows from Theorem~\ref{Theorsphericalmean} that
$$
\frac{\sqrt{n(p-1)}(\hat{\thetab}_n- \thetab)}{\sqrt{1-e_{n2}}}
\stackrel{\mathcal{D}}{\to}
\mathcal{N}
\big(
{\bf 0}
,
{\bf I}_p- \thetab \thetab\pr
\big)
\quad\textrm{ and } \quad
\frac{n(p-1)\big(1-(\thetab\pr\hat{\thetab}_n)^2\big)}{1-e_{n2}}  
\stackrel{\mathcal{D}}{\to}
\chi^2_{p-1}
$$ 
and from Theorem~\ref{TheorsphericalCap} that
$$
\frac{2n(p-1)(1-\thetab\pr\hat{\thetab}_n)}{1-e_{n2}}
\stackrel{\mathcal{D}}{\to}
\chi^2_{p-1},
$$
as~$n\to\infty$ under~${\rm P}\n_{\thetab,\kappa_n,f}$ (in this proof, all stochastic convergences are under this sequence of hypotheses). Therefore, it is sufficient to show that 
\begin{eqnarray*}
\frac{1-\hat{e}_{n2}}{1-e_{n2}} 
&\!\!\!=\!\!\!&
\frac{1-\frac{1}{n} \sum_{i=1}^n (\Xb_{ni}\pr\thetab)^2}{1-e_{n2}} 
-
\frac{\tilde{e}_{n2}^{1/2}}{1-e_{n2}}
\times
\frac{\frac{1}{n} \sum_{i=1}^n( (\Xb_{ni}\pr\hat{\thetab}_n)^2-(\Xb_{ni}\pr\thetab)^2) }{\tilde{e}_{n2}^{1/2}}
\\[2mm]
&\!\!\!=:\!\!\!&
Y_{n1}
-
\frac{\tilde{e}_{n2}^{1/2}}{1-e_{n2}}
\,
Y_{n2}
=
1+o_{\rm P}(1)
.
\end{eqnarray*}
Since Theorem~\ref{Theorexpansions} implies that
$$ 
{\rm Var}[Y_{n1}]
=
\frac{{\rm Var}[u_{n1}^2]}{n(1-e_{n2})^2} 
=
\frac{{\rm Var}[v_{n1}^2]}{n(1-e_{n2})^2}
\leq
\frac{{\rm E}[v_{n1}^4]}{n(1-e_{n2})^2}
=
O\Big(\frac{1}{n}\Big)
,
$$
we have 
$
{\rm E}[(Y_{n1}-1)^2]
=
\big({\rm E}[Y_{n1}]-1\big)^2+{\rm Var}[Y_{n1}]
=
{\rm Var}[Y_{n1}]
=
o(1)
,
$
so that~$Y_{n1}= 1+ o_{\rm P}(1)$. Since the same theorem also implies that~$\tilde{e}_{n2}^{1/2}/(1-e_{n2})=O(1)$, it is sufficient to prove that~$Y_{n2}=o_{\rm P}(1)$. 

To do so, write 
\begin{equation}
	\label{Y2write}
Y_{n2}
=
\frac{1}{\sqrt{n}}
(\hat{\thetab}_n+\thetab)'
\Bigg[
\frac{1}{n\tilde{e}_{n2}^{1/4}} 
\sum_{i=1}^n
\Xb_{ni}\Xb_{ni}\pr
\Bigg]
\frac{\sqrt{n}(\hat{\thetab}_n-\thetab)}{\tilde{e}_{n2}^{1/4}}
\cdot
\end{equation}
Using Lemma~\ref{LemS} (with~$\thetab_n\equiv \thetab$) and~(\ref{1storderDelta}), we then obtain
\begin{eqnarray*}
Y_{n2}
&\!\!\!=\!\!\!&
\frac{1}{\sqrt{n}}
(\hat{\thetab}_n+\thetab)'
\Bigg[
\frac{1}{n\tilde{e}_{n2}^{1/4}} 
\sum_{i=1}^n
u_{ni}^2
\thetab\thetab\pr
+o_{\rm P}(1)
\Bigg]
\\[2mm]
& &
\hspace{33mm} 
\times
\Bigg(
\frac{\sqrt{n}}{{\tilde{e}_{n2}^{1/4}}} ({\bf I}_p - \thetab \thetab\pr)(\bar{\Xb}_n- e_{n1} \thetab)
+
O_{\rm P}\bigg(\frac{{\tilde{e}_{n2}^{1/4}}}{\sqrt{n}}\bigg)
\Bigg)
\\[2mm]
&\!\!\!=\!\!\!&
\frac{1}{\sqrt{n}}
(\hat{\thetab}_n+\thetab)'
\Bigg[
\frac{1}{n\tilde{e}_{n2}^{1/4}} 
\sum_{i=1}^n
u_{ni}^2
\thetab\thetab\pr
+o_{\rm P}(1)
\Bigg]
O_{\rm P}\bigg(\frac{{\tilde{e}_{n2}^{1/4}}}{\sqrt{n}}\bigg)
+
o_{\rm P}(1)
,
\end{eqnarray*}
where we used~(\ref{sdf}). Since~$n^{-1} \sum_{i=1}^n
u_{ni}^2\leq 1$ almost surely, we conclude that~$Y_{n2}$ is $o_{\rm P}(1)$, which establishes the result. 
\cqfd
\vspace{3mm}


\subsection{Proofs of Lemma~\ref{Lemtest} and Theorem~\ref{Theortest}} 
\label{SupSec5}
\textrm{}
\vspace{2mm}

\noindent {\sc Proof of Lemma~\ref{Lemtest}.}
We start with~$\bar{\Xb}_{n}\pr\thetab_{0}$. Since~$
\bar{\Xb}_{n}\pr\thetab_{0}
=
\bar{\Xb}_{n}\pr\thetab_{n}
-
\nu_n \bar{\Xb}_{n}\pr\taub_n
$,
we have 
\begin{eqnarray*}
	{\rm E}[(\bar{\Xb}_n\pr\thetab_{0}-1)^2]
&\!\!\!\leq \!\!\!&
3 {\rm E}[(\bar{\Xb}_{n}\pr\thetab_{n}-e_{n1})^2]
+
3\nu_n^2 {\rm E}[(\bar{\Xb}_{n}\pr\taub_n)^2]
+
3 (e_{n1}-1)^2 
\\[2mm]
&\!\!\!= \!\!\!&
3 {\rm E}[(\bar{\Xb}_{n}\pr\thetab_{n}-e_{n1})^2] 
+
o(1)
,
\end{eqnarray*}
where we used the facts that~$\|\bar{\Xb}_{n}\|\leq 1$ almost surely and that~$e_{n1}=1+o(1)$. Since the tangent-normal decomposition with respect to~$\thetab_n$ further entails that
$$
{\rm E}[(\bar{\Xb}_{n}\pr\thetab_{n}-e_{n1})^2]
= 
{\rm E}\bigg[ \bigg( \frac{1}{n} \sum_{i=1}^n (u_{ni}-e_{n1}) \bigg)^2\bigg]
=
{\rm Var}\bigg[\frac{1}{n}\sum_{i=1}^n u_{ni}\bigg]
=
\frac{\tilde{e}_{n2}}{n}
=
o(1) 
,
$$
we conclude that~$\bar{\Xb}_n\pr\thetab_{0}$ converges to one in quadratic mean, hence also in probability.  

We turn to~$R_n$, which we decompose as 
$$
R_n
=
\frac{1-\frac{1}{n} \sum_{i=1}^n(\Xb_{ni}\pr\thetab_n-\nu_n\Xb_{ni}'\taub_n)^2}{\sqrt{2(p-1)} \tilde{e}_{n2}^{1/2}}
=
R_{n1}
+
R_{n2}
-
R_{n3}
,
$$
with 
$$
R_{n1}
:=
\frac{1-\frac{1}{n} \sum_{i=1}^n(\Xb_{ni}\pr\thetab_n)^2}{\sqrt{2(p-1)} \tilde{e}_{n2}^{1/2}}
,
\qquad
R_{n2}
:=
\frac{\sqrt{2}\nu_n}{\sqrt{p-1}\, n\tilde{e}_{n2}^{1/2}} 
\thetab_n' \bigg(  \sum_{i=1}^n \Xb_{ni} \Xb_{ni}' \bigg)\taub_n
$$
and
$$
R_{n3}
:=
\frac{\nu_n^2}{\sqrt{2(p-1)} n\tilde{e}_{n2}^{1/2}} \sum_{i=1}^n (\Xb_{ni}'\taub_n)^2
.
$$
Since~(\ref{sjffl}) entails that
$$
{\rm E}[R_{n1}]
=
\frac{1-e_{n2}}{\sqrt{2(p-1)}\tilde{e}_{n2}^{1/2}}
=
1+o(1)
$$
and
$$ 
{\rm Var}[R_{n1}]
=
\frac{{\rm Var}[u_{n1}^2]}{2(p-1)n\tilde{e}_{n2}} 
=
\frac{{\rm Var}[v_{n1}^2]}{2(p-1)n\tilde{e}_{n2}} 
\leq
\frac{{\rm E}[v_{n1}^4]}{2(p-1)n\tilde{e}_{n2}} 
=
o(1)
,
$$
we have that~$R_{n1}$ converges to one in quadratic mean, hence also in probability. As for~$R_{n2}$, Lemma~\ref{LemS} and Theorem~\ref{Theorexpansions}(ii) yield
$$
R_{n2}
=
\frac{\sqrt{2}\nu_n}{\sqrt{p-1} \tilde{e}_{n2}^{1/4}} 
\thetab_n' \bigg( 
\frac{U_n}{\tilde{e}_{n2}^{1/4}} 
\thetab_n\thetab_n\pr
+
o_{\rm P}(1)
 \bigg)\taub_n
=
\frac{\sqrt{2}\nu_n (\thetab_n\pr\taub_n) U_n}{\sqrt{p-1} \, \tilde{e}_{n2}^{1/2}} 
+
o_{\rm P}(1)
,
$$
where we let~$U_n:=(1/n) \sum_{i=1}^nu_{ni}^2$. Since~$\thetab_n=\thetab+\nu_n\taub_n$ is a unit $p$-vector, we have~$\thetab'\taub_n=-\nu_n\|\taub_n\|^2/2$, which yields~$\thetab_n\pr\taub_n=(\thetab+\nu_n\taub_n)\pr\taub_n=\nu_n\|\taub_n\|^2/2$. Thus, using Theorem~\ref{Theorexpansions}(ii) and the fact that~$U_n\leq 1$ almost surely, we obtain
$$
R_{n2}
=
\frac{\nu_n^2\|\taub_n\|^2 U_n}{\sqrt{2(p-1)} \, \tilde{e}_{n2}^{1/2}} 
+
o_{\rm P}(1)
=
o_{\rm P}(1)
.
$$
Finally, since~$(\Xb_{ni}'\taub_n)^2\leq \|\taub_n\|^2$ almost surely, Theorem~\ref{Theorexpansions}(ii) also entails that~$
R_{n3}
=
o_{\rm P}(1)
.
$
Therefore,~$R_n=1+o_{\rm P}(1)$, as was to be proved. 
\cqfd
\vspace{3mm}

{\sc Proof of Theorem~\ref{Theortest}.}
Since Part~(i) of the result is actually a particular case of Part~(ii), we only prove the latter. Accordingly, all stochastic convergences in this proof will be as~$n\to\infty$ under~${\rm P}\n_{\thetab_n,\kappa_n,f}$, with~$\thetab_n=\thetab_0+\nu_n\taub_n$, $\nu_n:=1/\sqrt{n\kappa_n\varphi_f(\kappa_n)}$ and~$\taub_n\to\taub$. 
Consider then
\begin{eqnarray*}
\Tb^W_n
&\!\!\!:=\!\!\!&
\frac{(p-1)^{1/4}\sqrt{n}(\bar{\Xb}_n-e_{n1}\thetab_0)}{2^{1/4}{\tilde{e}_{n2}^{1/4}}}
\\[2mm]
&\!\!\!=\!\!\!&
\frac{(p-1)^{1/4}\sqrt{n}(\bar{\Xb}_n-e_{n1}\thetab_n)}{2^{1/4}{\tilde{e}_{n2}^{1/4}}}
+
\frac{(p-1)^{1/4}}{2^{1/4}\tilde{e}_{n2}^{1/4}\sqrt{\kappa_n\varphi_f(\kappa_n)}}
e_{n1}\taub_n
\\[2mm]
&\!\!\!=\!\!\!&
\frac{(p-1)^{1/4}\sqrt{n}(\bar{\Xb}_n-e_{n1}\thetab_n)}{2^{1/4}{\tilde{e}_{n2}^{1/4}}}
+
\taub
+
o(1)
,
\end{eqnarray*}
where we used Theorem~\ref{Theorexpansions}(ii) and the fact that~$e_{n1}=1+o(1)$. Now, proceeding exactly as in the proof of Theorem~\ref{Theorsphericalmean}, it can be shown that
$$
\frac{\sqrt{n}(\bar{\Xb}_n- e_{n1} \thetab_n)}{{\tilde{e}_{n2}^{1/4}}}
\stackrel{\mathcal{D}}{\to}
\mathcal{N}
\bigg(
{\bf 0}
,
\frac{\sqrt{2}}{\sqrt{p-1}}( {\bf I}_p- \thetab_0 \thetab_0\pr) 
\bigg)
.
$$
It follows that 
$
\Tb^W_n
\stackrel{\mathcal{D}}{\to}
\mathcal{N}
\big(
\taub
,
 {\bf I}_p- \thetab_0 \thetab_0\pr
\big)
, 
$
so that Lemma~\ref{Lemtest} entails that
$$
W_n = \tilde{W}_n+o_{\rm P}(1)=
(\Tb^W_n)'( {\bf I}_p- \thetab_0 \thetab_0\pr)\Tb^W_n
+o_{\rm P}(1)
\stackrel{\mathcal{D}}{\to}
\chi^2_{p-1}
\big(
\|\taub\|^2
\big)
.
$$
Turning then to the Wald test, consider now 
$$
\Tb^S_n
:=
\frac{(p-1)^{1/4}\sqrt{n}(\hat{\thetab}_n- \thetab_0)}{2^{1/4}{\tilde{e}_{n2}^{1/4}}}
=
\frac{(p-1)^{1/4}\sqrt{n}(\hat{\thetab}_n- \thetab_n)}{2^{1/4}{\tilde{e}_{n2}^{1/4}}}
+
\taub
+
o(1)
.
$$
Since, under the sequence of hypotheses considered, Lemma~\ref{Lemseconddelta} implies that
$$
\frac{\sqrt{n}}{{\tilde{e}_{n2}^{1/4}}}
(\hat{\thetab}_n- \thetab_n)
=
\frac{\sqrt{n}}{{\tilde{e}_{n2}^{1/4}}}
 ({\bf I}_p - \thetab_0 \thetab_0\pr)
(\bar{\Xb}_n- e_{n1} \thetab_n)
+ 
o_{\rm P}(1)
,
$$
we have that
$
\Tb^S_n
=
({\bf I}- \thetab_0 \thetab_0\pr) \Tb^W_n
+
o_{\rm P}(1)
.
$
Using Lemma~\ref{Lemtest} again, this yields that~$S_n=\tilde{S}_n+o_{\rm P}(1)=(\Tb^S_n)'( {\bf I}_p- \thetab_0 \thetab_0\pr)
\Tb^S_n+o_{\rm P}(1)=(\Tb^W_n)'( {\bf I}_p- \thetab_0 \thetab_0\pr)
\Tb^W_n+o_{\rm P}(1)=W_n+o_{\rm P}(1)$, which establishes the result. 
\cqfd
\vspace{3mm}


\subsection{Proofs of Proposition~\ref{PropCondLANexpb} and Theorem~\ref{TheorLAN}}
\label{SupSec6}

The proof of Proposition~\ref{PropCondLANexpb} requires the following result.

\begin{Lem}
\label{chipo}
Fix~$b>0$. Then there exists~$C_b$ such that for any~$x,y\in\R$ with~$x,y>0$, one has~$|y^b-x^b-b(y-x)x^{b-1}|\leq C_b (y-x)^2 (|x|^{b-2}+|y|^{b-2})$.  
\end{Lem}

{\sc Proof}. 
Since~$x,y>0$, the mapping~$z\mapsto z^b$ is continuous on the interval with end points~$x$ and~$y$, and it is differentiable on the interior of this interval. The mean value theorem then yields that, for some~$c$ between~$x$ and~$y$,
\begin{eqnarray*}
	\lefteqn{
|y^b-x^b-b(y-x)x^{b-1}|
=
|b(b-1)(y-x)^2c^{b-2}/2|
}
\\[2mm]
& & 
\hspace{-7mm} 
\leq
|b(b-1)| (y-x)^2 \max(|x|^{b-2},|y|^{b-2})
\leq
|b(b-1)| (y-x)^2 (|x|^{b-2}+|y|^{b-2})
,
\end{eqnarray*}
which establishes the result. 
\cqfd
\vspace{3mm}

{\sc Proof of Proposition~\ref{PropCondLANexpb}.}
With~$f(z)=\exp(z^b)$, 
\begin{eqnarray}
	\lefteqn{
\hspace{-1mm} 
\kappa^{(p+1)/2} (\varphi_f(\kappa))^{(p-3)/2}
\frac{1}{f(\kappa)}
\int_{-1}^1
\big(\varphi_f(\kappa s)-\varphi_f(\kappa)\big)^2
(1-s^2)^{(p-3)/2}
f(\kappa s)
\,
ds
}
\nonumber
\\[2mm]
& & 
\hspace{2mm} 
=
b^{(p+1)/2} \kappa^{b(p+1)/2}
\int_{-1}^1
(s^{b-1}-1)^2
(1-s^2)^{(p-3)/2}
e^{\kappa^b (s^b-1)}
\,ds
\nonumber
\\[2mm]
& & 
\hspace{2mm} 
=
b^{(p+1)/2} \kappa^{b(p+1)/2}
\int_{0}^2
((1-r)^{b-1}-1)^2
(r(2-r))^{(p-3)/2}
e^{\kappa^b\{(1-r)^b-1\}}
\,dr
,
\label{Freddy}
\end{eqnarray}
where we let~$r=1-s$. Let~$\delta=1$ if~$b\in(0,1)$ and~$0$ otherwise. Then, by using Lemma~\ref{lemtec}(i)--(ii) and the fact that there exists some constant~$C$ such that~$|(1-r)^{b-1}-1|\leq C r |1-r|^{\delta(b-1)}$ for any~$r\in[0,2]$, we obtain that~(\ref{Freddy}) is upper-bounded by
\begin{eqnarray*}
	\lefteqn{
\hspace{-2mm} 
C
\kappa^{b(p+1)/2}
\!
\int_{0}^{1/2}
\!
r^{(p+1)/2}
e^{-K\kappa^b r}
\,dr
}
\\[2mm]
& & 
\hspace{13mm} 
+
C
\kappa^{b(p+1)/2}
e^{-K\kappa^b/2}
\!
\int_{1/2}^{2}
(1-r)^{2(b-1)}
(2-r)^{(p-3)/2}
\,dr
\\[2mm]
& & 
\hspace{-3mm} 
\leq
C
\frac{1}{\kappa^{b}}
\int_{0}^{K\kappa^b/2}
z^{(p+1)/2}
e^{-z}
\,dz
+
o(1)
=
O\Big(\frac{1}{\kappa^{b}}\Big)
+
o(1)
=
o(1)
.
\end{eqnarray*}

We may therefore focus on~(\ref{assudeath}). Fix then a positive sequence~$(\kappa_n)$ diverging to infinity (which, for~$b\in(\frac{1}{2},1)$, is assumed to satisfy the assumption stated in the proposition), a bounded positive sequence~$(t_n)$, and consider the quantities~$h^\pm_n(s,w)$ appearing in~(\ref{assudeath}). First note that 
$$
\frac{\kappa_n s+h^\pm_n(s,w)}{\kappa_n s}
=
1-{\textstyle{\frac{1}{2}}}t_n^2 \nu_n^2
\pm \frac{c_n t_n \nu_n (1-s^2)^{1/2} w^{1/2}}{s}
,
$$
so that, for~$n$ large enough,
$$
\bigg|1-\frac{\kappa_n s+h^\pm_n(s,w)}{\kappa_n s}\bigg|
\leq
\frac{1}{4}
+
\frac{M\nu_n}{|s|}
,
$$
where we let~$M:=\sup_n t_n$. Hence, for~$s\notin \mathcal{I}_n :=(-4M\nu_n,4M\nu_n)$ and~$n$ large enough,  
\begin{equation}
	\label{espieglerie}
\frac{\kappa_n s+h^\pm_n(s,w)}{\kappa_n s}
\in 
\bigg[ \frac{1}{2},\frac{3}{2} \bigg]
.
\end{equation}
For~$n$ large enough, we then have
\begin{eqnarray*}
	\lefteqn{
\hspace{-1mm} 
\frac{n(\kappa_n \varphi_f(\kappa_n))^{(p-1)/2}}{f(\kappa_n)}
}
\\[2mm]
& & 
\hspace{3mm} 
\times
\int_{-1}^1
\int_{0}^1
\Big|
\log f(\kappa_n s+h^\pm_n(s,w))
-
\log f(\kappa_n s)
-
h^\pm_n(s,w) \varphi_f(\kappa_n s)
\Big| 
\\[2mm]
& & 
\hspace{33mm} 
\times f(\kappa_n s) (1-s^2)^{(p-3)/2}w^{-1/2}
\,
dG_p(w)
ds
\\[2mm]
& & 
\hspace{-3mm} 
=
n\kappa_n^{b(p-1)/2}
\int_{-1}^1
\int_{0}^1
\big|
(\kappa_n s+h^\pm_n(s,w))^b
-
(\kappa_n s)^b
-
b h^\pm_n(s,w) (\kappa_n s)^{b-1}
\big| 
\\[2mm]
& & 
\hspace{33mm} 
\times 
(1-s^2)^{(p-3)/2}
e^{\kappa_n^b(s^b-1)}
\,
dG_p(w)
ds
\\[2mm]
& & 
\hspace{-3mm} 
\leq T_{n1}+T_{n2}
,
\end{eqnarray*}
where, with~$\mathcal{I}_n^c:=[-1,1]\setminus \mathcal{I}_n$, we let
\begin{eqnarray*}
	\lefteqn{
\hspace{-16mm} 
T_{n1}
:=
n\kappa_n^{b(p-1)/2}
\int_{\mathcal{I}_n^c}
\int_{0}^1
\big|
(\kappa_n s+h^\pm_n(s,w))^b
-
(\kappa_n s)^b
-
b h^\pm_n(s,w) (\kappa_n s)^{b-1}
\big| 
}
\\[2mm]
& &
\hspace{25mm} 
\times 
(1-s^2)^{(p-3)/2}
e^{\kappa_n^b(s^b-1)}
\,
dG_p(w)
ds
\end{eqnarray*}
and
\begin{eqnarray*}
	\lefteqn{
\hspace{-1mm} 
T_{n2}
:=
n\kappa_n^{b(p-1)/2} e^{-(1-C_{\varepsilon})\kappa_n^b}
}
\\[2mm]
& &
\hspace{13mm} 
\times 
\int_{\mathcal{I}_n}
\int_{0}^{1}
\big|
(\kappa_n s+h^\pm_n(s,w))^b
-
(\kappa_n s)^b
-
b h^\pm_n(s,w) (\kappa_n s)^{b-1}
\big| 
\\[2mm]
& &
\hspace{43mm} 
\times 
(1-s^2)^{(p-3)/2}
\,
dG_p(w)
ds
;
\end{eqnarray*}
here, $C_{\varepsilon}:=1/2$ if~$b\geq 1$ and  $C_{\varepsilon}:=\varepsilon/4$ if~$b\in(0,1)$, where~$\varepsilon>0$ is as in the statement of the proposition. 

Let us first consider~$T_{1n}$. It directly follows from~(\ref{espieglerie}) that, for~$n$ large enough,~$\kappa_n s+h^\pm_n(s,w)$ and~$\kappa_n s$ share the same sign in the integrand of~$T_{1n}$. Consequently, using Lemma~\ref{chipo} then~(\ref{espieglerie}) yields 
\begin{eqnarray*}
\lefteqn{
\big|
(\kappa_n s+h^\pm_n(s,w))^b
-
(\kappa_n s)^b
-
b h^\pm_n(s,w) (\kappa_n s)^{b-1}
\big| 
}
\\[2mm]
& & 
\hspace{13mm} 
\leq
C
(h^\pm_n(s,w))^2
\Big(
|\kappa_n s+h^\pm_n(s,w)|^{b-2}
+
|\kappa_n s|^{b-2}
\Big) 
\\[2mm]
& & 
\hspace{13mm} 
\leq 
C (\kappa_n^2 \nu_n^4 + \kappa_n^2 \nu_n^2 (1-s^2) ) \kappa_n^{b-2} 
\\[2mm]
& & 
\hspace{13mm} 
\leq 
C
 \big( \kappa_n^{-b} +  1-s^2 \big)  
\end{eqnarray*}
for~$n$ large enough. Therefore, by using again Lemma~\ref{lemtec}(i)--(ii), we obtain that, still for~$n$ large enough, 
\begin{eqnarray*}
T_{n1}
&\!\!\leq\!\!&
C\kappa_n^{b(p-1)/2} 
\int_{-1}^1
 \big( \kappa_n^{-b} +  1-s^2 \big)  
e^{\kappa_n^b(s^b-1)}
(1-s^2)^{(p-3)/2}
\,
ds
\\[2mm]
&\!\!\leq\!\!&
C\kappa_n^{b(p-1)/2} 
\int_{0}^{2}
 \big( \kappa_n^{-b} +  r(2-r) \big)  
(r(2-r))^{(p-3)/2}
e^{\kappa_n^b((1-r)^b-1)}
\,
dr
\\[2mm]
&\!\!\leq\!\!&
C\kappa_n^{b(p-1)/2} 
\int_{0}^{2}
 \big( \kappa_n^{-b} +  r(2-r) \big)  
(r(2-r))^{(p-3)/2}
e^{-K\kappa_n^b r}
\,
dr
.
\end{eqnarray*}
Letting~$z=K\kappa_n^b r$, we obtain
\begin{eqnarray*}
T_{n1}
&\!\!\!\!\leq\!\!\!\!&
C\kappa_n^{b(p-3)/2}
\!
\int_{0}^{2K\kappa_n^b}
\!
 \big(\kappa_n^{-b} + {\textstyle{\frac{z}{K\kappa_n^b}\big(2-\frac{z}{K\kappa_n^b}\big)}} \big)  
\big({\textstyle{\frac{z}{K\kappa_n^b}\big(2-\frac{z}{K\kappa_n^b}\big)}}\big)^{(p-3)/2}
e^{-z}
\,
dz
\\[2mm]
&\!\!\!\!\leq\!\!\!\!&
C\kappa_n^{-b}
\int_{0}^{2K\kappa_n^b}
 \big(1 + {\textstyle{\frac{z}{K}\big(2-\frac{z}{K\kappa_n^b}\big)}} \big)  
\big({\textstyle{\frac{z}{Kb}\big(2-\frac{z}{K\kappa_n^b}\big)}}\big)^{(p-3)/2}
e^{-z}
\,
dz,
\end{eqnarray*}
which shows using Lemma \ref{lemDCT} that~$T_{n1}$ is~$O(\kappa_n^{-b})$, hence~$o(1)$. 

Turning to~$T_{n2}$, we have 
\begin{eqnarray*}
\lefteqn{
\big|
(\kappa_n s+h^\pm_n(s,w))^b
-
(\kappa_n s)^b
-
b h^\pm_n(s,w) (\kappa_n s)^{b-1}
\big| 
}
\\[2mm]
& & 
\hspace{13mm} 
\leq
|\kappa_n s+h^\pm_n(s,w)|^b
+
|\kappa_n s|^b
+
b |h^\pm_n(s,w)| |\kappa_n s|^{b-1}
\\[2mm]
& & 
\hspace{13mm} 
\leq 
C
\kappa_n^b
(
|s|^b
+
\nu_n^b
)
+
\kappa_n^b
|s|^b
+
C
\kappa_n^b
(|s|+
\nu_n)
|s|^{b-1}
)
\\[2mm]
& & 
\hspace{13mm} 
\leq
C
\kappa_n^b
 (|s|^{b}+\nu_n|s|^{b-1}+\nu_n^b)
,
\end{eqnarray*}
which yields
\begin{eqnarray*}
T_{n2}
&\!\!\!=\!\!\!&
C n\kappa_n^{b(p+1)/2} e^{-(1-C_{\varepsilon})\kappa_n^b}\int_{-4M\nu_n}^{4M\nu_n}
(|s|^{b}+\nu_n|s|^{b-1}+\nu_n^b)
\,
ds
\\[2mm]
&\!\!\!=\!\!\!&
Cn\kappa_n^{b(p+1)/2} \nu_n^{b+1} e^{-(1-C_{\varepsilon})\kappa_n^b}
\\[3mm]
&\!\!\!=\!\!\!&
Cn^{(1-b)/2} \kappa_n^{b(p-b)/2} e^{-(1-C_{\varepsilon})\kappa_n^b}
.
\end{eqnarray*}
Consequently, if~$b\geq 1$, then~$T_{n2}$ is~$o(1)$, as was to be shown. Focus then on the case~$b\in(\frac{1}{2},1)$. By assumption, for~$n$ large enough,
$$
\frac{\log n}{\kappa_n^b}\leq \frac{2-\varepsilon}{1-b},
\quad
\textrm{ or equivalently},
\quad
n\leq e^{(\frac{2-\varepsilon}{1-b}) \kappa_n^b}
,
$$
which yields
\begin{eqnarray*}
	\lefteqn{
T_{n2}
\leq
Cn^{(1-b)/2} \kappa_n^{b(p-b)/2} e^{-(1-(\varepsilon/4))\kappa_n^b}
}
\\[2mm]
& & 
\hspace{13mm} 
\leq
C \kappa_n^{b(p-b)/2} e^{\{(2-\varepsilon)/2-(1-\varepsilon/4)\}\kappa_n^b}
=
C \kappa_n^{b(p-b)/2} e^{-\varepsilon\kappa_n^b/4}
,
\end{eqnarray*}
so that~$T_{n2}=o(1)$. The result follows. 
\cqfd
\vspace{3mm}

{\sc Proof of Theorem~\ref{TheorLAN}.}
Write
\begin{eqnarray*}
\log \frac{d{\rm P}\n_{\thetab+\nu_n\taub_n,\kappa_n,f}}{d{\rm P}\n_{\thetab,\kappa_n,f}} 
&\!\!\!=\!\!\!&
\sum_{i=1}^n 
\big(
\log f(\kappa_n u_{ni}+\kappa_n \nu_n\taub_n'\Xb_{ni})
-
\log f(\kappa_n u_{ni})
\big)
\\[2mm]
&\!\!\!=\!\!\!&
 L_{n1}+L_{n2}+L_{n3},
\end{eqnarray*}
with
$$
L_{n1}
:=
\nu_n^{-1} \taub_n' \bar{\Xb}_n
=
\sqrt{n\kappa_n\varphi_f(\kappa_n)}\, \taub_n' \bar{\Xb}_n
,
$$
$$
L_{n2}
:=
\frac{\sqrt{\kappa_n}}{\sqrt{n\varphi_f(\kappa_n)}}
\sum_{i=1}^n 
\big(
\varphi_f(\kappa_n u_{ni})-\varphi_f(\kappa_n)
\big)
 \taub_n' \Xb_{ni}
,
$$
and
$$
L_{n3}
:=
\sum_{i=1}^n 
\bigg(\!
\log f(\kappa_n u_{ni}+\kappa_n \nu_n\taub_n'\Xb_{ni})
-
\log f(\kappa_n u_{ni})
-
\kappa_n \nu_n
  \varphi_f(\kappa_n u_{ni}) \taub_n' \Xb_{ni}
\bigg)
.
$$

Using the identity~$\taub_n'\thetab=-\frac{1}{2}\nu_n\|\taub\|^2$ and Lemma~\ref{Lemtest}, we readily obtain
\begin{eqnarray*}
\lefteqn{
\hspace{-3mm} 
L_{n1}
=
\nu_n^{-1} \taub_n' (\mathbf{I}_p-\thetab\thetab')\bar{\Xb}_n
+
\nu_n^{-1} (\taub_n'\thetab) (\thetab' \bar{\Xb}_n)
=
\taub_n' \Deltab\n_{\thetab,f}
-
\frac{1}{2} 
\|\taub_n\|^2 (\thetab' \bar{\Xb}_n)
}
\\[2mm]
& & 
\hspace{3mm} 
=
\taub_n' \Deltab\n_{\thetab,f}
-
\frac{1}{2} 
\|\taub_n\|^2 
+
o_{\rm P}(1)
=
\taub_n' 
\Deltab\n_{\thetab,f}
-
\frac{1}{2} 
\taub_n\pr 
\Gamb_{\thetab}
\taub_n
+
o_{\rm P}(1)
,
\end{eqnarray*}
so that we only need to show that both~$L_{n2}$ and~$L_{n3}$ are~$o_{\rm P}(1)$.  

We start with~$L_{n2}$. Using the tangent-normal decomposition of~$\Xb_{ni}$ with respect to~$\thetab$, write~$L_{n2}=L_{n2a}+\taub_n\pr {\bf L}_{n2b}$, where we let
$$
L_{n2a}
:=
\frac{\sqrt{\kappa_n}}{\sqrt{n\varphi_f(\kappa_n)}} 
\sum_{i=1}^n  
\big(\varphi_f(\kappa_n u_{ni})-\varphi_f(\kappa_n)\big)
u_{ni} (\taub_n\pr\thetab)
$$
and
$$
L_{n2b} 
:=
\frac{\sqrt{\kappa_n}}{\sqrt{n\varphi_f(\kappa_n)}} 
\sum_{i=1}^n  
\big(\varphi_f(\kappa_n u_{ni})-\varphi_f(\kappa_n)\big)
v_{ni} \Sb_{ni}
.
$$
We have
\begin{eqnarray*}
{\rm E}[L_{n2a}^2]
&\!\!\!\leq \!\!\!&
\frac{n\kappa_n}{\varphi_f(\kappa_n)} 
\,
{\rm E}\big[  \big(\varphi_f(\kappa_n u_{n1})-\varphi_f(\kappa_n)\big)^2 u_{n1}^2 \big] \,(\taub_n\pr\thetab)^2
\\[1mm]
&\!\!\!\leq \!\!\!&
\frac{n\kappa_n\nu_n^2}{\varphi_f(\kappa_n)} 
\|\taub_n\|^2
\,
{\rm E}\big[  \big(\varphi_f(\kappa_n u_{n1})-\varphi_f(\kappa_n)\big)^2 u_{n1}^2 \big] 
\\[1mm]
&\!\!\!\leq \!\!\!&
\frac{\|\taub_n\|^2}{\kappa_n\varphi_f(\kappa_n)} 
\,
\bigg(
\frac{\kappa_n}{\varphi_f(\kappa_n)} 
\,
{\rm E}\big[  \big(\varphi_f(\kappa_n u_{n1})-\varphi_f(\kappa_n)\big)^2 \big] 
\bigg)
\end{eqnarray*}
and
\begin{eqnarray*}
{\rm E}[\|{\bf L}_{n2b}\|^2]
&\!\!\!= \!\!\!&
\frac{\kappa_n}{\varphi_f(\kappa_n)} 
\,
{\rm E}\big[  \big(\varphi_f(\kappa_n u_{n1})-\varphi_f(\kappa_n)\big)^2 v_{n1}^2 \big] 
\\[1mm]
&\!\!\!\leq \!\!\!&
\frac{\kappa_n}{\varphi_f(\kappa_n)} 
\,
{\rm E}\big[  \big(\varphi_f(\kappa_n u_{n1})-\varphi_f(\kappa_n)\big)^2 \big] 
.
\end{eqnarray*}
Now, by using Lemma~\ref{LemExpansionConstant}(i) and the fact that~$f\in\mathcal{F}_{\rm LAN}(p,\kappa_n,\|\taub_n\|)$, we obtain
\begin{eqnarray*}
	\lefteqn{
\frac{\kappa_n}{\varphi_f(\kappa_n)} 
\,
{\rm E}\big[  \big(\varphi_f(\kappa_n s)-\varphi_f(\kappa_n)\big)^2 \big]
}
\\[2mm]
& & 
\hspace{-3mm} 
=
\frac{\kappa_n c_{p,\kappa_n,f}}{\varphi_f(\kappa_n)} 
\int_{-1}^1
\int_{\mathcal{S}^\perp_{\thetab}}
\big(\varphi_f(\kappa_n s)-\varphi_f(\kappa_n)\big)^2
(1-s^2)^{(p-3)/2}
f(\kappa_n s)
\,
d\sigma(\ub)
ds
\\[2mm]
& & 
\hspace{-3mm} 
=
O\Big(
\kappa_n^{(p+1)/2} (\varphi_f(\kappa_n))^{(p-3)/2}
\Big)
\!
\\[2mm]
& & 
\hspace{13mm} 
\times
\int_{-1}^1
\big(\varphi_f(\kappa_n s)-\varphi_f(\kappa_n)\big)^2
(1-s^2)^{(p-3)/2}
\frac{f(\kappa_n s)}{f(\kappa_n)}
\,
ds
=
o(1)
.
\end{eqnarray*}
Therefore, ${\rm E}[L_{n2a}^2]$ and ${\rm E}[\|{\bf L}_{n2b}\|^2]$ are~$o(1)$, which implies that~$L_{n2a}$ and~${\bf L}_{n2b}$, hence also~$L_{n2}$, are~$o_{\rm P}(1)$.

Let us turn to~$L_{n3}$. 
Since~$\taub_n'\Xb_{n1}=u_{n1}\taub_n'\thetab+v_{n1} \taub_n'\Sb_{n1}=-\frac{1}{2}\nu_n u_{n1}\|\taub_n\|^2+v_{n1} \taub_n'\Sb_{n1}$ and $
\|(\mathbf{I}_p-\thetab\thetab')\taub_n\|^2
=
\|\taub_n\|^2-(\thetab'\taub_n)^2
=
\|\taub_n\|^2-\frac{1}{4}\nu_n^2\|\taub_n\|^4
=
c_n^2\|\taub_n\|^2
$, rotation invariance yields that~${\rm E}[|L_{n3}|]$ is upper-bounded by
\begin{eqnarray*}
\lefteqn{
n
{\rm E}
\Big[
\big|
\log f(\kappa_n u_{n1}+\kappa_n \nu_n\taub_n'\Xb_{n1})
-
\log f(\kappa_n u_{n1})
-
 \kappa_n \nu_n
  \varphi_f(\kappa_n u_{n1}) \taub_n' \Xb_{n1}
\big|
\Big]
}
\\[2mm]
& & 
\hspace{-5mm}  
=
n
{\rm E}
\Big[
\big|
\!
\log f(\kappa_n u_{n1}-{\textstyle{\frac{1}{2}}}\kappa_n\nu_n^2 u_{n1}\|\taub_n\|^2+c_n \kappa_n\nu_n v_{n1} \|\taub_n\| U_{n1})
-
\log f(\kappa_n u_{n1}) 
\\[2mm]
& & 
\hspace{-2mm} 
-
\big(
-{\textstyle{\frac{1}{2}}}\kappa_n\nu_n^2\varphi_f(\kappa_n u_{n1}) u_{n1}\|\taub_n\|^2+c_n \kappa_n\nu_n \varphi_f(\kappa_n u_{n1}) v_{n1} \|\taub_n\| U_{n1}
\big)
\big| 
\Big] 
,
\end{eqnarray*}
where~$U_{n1}=\thetab_\perp'\Sb_n$, with~$\thetab_\perp$ an arbitrary unit vector orthogonal to~$\thetab$. Clearly, $U_{n1}$ is equal in distribution to any marginal of a random vector that is uniformly distributed over~$\mathcal{S}^{p-2}$. Therefore, $-U_{n1}\stackrel{\mathcal{D}}{=} U_{n1}$ and~$W:=U_1^2$ has the cumulative distribution function~$G_p$ in page~\pageref{pageGp}, so that conditioning with respect to the sign of~$U_{n1}$ yields that~${\rm E}[|L_{n3}|]$ is~$o(1)$  if and only if
\begin{eqnarray*}
	\lefteqn{
n
{\rm E}
\big[
\big|
\log f(\kappa_n u_{n1}-{\textstyle{\frac{1}{2}}}\kappa_n\nu_n^2 u_{n1}\|\taub_n\|^2\pm c_n \kappa_n\nu_n v_{n1} \|\taub_n\| W^{1/2})
-
\log f(\kappa_n u_{n1}) 
}
\\[2mm]
& & 
\hspace{-5mm} 
-
\big(
-{\textstyle{\frac{1}{2}}}\kappa_n\nu_n^2\varphi_f(\kappa_n u_{n1}) u_{n1}\|\taub_n\|^2\pm c_n \kappa_n\nu_n \varphi_f(\kappa_n u_{n1}) v_{n1} \|\taub_n\| W^{1/2}
\big)
\big|
\big]
=o(1)
.
\end{eqnarray*}
In view of Lemma~\ref{LemExpansionConstant}(i), this is the case if and only if
\begin{eqnarray*}
	\lefteqn{
	\hspace{-1mm} 
\frac{1}{f(\kappa_n)}
\int_{-1}^1
\int_{0}^1
\Big|
\log f(\kappa_n s-{\textstyle{\frac{1}{2}}}\kappa_n\nu_n^2 s\|\taub_n\|^2\pm c_n \kappa_n\nu_n (1-s^2)^{1/2} \|\taub_n\| w^{1/2})  
}
\\[2mm]
& & 
\hspace{-8mm} 
-
\!
 \log f(\kappa_n s)
\!
-
\!
\big(
\!
\!
-
\!
{\textstyle{\frac{1}{2}}}\kappa_n\nu_n^2\varphi_f(\kappa_n s) s\|\taub_n\|^2
\!
\pm
\!
c_n \kappa_n\varphi_f(\kappa_n s) \nu_n (1-s^2)^{1/2} \|\taub_n\|w^{1/2}
\big)
\Big| 
\\[2mm]
& & 
\hspace{13mm} 
\times
f(\kappa_n s)
(1-s^2)^{(p-3)/2}
\,
dG_p(w) 
ds
=
o\bigg(
\frac{1}{n(\kappa_n \varphi_f(\kappa_n))^{(p-1)/2}}
\bigg)
.
\end{eqnarray*}
Since~$f$ belongs to~$\mathcal{F}_{\rm LAN}(p,\kappa_n,\|\taub_n\|)$, the result then follows.  
\cqfd
\vspace{-1mm}



\bibliographystyle{imsart-nameyear.bst} 
\bibliography{Paper.bib}           
\vspace{3mm} 



\end{document}